\documentclass{amsart}
\usepackage{tabu}
\usepackage[margin=95pt]{geometry}
\usepackage{enumitem}
\allowdisplaybreaks

\usepackage{amssymb}
\usepackage{amsmath}
\usepackage{amscd}
\usepackage{amsbsy}
\usepackage{comment}
\usepackage{enumitem}
\usepackage[matrix,arrow]{xy}
\usepackage{hyperref}
\DeclareMathOperator{\Cl}{Cl}
\DeclareMathOperator{\Norm}{Norm}
\DeclareMathOperator{\Rad}{Rad}

\DeclareMathOperator{\GL}{GL}

\DeclareMathOperator{\Gal}{Gal}

\DeclareMathOperator{\ord}{ord}

\DeclareMathOperator{\lcm}{lcm}

\newcommand{\Q}{{\mathbb Q}}
\newcommand{\Z}{{\mathbb Z}}

\newcommand{\F}{{\mathbb F}}
\newcommand{\cA}{\mathcal{A}}
\newcommand{\cD}{\mathcal{D}}
\newcommand{\cB}{\mathcal{B}}

\newcommand{\cL}{\mathcal{L}}

\newcommand{\cX}{\mathcal{X}}
\newcommand{\cY}{\mathcal{Y}}
\newcommand{\cE}{\mathcal{E}}

\newcommand{\cN}{\mathcal{N}}

\newcommand{\cT}{\mathcal{T}}
\newcommand{\cS}{\mathcal{S}}
\newcommand{\cI}{\mathcal{I}}
\newcommand{\cZ}{\mathcal{Z}}

\newcommand{\cC}{\mathcal C}

\newcommand{\cW}{\mathcal W}
\newcommand{\cU}{\mathcal U}
\newcommand{\OO}{{\mathcal O}}
\newcommand{\gn}{\mathfrak{n}}

\newcommand{\fq}{\mathfrak{q}}
\newcommand{\fP}{\mathfrak{P}}
\newcommand{\fA}{\mathfrak{A}}
\newcommand{\fB}{\mathfrak{B}}
\newcommand{\fL}{\mathfrak{L}}
\newcommand{\fC}{\mathfrak{C}}

\def\mod#1{{\ifmmode\text{\rm\ (mod~$#1$)}
\else\discretionary{}{}{\hbox{ }}\rm(mod~$#1$)\fi}}

\begin {document}

\newtheorem{thm}{Theorem}
\newtheorem{lem}{Lemma}[section]
\newtheorem{prop}[lem]{Proposition}

\newtheorem{cor}[lem]{Corollary}

\theoremstyle{definition}

\theoremstyle{remark}

\title[Prime power gaps]
{Differences between perfect powers : prime power gaps}
%\author{Samir Siksek}
%\address{Mathematics Institute\\
%	University of Warwick\\
%	Coventry\\
%	CV4 7AL \\
%	United Kingdom}

%\email{s.siksek@warwick.ac.uk}
\author[Michael Bennett]{Michael A. Bennett}
\address{Department of Mathematics, University of British Columbia, Vancouver, B.C., V6T 1Z2 Canada}
\email{bennett@math.ubc.ca}

\author{Samir Siksek}
\address{Mathematics Institute, University of Warwick, Coventry CV4 7AL, United Kingdom}
\email{S.Siksek@warwick.ac.uk}
\thanks{The first-named author is supported by NSERC. The second-named author is
supported by an EPSRC Grant EP/S031537/1 \lq\lq Moduli of
elliptic curves and classical Diophantine problems\rq\rq.}

\date{\today}


\begin {abstract}
We develop machinery to explicitly determine, in many instances, when the difference $x^2-y^n$ is divisible only by powers of  a given fixed prime. This combines a wide variety of techniques from Diophantine approximation (bounds for linear forms in logarithms, both archimedean and non-archimedean, lattice basis reduction, methods for solving Thue-Mahler and $S$-unit equations, and the Primitive Divisor Theorem of Bilu, Hanrot and Voutier) and classical Algebraic Number Theory, with results derived from the modularity of Galois representations attached to Frey-Hellegoaurch elliptic curves. By way of example, we completely solve the equation
$$
x^2+q^\alpha = y^n, \; \; 
$$
where $2 \leq q < 100$ is prime, and $x, y, \alpha$ and $n$ are integers with $n \geq 3$ and $\gcd (x,y)=1$.
\end {abstract}
\maketitle

%------------------------------
\section{Introduction}
%------------------------------

The {\it Lebesgue-Nagell} equation
\begin{equation} \label{LebNag}
x^2+D = y^n
\end{equation}
has a very extensive literature, motivated, at least in part, by attempts to
extend Mih\u{a}ilescu's Theorem \cite{Mih} (n\'ee Catalan's Conjecture) to
larger gaps in the sequence of perfect powers. In equation (\ref{LebNag}), we
will suppose that $x$ and $y$ are coprime nonzero integers, and that the prime
divisors of $D$ belong to a fixed, finite set of primes $S$.  Under these
assumptions, bounds for linear forms in logarithms, $p$-adic and
complex, imply that the set of integer solutions $(x,y,n)$ to (\ref{LebNag}),
with $|y|>1$ and $n \geq 3$, is finite and effectively determinable. If, in
addition, we suppose that $D$ is positive and that $y$ is odd, then these
solutions may be explicitly determined, provided $|S|$ is not too large,
through appeal to the Primitive Divisor Theorem of Bilu, Hanrot and Voutier
\cite{BHV}, in conjunction with techniques from Diophantine approximation.

If either $D > 0$ and $y$ is even, or if $D<0$, the Primitive Divisor Theorem cannot be applied to solve equation (\ref{LebNag}) and we must work rather harder, appealing to either 
bounds for linear forms in logarithms or to  results based upon the modularity of Galois representations associated to certain Frey--Hellegouarch elliptic curves. In a companion paper \cite{BeSi35}, we develop machinery for handling (\ref{LebNag}) in the first difficult case where $D > 0$ and $y$ is even. 
%In the paper at hand, we will primarily focus our attention on the second case, where $D < 0$.
Though the techniques we discuss in the present paper are more widely
applicable, we will, for simplicity, restrict attention to the case where $D$
in equation (\ref{LebNag}) is divisible by a single prime $q$, whilst
treating both the cases $D < 0$ and $D > 0$.  
That is, we will concern ourselves
primarily with the equation
\begin{equation}\label{eqn:general}
x^2+(-1)^\delta q^\alpha=y^n, \qquad q \nmid x,
\end{equation}
where $\delta \in \{0,1\}$ and $\alpha$ is a nonnegative integer.
In case $\delta=0$, our main result is the following.
\begin{thm} \label{thm:everything}
If $x, y, q, \alpha$ and $n$ are positive integers with $q$ prime, $2 \leq q < 100$, $q \nmid x$, $n \geq 3$ and
\begin{equation} \label{covid-e}
x^2 +q^{\alpha} = y^n,
\end{equation}
then  $(q,\alpha,y,n)$ is one of
\begin{gather*}
(2,1,3,3), 
\; (2,2,5,3), 
\; (2,5,3,4), 
\; (3,5,7,3), 
\; (3,4,13,3), 
\; (7,1,2,3), 
\; (7,3,8,3), 
\; (7,1,32,3), \\
 (7,2,65,3), 
\; (7,1,2,4), 
\; (7,2,5,4), 
\; (7,1,2,5), 
\; (7,1,8,5), 
\; (7,1,2,7), 
\; (7,3,2,9), 
\; (7,1,2,15)\\
 (11,1,3,3),
\; (11,1,15,3), 
\;  (11,2,5,3), 
\; (11,3,443,3), 
\; (13,1,17,3), 
\; (17,1,3,4), 
\; (19,1,7,3), \\
 (19,1,55,5), 
\; (23,1,3,3), 
\; (23,3,71,3), 
\; (23,3,78,4), 
\; (23,1,2,5), 
\; (23,1,2,11), 
\; (29,2,5,7), \\
(31,1,4,4), 
\; (31,1,2,5), 
\; (31,1,2,8), 
\; (41,2,29,4), 
\; (41,2,5,5), 
\; (47,1,6,3), 
\; (47,1,12,3), \\
 (47,1,63,3), 
\; (47,2,17,3), 
\; (47,3,74,3), 
\; (47,1,3,5), 
\; (47,1,2,7), 
\; (53,1,9,3), 
\; (53,1,29,3),\\ 
(53,1,3,6), 
\; (61,1,5,3), 
\; (67,1,23,3), 
\; (71,1,8,3), 
\; (71,1,6,4), 
\; (71,1,3,7), 
\; (71,1,2,9), \\ 
(79,1,20,3),
\; (79,1,2,7), 
\; (83,1,27,3), 
\; (83,1,3,9), 
\; (89,1,5,3),
\; (97,2,12545,3),
\;  (97,1,7,4). \\
\end{gather*}
\end{thm}

One might note that the restriction $q \nmid x$ can be removed, with a modicum of effort, at least for certain values of $q$.
%, including $q \in \{ 2, 3, 11, 13, 19 \}$. 
The cases where primitive divisor arguments are inapplicable correspond to $q
\in \{ 7,23, 31, 47, 71, 79 \}$ and $y$ even (and this is where the great
majority of work lies in proving Theorem \ref{thm:everything}). If $q=7$,
Theorem \ref{thm:everything} generalizes recent work of Koutsianas \cite{Kou}, who
established a like result under certain conditions upon $\alpha$ and $q$, and,
in particular, showed that equation (\ref{covid-e}) has no solutions with $q=7$
and prime $n \equiv 13, 23 \mod{24}$. We note that the solution(s) with $q=83$
were omitted in the statement of Theorem 1 of Berczes and Pink \cite{BePi2}.

Our results for \eqref{eqn:general} with $\delta=1$
are less complete, at least when $\alpha$ is odd.
\begin{thm} \label{thm:odd}
Suppose that
\begin{equation}\label{eqn:goodq}
q \in \{ 7, 11, 13, 19, 23, 29, 31, 43, 47, 
53, 59, 61, 67, 71, 79, 83 \}.
\end{equation}
If $x$ and $n$ are positive integers, $q \nmid x$, $n \geq 3$ and
\begin{equation} \label{covid1}
x^2 - q^{2k+1} = y^n,
\end{equation}
where $y$ and $k$ are integers, then $(q,k,y,n)$ is one of
\begin{gather*}
(7,2,393,3),
\; (7, 1, -3, 5 ), 
\; (11,1,37,3)
\; ( 11, 0, 5, 5 ), 
\; ( 11, 1, 37, 3 ), 
\; ( 13, 0, 3, 5 ), \\
(19,0,5,3),
\; (19,2,-127,3),
\; (19, 0, -3, 4 ), 
\; ( 19, 0, 3, 4 ), 
\; (23,1,1177,3),
\; (31, 0, -3, 3 ),\\ 
 ( 43, 0, -3, 3 ),  
\; ( 71, 0, 5, 3 ), 
\; (71, 1, -23, 3 ) \; \mbox{ or } 
\; (79,0,45,3). 
\end{gather*}
\end{thm}
To the best of our knowledge, these are the first examples of primes $q$ for which equation (\ref{covid1}) has been completely solved (though the cases with $k=0$ are treated in the thesis of Barros \cite{Barr}).
There are eight other primes in the range
$3 \le q < 100$ for which we are unable to give
a similarly satisfactory statement.
For four of these, namely $q=3, 5, 17$ and $37$,
the equation \eqref{covid1} has a solution
with $y= \pm 1$. 
For such primes we are
unaware of any results that would enable us to completely treat fixed exponents
$n$ of moderate size; this difficulty
is well-known for the $D=-2$ case of (\ref{LebNag}). 
One should note that it is relatively easy to solve  \eqref{covid1} for $q \in \{ 3, 5, 37 \}$, under the additional assumption that $y$ is even (and somewhat harder if $q=17$ and $y$ is even).
For the other four primes, namely $q=41, 73, 89$ and $97$,
we give a method which appears theoretically capable of 
success, but is alas prohibitively expensive,
computationally-speaking. We content ourselves by
proving the following modest result for these primes.
\begin{thm}\label{thm:modest}
Let $q \in \{ 41, 73, 89, 97 \}$.
The only solutions to \eqref{covid1}
with $q \nmid x$ and $3 \le n \le 1000$ are with $(q, k, y, n)$ one of
\begin{gather*}
( 41, 0, -2, 5 ),\;
(41, 0, 2, 3 ),\;
( 41, 0, 2, 7 ),\;
(41, 1, 10, 5 ),\;
( 73, 0, -6, 4 ),\\
( 73, 0, -4, 3 ),\;  
( 73, 0, 2, 3 ),\; 
( 73, 0, 3, 3 ), \;
( 73, 0, 6, 3 ), \;
(73, 0, 6, 4 ),\; 
( 73, 0, 72, 3 ),\\ 
( 89, 0, -4, 3 ),\; 
( 89, 0, -2, 3 ),\; 
( 89, 0, 2, 5 ),\; 
( 89, 0, 2, 13 ) \; \mbox{ or } \;
(97,0,2,7).
\end{gather*}
There are no solutions to equation  \eqref{covid1} with $n> 1000$, $q \nmid x$ and either $q=73$ and $y \equiv 0 \mod{2}$, or with $q=97$ and $y \equiv 1 \mod{2}$.
\end{thm}

The additional assumption that the exponent of our prime $q$ is even simplifies matters considerably.
In the case of equation (\ref{covid-e}),  Berczes and Pink \cite{BePi} deduced Theorem \ref{thm:everything} for even values of $\alpha$ (whence primitive divisor technology works efficiently). For completeness, we extend this to $q < 1000$.
\begin{thm} \label{thm:even+}
If $x, y, q, k$ and $n$ are positive integers with $q$ prime, $2 \leq q < 1000$, $q \nmid x$, $n \geq 3$ and
\begin{equation} \label{covid-plus}
x^2 +q^{2k} = y^n,
\end{equation}
then $(q,k,y,n)$ is one of
	\begin{gather*}
 (2,1,5,3), (3,2,13,3), (7,1,65,3), (7,1,5,4),  (11,1,5,3), (29,1,5,7), (41,1 29,4), (41,1,5,5), \\
 (47,1,17,3),  (97,1,12545,3),  (107,1,37,3), (191,1,65,3), (239,1,169,4), (239,1,13,8), \\
 (431,1,145,3), (587,1,197,3) \; \mbox{ or } \;  (971,1,325,3). \\
\end{gather*}
\end{thm}

More interesting for us is the case where the difference $x^2-y^n$ is positive (so that primitive divisor arguments are inapplicable and there are no prior results available in the literature). Here, we prove the following.

\begin{thm} \label{thm:even}
If $x,  q, k$ and $n$ are positive integers with $q$ prime,  $2 \leq q < 1000$, $q \nmid x$, $n \geq 3$ and
\begin{equation} \label{covid}
x^2 - q^{2k} = y^n,
\end{equation}
where $y$ is an integer, then $(q,k,y,n)$ is one of
$$
\begin{array}{c}
(3,1,-2,3), (3,1,40,3), (3,1, \pm 2,4), (3,2,-2,5), (5,2,6,3), (7,2,15,3), 
( 7, 1, 2, 5 ),
 (11,1,12,3), \\ 
(11,2,3,5), (13,1,3,3), (13,1,12,5), (17,1,-4,3), (17,1,\pm 12, 4), 
( 17, 2, 42, 3 ),
(29,1,-6,3),  \\
 (31,1,2,7), (43,1,-12,3),  (43,1,126,3),  (43,4,96222,3),  (47,1,6300,3), (53,1,6,3), \\
  (71,1,30,3), (71,2,-136,3), (89,1,84,3), (97, 2, 3135,3), (101,1,24,3), (109,1,20,3),   \\
(109,1,35,3), (109,1,570,3), (127,1,-10,3),  (127,1,8,3),  (127,1,198,3),  (127,1,2,9), \\
 (179,1,-30,3), (193,1,63,3), (197,1,260,3), (223,1,30,3),    (251,1,-10,3),  (251,1,-6,5), \\
  (257,1,-4,5),  (263,1,2418,3),  (277,1,-30,3), (307,1,60,3),   (307,1,176,3), (307,2,2262,3),  \\
(359,1,-28,3),  (383,2,25800,3),  (397,1,-42,3), (431,1,12,3), (433,1,-12,3), \\
 (433,1,143,3),  (433,2,26462,3),  (479,1,90,3),  (499,1,-12,5), (503,1,828,3), 
(557,1,-60,3),\\  (577,1, \pm 408,4), (593,1,-70,3), (601,1,72,3), (659,1,42,3),  
( 683, 1, 193346, 3 ),
(701,1,4452,3), \\
(727,1,18,3), (739,1,234,3), (769,1,255,3), (811,1,-70,3), (857,1,-72,3) \mbox{ or } (997,1,48,3).
\end{array}
$$
\end{thm}
% added three missing soluntions 2/8/2021!
% [ 7, 1, 2, 5 ],
%    [ 17, 2, 42, 3 ],
%    [ 683, 1, 193346, 3 ]

We note that, with sufficient computational power, there is no obstruction to extending the results of Theorems \ref{thm:even+} and \ref{thm:even} to larger prime values $q$. Without fundamentally new ideas, it is not clear that the same may be said of, for example,  Theorem
\ref{thm:everything}. In this case,  the bounds we obtain upon the exponent $n$ via linear forms in logarithms, even for relatively small $q$, leave us with a  computation which, while finite, is barely tractable. 

\bigskip

Equation~\eqref{Main-eq} has been completely resolved
(\cite{Ivorra}, \cite{SiksekDio})
for $q=2$, except for the case $(\alpha,\delta)=(1,1)$
which corresponds to $D=-2$ in \eqref{LebNag}.
The solutions for $q=2$
in our theorems are included for completeness.
For the remainder
of the paper, we suppose that $q$ is an odd prime.
In particular, we are concerned with the equation
\begin{equation} \label{Main-eq}
x^2 + (-1)^\delta q^\alpha = y^n, \qquad \gcd(x,y)=1, \quad \alpha>0, 
\end{equation}
where $q$ is a fixed odd prime, $n \ge 3$,
and $\delta \in \{ 0, 1 \}$.

\bigskip

Our proofs will use a broad combination of techniques,
which include:
\begin{itemize}
\item Lower bounds for linear forms in complex and $p$-adic
logarithms which yield bounds for the exponent $n$
in \eqref{Main-eq}.
\item Frey--Hellegouarch curves
and their Galois representations which provide a wealth of 
local information regarding solutions to \eqref{Main-eq}.
\item The celebrated Primitive Divisor Theorem
of Bilu, Hanrot and Voutier, that can be
used to treat most cases of \eqref{Main-eq}
when $y$ is odd and $\delta=0$.
\item Elementary descent arguments that reduce \eqref{Main-eq}
for a fixed exponent $n$ to Thue--Mahler equations,
which are possible to resolve thanks to the Thue--Mahler
solver associated to \cite{GKMS}.
\end{itemize}

\bigskip

The outline of this paper is as follows. 
In Section \ref{sec:tiny}, we solve the 
equation $x^2+(-1)^\delta q^\alpha=y^n$ for $n \in \{ 3, 4 \}$ 
and $3 \le q <100$ by reducing the problem to the 
determination of $S$-integral points on elliptic curves.
In Section \ref{sec:elementary}, we solve the
equation $x^2-q^{2k}=y^n$, for $q$ in the
range $3 \le q < 1000$,  with $y$ odd using an elementary
sieving argument; this completes the proof of Theorem~\ref{thm:even} in case $y$ is odd.
Next, Section~\ref{sec:BHV} provides a short overview of Lucas
sequences, their ranks of apparition, and the Primitive Divisor
Theorem of Bilu, Hanrot and Voutier.
We make use of this machinery in Section~\ref{Thm4}
to solve the equation $x^2+q^{2k}=y^n$
for $q$ in the range $3 \le q<1000$, thereby proving
Theorem~\ref{thm:even+}. Section~\ref{sec:reduction}
reduces the equation $x^2-q^{2k}=y^n$,
for even values of  $y$, to Thue--Mahler equations
of the form 
\begin{equation}\label{eqn:binary}
y_1^n-2^{n-2} y_2^n=q^k. 
\end{equation}
In Section~\ref{sec:modular}, we give a brief outline
of the modular approach to Diophantine equations.
Section~\ref{sec:via} applies this modular approach,
particularly the $(n,n,n)$ Frey--Hellegouarch elliptic
curves of Kraus \cite{Krausppp}
to equation (\ref{eqn:binary}); this allows us
to deduce that there are no solutions
for $3 \le q <1000$ except for possibly $q \in \{ 31, 127, 257 \}$,
where the mod $n$ representation of the Frey--Hellegouarch
curve arises from that of an elliptic curve
with full $2$-torsion and conductor $2q$.
Before we can complete the proof of Theorem~\ref{thm:even},
we need an upper bound for the exponent $n$.
We give a sharpening of Bugeaud's bound \cite{Bu-Acta}
for the equation $x^2-q^{2k}=y^n$, which 
uses \eqref{eqn:binary}
and the theory of linear forms in real
and $p$-adic logarithms.
In Section~\ref{sec:sieve}, we complete the proof
of Theorem~\ref{thm:even};
our approach makes use of a sieving
technique that builds on the information
obtained from the modular approach
in Section~\ref{sec:via} and the upper bound
for $n$ established in Section~\ref{sec:upper}.
The remainder of the paper is concerned
with \eqref{Main-eq} where $\alpha=2k+1$,
and for $3 \le q < 100$.
In Section~\ref{sec:oddodd+} , we resolve
$x^2+q^{2k+1}=y^n$ with $y$ odd 
with the help of the Primitive Divisor Theorem,
and in Section~\ref{sec:n=5} we solve
$x^2-q^{2k+1}=y^5$ by reducing to Thue--Mahler
equations.

It remains, then, to handle the equations $x^2-q^{2k+1}=y^n$ and $x^2+q^{2k+1}=y^n$ where, in the latter case, we may additionally assume that $y$ is even.
In Section \ref{sec:FreyH}, we study the more general equation
\begin{equation} \label{bella}
y^n + q^\alpha z^n = x^2, \qquad \qquad \gcd(x,y)=1
\end{equation}
where $q$ is prime,
using Galois representations of Frey--Hellegouarch curves. Our approach builds on previous work
of the first author and Skinner \cite{BenS}, and also on the work of Ivorra and Kraus \cite{IvKr}. 
We then restrict ourselves in Section \ref{Sec14} to the case $z=\pm 1$  and $\alpha$ odd in (\ref{bella}). In this section, we develop a variety of sieves based upon local information coming from the Frey--Hellegouarch curves that allows us, in many situations, to eliminate values of $q$ from consideration completely and, in the more difficult cases, to solve 
equation (\ref{Main-eq}) for a fixed pair $(q,n)$. In particular, we employ this strategy to complete the proofs of  Theorems~\ref{thm:odd} and \ref{thm:modest}. Finally, in Section \ref{Thm2}, we return to bounds for linear forms in $p$-adic and complex logarithms to derive explicit upper bounds upon $n$ in (\ref{Main-eq}), and then report upon a (somewhat substantial) computation to use the arguments of Section \ref{Sec14} to solve (\ref{Main-eq}) for all remaining pairs $(q,n)$ required to finish the proof of Theorem \ref{thm:everything}.

%To prove Theorems \ref{thm:everything}, \ref{thm:odd}, \ref{thm:even+} and \ref{thm:even}, we will focus on the special case of equation (\ref{eqn:major}) with $z = \pm 1$. In particular, for the remainder of the paper, we will 
%consider the equation \eqref{Main-eq}

%%%%%%%%%%%%%%%%%%%%%%%%%%%%%%%%%%%%%%%%%%%%%%%%%%%%%%%%%%%%%%%%
%\section{A first aside : Theorems \ref{thm:everything}, \ref{thm:odd}, \ref{thm:even+} and \ref{thm:even}  for  $n \in \{ 3, 4 \}$} \label{tiny}
\section{Reduction to $S$-integral points on elliptic curves for $n \in \{3,4\}$} \label{sec:tiny}
%%%%%%%%%%%%%%%%%%%%%%%%%%%%%%%%%%%%%%%%%%%%%%%%%%%%%%%%%%%%%%%%

In the following sections, it will be of value to us to assume that the exponent $n$ in equation (\ref{Main-eq}) is not too small. This is primarily to ensure that the Frey--Hellegouarch curve we attach to a putative solution
has a corresponding 
 mod $n$ Galois representation that is irreducible.
For suitably large prime values of $n$ (typically, $n \geq 7$), the desired irreducibility follows
from Mazur's isogeny theorem. In Section \ref{sec:BHV}, such an assumption allows us to (mostly)  ignore so-called \emph{defective} Lucas sequences.

 In this section, we treat separately the cases $n=3$ and $n=4$ for $q<100$,
where the problem of solving equation (\ref{Main-eq})  reduces immediately to
one of determining $S$-integral points on specific models of genus one curves;
here $S=\{q\}$. This approach falters for many values of $q$ in the range
$100<q<1000$ as we are often unable to compute the Mordell--Weil groups of the 
relevant elliptic curves. Thus for the proofs of Theorems~\ref{thm:even+} and~\ref{thm:even}
for exponents $n=3$, $n=4$, where we treat values of $q$ less than $1000$,
we shall employ different techniques including sieving arguments and
reduction to Thue--Mahler equations.
%If $\delta=1$, we will restrict attention to those primes $q \neq t^2 \pm 1 $
%for any integers $t > 1$. In particular, for $\delta=1$ and $q < 100$, we will
%not consider
%$$
%q \in \{ 2, 3, 5, 17, 37 \}.
%$$

\subsection{The case $n=3$}\label{subsec:peq3}
Supposing that we have a coprime solution to (\ref{Main-eq}) with $n=3$, we can write $\alpha=6b+c$, where $0 \le c \le 5$.
Taking $X=y/q^{2b}$ and $Y=x/q^{3b}$, it follows that $(X,Y)$ is an $S$-integral point on the
elliptic curve
\begin{equation}\label{eqn:sint}
Y^2=X^3+(-1)^{\delta+1} q^c,
\end{equation}
where $S=\{q\}$. Here, for a particular choice of $\delta \in \{0,1\}$
and prime $q$ we may use the standard method for computing
$S$-integral points on elliptic curves based on lower
bounds for linear forms in elliptic logarithms (e.g. \cite{PZG}).
We made use of the built-in \texttt{Magma} (\cite{magma})
implementation of this method to compute these $S$-integral points
on \eqref{eqn:sint} for $\delta \in \{0,1\}$ and $2 \le q<100$.
We obtained a total of $83$ solutions to \eqref{Main-eq} for these
values of $q$ with $\alpha>0$, $x>0$, $y \ne 0$ and $\gcd(x,y)=1$.
These are given in Table~\ref{table:p=3}.
%
%Here are some interesting ones:
%\[
%21063928^2-17^7=76271^3, \quad
%30042907^2-43^8=96222^3, \quad
%1405096^2+97^2=12545^3.
%\]

\begin{table}
\caption{
Solutions to the equation
$x^2 + (-1)^\delta q^\alpha = y^3$
for primes $2 \le q < 100$, $\delta \in \{0,1\}$
and $x$, $y$, $\alpha$ integers satisfying
$\alpha>0$, $x>0$, $y \ne 0$, and $\gcd(x,y)=1$.
}
\label{table:p=3}
\centering
\begin{minipage}[t]{0.25\textwidth}%\centering
{
\tabulinesep=1.2mm
\begin{tabu}{||c|c|c|c|c||}
\hline
$q$ & $\delta$ & $\alpha$ & $x$ & $y$\\
\hline\hline
$ 2 $ & $ 0 $ & $ 1 $ & $ 5 $ & $ 3 $\\
\hline
$ 2 $ & $ 0 $ & $ 2 $ & $ 11 $ & $ 5 $\\
\hline
$ 2 $ & $ 1 $ & $ 1 $ & $ 1 $ & $ -1 $\\
\hline
$ 2 $ & $ 1 $ & $ 7 $ & $ 71 $ & $ 17 $\\
\hline
$ 2 $ & $ 1 $ & $ 9 $ & $ 13 $ & $ -7 $\\
\hline
$ 2 $ & $ 1 $ & $ 3 $ & $ 3 $ & $ 1 $\\
\hline
$ 3 $ & $ 0 $ & $ 4 $ & $ 46 $ & $ 13 $\\
\hline
$ 3 $ & $ 0 $ & $ 5 $ & $ 10 $ & $ 7 $\\
\hline
$ 3 $ & $ 1 $ & $ 1 $ & $ 2 $ & $ 1 $\\
\hline
$ 3 $ & $ 1 $ & $ 2 $ & $ 1 $ & $ -2 $\\
\hline
$ 3 $ & $ 1 $ & $ 2 $ & $ 253 $ & $ 40 $\\
\hline
$ 5 $ & $ 1 $ & $ 1 $ & $ 2 $ & $ -1 $\\
\hline
$ 5 $ & $ 1 $ & $ 4 $ & $ 29 $ & $ 6 $\\
\hline
$ 7 $ & $ 0 $ & $ 1 $ & $ 1 $ & $ 2 $\\
\hline
$ 7 $ & $ 0 $ & $ 1 $ & $ 181 $ & $ 32 $\\
\hline
$ 7 $ & $ 0 $ & $ 2 $ & $ 524 $ & $ 65 $\\
\hline
$ 7 $ & $ 0 $ & $ 3 $ & $ 13 $ & $ 8 $\\
\hline
$ 7 $ & $ 1 $ & $ 4 $ & $ 76 $ & $ 15 $\\
\hline
$ 7 $ & $ 1 $ & $ 5 $ & $ 7792 $ & $ 393 $\\
\hline
$ 11 $ & $ 0 $ & $ 1 $ & $ 4 $ & $ 3 $\\
\hline
$ 11 $ & $ 0 $ & $ 1 $ & $ 58 $ & $ 15 $\\
\hline
$ 11 $ & $ 0 $ & $ 2 $ & $ 2 $ & $ 5 $\\
\hline
$ 11 $ & $ 0 $ & $ 3 $ & $ 9324 $ & $ 443 $\\
\hline
$ 11 $ & $ 1 $ & $ 2 $ & $ 43 $ & $ 12 $\\
\hline
$ 11 $ & $ 1 $ & $ 3 $ & $ 228 $ & $ 37 $\\
\hline
$ 13 $ & $ 0 $ & $ 1 $ & $ 70 $ & $ 17 $\\
\hline
$ 13 $ & $ 1 $ & $ 2 $ & $ 14 $ & $ 3 $\\
\hline
$ 17 $ & $ 1 $ & $ 1 $ & $ 3 $ & $ -2 $\\
\hline
$ 17 $ & $ 1 $ & $ 1 $ & $ 4 $ & $ -1 $\\
\hline
\end{tabu}
}
\end{minipage}
\hfill
\begin{minipage}[t]{0.31\textwidth}%\centering
{
\tabulinesep=1.2mm
\begin{tabu}{||c|c|c|c|c||}
\hline
$q$ & $\delta$ & $\alpha$ & $x$ & $y$\\
\hline\hline
$ 17 $ & $ 1 $ & $ 1 $ & $ 5 $ & $ 2 $\\
\hline
$ 17 $ & $ 1 $ & $ 1 $ & $ 9 $ & $ 4 $\\
\hline
$ 17 $ & $ 1 $ & $ 1 $ & $ 23 $ & $ 8 $\\
\hline
$ 17 $ & $ 1 $ & $ 1 $ & $ 282 $ & $ 43 $\\
\hline
$ 17 $ & $ 1 $ & $ 1 $ & $ 375 $ & $ 52 $\\
\hline
$ 17 $ & $ 1 $ & $ 7 $ & $ 21063928 $ & $ 76271 $\\
\hline
$ 17 $ & $ 1 $ & $ 1 $ & $ 378661 $ & $ 5234 $\\
\hline
$ 17 $ & $ 1 $ & $ 2 $ & $ 15 $ & $ -4 $\\
\hline
$ 17 $ & $ 1 $ & $ 4 $ & $ 397 $ & $ 42 $\\
\hline
$ 19 $ & $ 0 $ & $ 1 $ & $ 18 $ & $ 7 $\\
\hline
$ 19 $ & $ 1 $ & $ 1 $ & $ 12 $ & $ 5 $\\
\hline
$ 19 $ & $ 1 $ & $ 5 $ & $ 654 $ & $ -127 $\\
\hline
$ 23 $ & $ 0 $ & $ 1 $ & $ 2 $ & $ 3 $\\
\hline
$ 23 $ & $ 0 $ & $ 3 $ & $ 588 $ & $ 71 $\\
\hline
$ 23 $ & $ 1 $ & $ 3 $ & $ 40380 $ & $ 1177 $\\
\hline
$ 29 $ & $ 1 $ & $ 2 $ & $ 25 $ & $ -6 $\\
\hline
$ 31 $ & $ 1 $ & $ 1 $ & $ 2 $ & $ -3 $\\
\hline
$ 37 $ & $ 1 $ & $ 1 $ & $ 6 $ & $ -1 $\\
\hline
$ 37 $ & $ 1 $ & $ 1 $ & $ 8 $ & $ 3 $\\
\hline
$ 37 $ & $ 1 $ & $ 1 $ & $ 3788 $ & $ 243 $\\
\hline
$ 37 $ & $ 1 $ & $ 3 $ & $ 228 $ & $ 11 $\\
\hline
$ 41 $ & $ 1 $ & $ 1 $ & $ 7 $ & $ 2 $\\
\hline
$ 43 $ & $ 1 $ & $ 1 $ & $ 4 $ & $ -3 $\\
\hline
$ 43 $ & $ 1 $ & $ 2 $ & $ 11 $ & $ -12 $\\
\hline
$ 43 $ & $ 1 $ & $ 8 $ & $ 30042907 $ & $ 96222 $\\
\hline
$ 43 $ & $ 1 $ & $ 2 $ & $ 1415 $ & $ 126 $\\
\hline
$ 47 $ & $ 0 $ & $ 1 $ & $ 13 $ & $ 6 $\\
\hline
$ 47 $ & $ 0 $ & $ 1 $ & $ 41 $ & $ 12 $\\
\hline
$ 47 $ & $ 0 $ & $ 1 $ & $ 500 $ & $ 63 $\\
\hline
$ 47 $ & $ 0 $ & $ 2 $ & $ 52 $ & $ 17 $\\
\hline
\end{tabu}
}
\end{minipage}
\hfill
\begin{minipage}[t]{0.31\textwidth}%\centering
{
\tabulinesep=1.2mm
\begin{tabu}{||c|c|c|c|c||}
\hline
$q$ & $\delta$ & $\alpha$ & $x$ & $y$\\
\hline\hline
$ 47 $ & $ 0 $ & $ 3 $ & $ 549 $ & $ 74 $\\
\hline
$ 47 $ & $ 1 $ & $ 2 $ & $ 500047 $ & $ 6300 $\\
\hline
$ 53 $ & $ 0 $ & $ 1 $ & $ 26 $ & $ 9 $\\
\hline
$ 53 $ & $ 0 $ & $ 1 $ & $ 156 $ & $ 29 $\\
\hline
$ 53 $ & $ 1 $ & $ 2 $ & $ 55 $ & $ 6 $\\
\hline
$ 61 $ & $ 0 $ & $ 1 $ & $ 8 $ & $ 5 $\\
\hline
$ 67 $ & $ 0 $ & $ 1 $ & $ 110 $ & $ 23 $\\
\hline
$ 71 $ & $ 0 $ & $ 1 $ & $ 21 $ & $ 8 $\\
\hline
$ 71 $ & $ 1 $ & $ 1 $ & $ 14 $ & $ 5 $\\
\hline
$ 71 $ & $ 1 $ & $ 2 $ & $ 179 $ & $ 30 $\\
\hline
$ 71 $ & $ 1 $ & $ 3 $ & $ 588 $ & $ -23 $\\
\hline
$ 71 $ & $ 1 $ & $ 4 $ & $ 4785 $ & $ -136 $\\
\hline
$ 73 $ & $ 1 $ & $ 1 $ & $ 3 $ & $ -4 $\\
\hline
$ 73 $ & $ 1 $ & $ 1 $ & $ 9 $ & $ 2 $\\
\hline
$ 73 $ & $ 1 $ & $ 1 $ & $ 10 $ & $ 3 $\\
\hline
$ 73 $ & $ 1 $ & $ 1 $ & $ 17 $ & $ 6 $\\
\hline
$ 73 $ & $ 1 $ & $ 1 $ & $ 611 $ & $ 72 $\\
\hline
$ 73 $ & $ 1 $ & $ 1 $ & $ 6717 $ & $ 356 $\\
\hline
$ 79 $ & $ 0 $ & $ 1 $ & $ 89 $ & $ 20 $\\
\hline
$ 79 $ & $ 1 $ & $ 1 $ & $ 302 $ & $ 45 $\\
\hline
$ 83 $ & $ 0 $ & $ 1 $ & $ 140 $ & $ 27 $\\
\hline
$ 89 $ & $ 0 $ & $ 1 $ & $ 6 $ & $ 5 $\\
\hline
$ 89 $ & $ 1 $ & $ 1 $ & $ 5 $ & $ -4 $\\
\hline
$ 89 $ & $ 1 $ & $ 1 $ & $ 9 $ & $ -2 $\\
\hline
$ 89 $ & $ 1 $ & $ 1 $ & $ 33 $ & $ 10 $\\
\hline
$ 89 $ & $ 1 $ & $ 1 $ & $ 408 $ & $ 55 $\\
\hline
$ 89 $ & $ 1 $ & $ 2 $ & $ 775 $ & $ 84 $\\
\hline
$ 97 $ & $ 0 $ & $ 2 $ & $ 1405096 $ & $ 12545 $\\
\hline
$ 97 $ & $ 1 $ & $ 1 $ & $ 77 $ & $ 18 $\\
\hline
$ 97 $ & $ 1 $ & $ 4 $ & $ 175784 $ & $ 3135 $\\
%& & & &  \\
\hline
\end{tabu}
}
\end{minipage}
\end{table}

%For larger values of $q$, we run into some problems with this approach (due to the difficulty of finding generators of the corresponding Mordell-Weil groups). We can instead, consider the curve
%$$
%E \; \; : \; \; Y^2=X^3 -3yX+2x,
%$$
%with discriminant
%$$
%\Delta_E  = 1728 (y^3-x^2) = 1728 (-1)^\delta q^\alpha.
%$$
%A short calculation with Tate's algorithm shows that the conductor of $E$ satisfies
%$$
%N_{E} = 2^{\alpha_2} 3^{\alpha_3} q,  \mbox{ where } \alpha_2 \in \{ 5, 6 \} \mbox{ and } \alpha_3 \in \{ 2, 3 \}.
%$$

%=====================================================
\subsection{The case $n=4$}
%=====================================================
Next we consider the case $n=4$ separately.
Write $\alpha=4b+c$ where $0 \le c \le 3$.
Let $X=(y/q^{b})^2$, $Y=xy/q^{3b}$. Then $(X,Y)$ is an $S$-integral point on the
elliptic curve
\begin{equation}\label{eqn:sint2}
Y^2=X(X^2+(-1)^{\delta+1} q^c),
\end{equation}
where $S=\{q\}$.
%Here, for a particular choice of $\delta \in \{0,1\}$
%and prime $q$ we may use the standard method for computing
%$S$-integral points on elliptic curves based on lower
%bounds for linear forms in elliptic logarithms (e.g. \cite{PZG}).
We again appealed to the built-in \texttt{Magma} (\cite{magma})
implementation of this method to compute these $S$-integral points
on \eqref{eqn:sint2} for $\delta \in \{0,1\}$ and $2 \le q<100$.
We obtained a total  of $16$ solutions to \eqref{Main-eq} for these
values of $q$ with $\alpha>0$, $x>0$, $y >0$ and $\gcd(x,y)=1$.
These are given in Table~\ref{table:p=4}.

\begin{table}
\caption{
Solutions to the equation
$x^2 + (-1)^\delta q^\alpha = y^4$
for primes $2 \le q < 100$, $\delta \in \{0,1\}$
and $x$, $y$, $\alpha$ integers satisfying
$\alpha>0$, $x>0$, $y > 0$, and $\gcd(x,y)=1$.
}
\label{table:p=4}
\centering
\begin{minipage}[t]{0.31\textwidth}%\centering
{
\tabulinesep=1.2mm
\begin{tabu}{||c|c|c|c|c||}
\hline
$q$ & $\delta$ & $\alpha$ & $x$ & $y$\\
\hline\hline
$ 2 $ & $ 0 $ & $ 5 $ & $ 7 $ & $ 3 $\\
\hline
$ 2 $ & $ 1 $ & $ 3 $ & $ 3 $ & $ 1 $\\
\hline
$ 3 $ & $ 1 $ & $ 1 $ & $ 2 $ & $ 1 $\\
\hline
$ 3 $ & $ 1 $ & $ 5 $ & $ 122 $ & $ 11 $\\
\hline
$ 3 $ & $ 1 $ & $ 2 $ & $ 5 $ & $ 2 $\\
\hline
\end{tabu}
}
\end{minipage}
\hfill
\begin{minipage}[t]{0.31\textwidth}%\centering
{
\tabulinesep=1.2mm
\begin{tabu}{||c|c|c|c|c||}
\hline
$q$ & $\delta$ & $\alpha$ & $x$ & $y$\\
\hline\hline
$ 7 $ & $ 0 $ & $ 1 $ & $ 3 $ & $ 2 $\\
\hline
$ 7 $ & $ 0 $ & $ 2 $ & $ 24 $ & $ 5 $\\
\hline
$ 17 $ & $ 0 $ & $ 1 $ & $ 8 $ & $ 3 $\\
\hline
$ 17 $ & $ 1 $ & $ 2 $ & $ 145 $ & $ 12 $\\
\hline
$ 19 $ & $ 1 $ & $ 1 $ & $ 10 $ & $ 3 $\\
\hline
\end{tabu}
}
\end{minipage}
\hfill
\begin{minipage}[t]{0.31\textwidth}%\centering
{
\tabulinesep=1.2mm
\begin{tabu}{||c|c|c|c|c||}
\hline
$q$ & $\delta$ & $\alpha$ & $x$ & $y$\\
\hline\hline
$ 23 $ & $ 0 $ & $ 3 $ & $ 6083 $ & $ 78 $\\
\hline
$ 31 $ & $ 0 $ & $ 1 $ & $ 15 $ & $ 4 $\\
\hline
$ 41 $ & $ 0 $ & $ 2 $ & $ 840 $ & $ 29 $\\
\hline
$ 71 $ & $ 0 $ & $ 1 $ & $ 35 $ & $ 6 $\\
\hline
$ 73 $ & $ 1 $ & $ 1 $ & $ 37 $ & $ 6 $\\
\hline
$ 97 $ & $ 0 $ & $ 1 $ & $ 48 $ & $ 7 $\\
\hline
\end{tabu}
}
\end{minipage}
\end{table}

\section{An elementary approach to $x^2-q^{2k}=y^n$ with $y$ odd}\label{sec:elementary}
In this section, we apply an elementary factorisation argument to prove Theorem~\ref{thm:even}
for $y$ odd. In other words, we consider the following equation
\begin{equation}\label{eqn:evenodd}
x^2-q^{2k} = y^n, \qquad \text{$x$, $k$, $n$ positive integers}, \quad n \ge 3, \quad \gcd(x,y)=1, \quad \text{$y$ an odd integer}.
\end{equation}
Here $q \ge 3$ is a prime.
From this, we immediately see that
\begin{equation}\label{eqn:beginning}
x+q^k=y_1^n \; \; \mbox{ and } \; \; x-q^k=y_2^n,
\end{equation}
with $y=y_1y_2$,
so that we have
\begin{equation}\label{easy2}
y_1^n-y_2^n=2q^k.
\end{equation}
If $2 \mid n$, then $y_1^n \equiv y_2^n \equiv 1 \mod{4}$, a contradiction.
We may suppose henceforth, without loss of generality, that $n$ is an odd prime.
%We prove an elementary result that greatly simplifies matters if $y$ is odd, i.e. in case (\ref{eqn:beginning}).
%Observe that, from (\ref{easy2}),
Observe that
\begin{equation}\label{eqn:factored}
(y_1-y_2)(y_1^{n-1}+y_1^{n-2} y_2+\cdots+y_2^{n-1})=
y_1^n-y_2^n=2q^k.
\end{equation}
Clearly $y_1>y_2$ and, as they are both odd, $y_1 -y_2 \ge 2$
and $2 \mid (y_1-y_2)$.
Write
\[
d=\gcd(y_1-y_2, y_1^{n-1}+y_1^{n-2} y_2+\cdots+y_2^{n-1})
\]
so that $y_2 \equiv y_1 \mod{d}$ and
\[
0 \equiv y_1^{n-1}+y_1^{n-2} y_2+\cdots+y_2^{n-1} \equiv n y_1^{n-1} \mod{d}.
\]
Similarly, we have $n y_2^{n-1} \equiv 0 \mod{d}$ and so $d \in \{ 1, n \}$.

We first deal with the case $d=n$, whereby, from (\ref{eqn:factored}),   $q=n$.
Let $r=\ord_n(y_1-y_2) \ge 1$ and write $y_1=y_2+n^r \kappa$ where
$n \nmid \kappa$. Then
\[
\ord_n(y_1^{n-1}+y_1^{n-2} y_2+\cdots+y_2^{n-1})
=\ord_n \left(\frac{(y_2+n^r \kappa)^n \, -\, y_2^n}{n^r \kappa} \right)=1.
\]
Hence
\begin{equation}\label{eqn:pair}
y_1-y_2=2 n^{k-1} \; \mbox{ and } \;  y_1^{n-1}+y_1^{n-2} y_2+\cdots+y_2^{n-1}=n,
\end{equation}
and so
\[
n \; = \; \prod_{i=1}^{n-1} \lvert y_1-\zeta_n^i y_2 \rvert \; \ge \;
\lvert \,  \lvert y_1 \rvert \, -\, \lvert y_2 \rvert \, \rvert^{n-1}.
\]
Recall that $y_1$ and $y_2$ are both odd. If $y_2 \ne \pm y_1$, then
the right hand-side of this last inequality is at least $2^{n-1}$, which is impossible.
Thus $y_2=\pm y_1$, so that, from (\ref{eqn:pair}),  $y_1^{n-1} \mid n$.
It follows that $|y_1|=|y_2|= 1$, contradicting
\eqref{eqn:beginning}.

Thus $d=1$, whence
\begin{equation}\label{eqn:sys}
y_1-y_2=2 \; \; \mbox{ and } \; \; y_1^{n-1}+y_1^{n-2} y_2+\cdots+y_2^{n-1}=q^k.
\end{equation}
Since the polynomial $X^{n-1}+X^{n-2}+\cdots+1$ has a root
modulo $q$, the Dedekind--Kummer theorem tells us that
$q$ splits in $\Z[\zeta_n]$ and so $q \equiv 1 \mod{n}$. We therefore have the following.
\begin{prop} \label{prop:even-odd}
If $x, y, q, k$ and $n$ are positive integers satisfying \eqref{eqn:evenodd} with $n$ and $q$ prime, then $n \mid q-1$ and there exists an odd positive integer $X$ such that $y=X(X+2)$ and
\begin{equation}\label{eqn:X+2}
(X+2)^n-X^n=2 q^k.
\end{equation}
\end{prop}

This last result makes it an extremely straightforward matter to solve equation
(\ref{covid}) in case $y$ is odd. 

%Of course, in the case where $y$ is odd in equation \eqref{covid}, Proposition \ref{even-odd} provides us with much more precise information than Theorem \ref{thm:upper}. In this case, restricting attention to relatively small primes, we have the following.

\begin{lem}\label{lem:sieve}
The only solutions to \eqref{eqn:evenodd} with $3 \le q < 1000$
prime
correspond to the identities
\begin{gather*}
76^2-7^4=15^3, \quad
122^2- 11^4=3^5. \quad
14^2-13^2=3^3, \quad
175784^2- 97^4=3135^3, \\
234^2-109^2=35^3, \quad
536^2-193^2=63^3, \quad
1764^2-433^2=143^3, \quad
4144^2- 769^2=255^3.
\end{gather*}
\end{lem}
\begin{proof}
Suppose first that $n=3$, where 
equation \eqref{eqn:X+2} becomes
\begin{equation}\label{eqn:1or2}
3(X+1)^2+1=q^k.
\end{equation} 
From \cite{Cohn97} and \cite{Cohn2003}, we know that the equation $3u^2+1=y^m$ has no solutions with $m \ge 3$. We conclude that $k=1$ or $2$.
Solving equation \eqref{eqn:1or2} with $k=1$ or $2$ and $3 \le q <1000$ leads to the seven solutions with $n=3$.

We now suppose that $n \ge 5$ is prime.
By a theorem of the first author and Skinner \cite[Theorem 2]{BenS},
the only solutions to the
equation $X^n+Y^n=2Z^2$ with $n \ge 5$ prime and $\gcd(X,Y)=1$
are with either $|XY|=1$ or $(n,X,Y,Z)=(5,3,-1,\pm 11)$.
We note that if $k$ is even
then %\eqref{eqn:factored}
\eqref{eqn:X+2}
can be rewritten as $(X+2)^n-X^n=2(q^{k/2})^2$, and therefore
$n=5$, $X=1$ and $q^{k/2}=11$. This yields
the solution %\eqref{eqn:sol}.
$122^2- 11^4=3^5$.

We may therefore suppose that $k$
is odd.
%The pair of equations \eqref{eqn:pair} give
%\begin{equation}\label{eqn:presieve}
%(y_2+2)^n-y_2^n=2 q^k.
%\end{equation}
Recalling that $n \mid (q-1)$ leaves us with precisely
$191$ pairs $(q,n)$ to consider, ranging from $(11,5)$ to $(997,83)$.
%\[
%\begin{gathered}
%(11, 5),\; (23, 11),\; (29, 7),\; (31, 5),\; (41, 5),\; (43, 7), \;
% (47, 23),\; (53, 13), \;
%(59, 29),\; (61, 5),\; (67, 11),\\
% (71, 5),\;
%(71, 7),\; (79, 13),\; (83, 41),\; (89, 11), \; (101,5), \; (101,11), \; 
%\end{gathered}
%\]
Fix one of these pairs $(q,n)$ and let $\ell \nmid nq$ be an odd prime.
Let $\cZ_\ell$ be the set of $\beta \in \Z/(\ell-1)\Z$ such that
$\beta$ is odd and the polynomial
\[
(X+2)^n-X^n-2 q^\beta
\]
has a root in $\F_{\ell}$. We note that the value of
$q^k$ modulo $\ell$ depends only on the residue class of $k$
modulo $\ell-1$. From \eqref{eqn:X+2}, we deduce
that $(k \bmod{\ell}) \in \cZ_\ell$. Now let $\ell_1,\ell_2,\dotsc,\ell_m$
be a collection of odd primes with $\ell_i \nmid nq$ for $1 \leq i \leq m$. Let
\begin{equation}\label{eqn:M}
M=\lcm(\ell_1-1,\ell_2-1,\ldots,\ell_m-1)
\end{equation}
and set
\begin{equation}\label{eqn:cZ}
\cZ_{\ell_1,\dotsc,\ell_m}=\{\beta \in \Z/M\Z \; : \;
\text{$(\beta \bmod{\ell_i}) \in \cZ_{\ell_i}$ for $i=1,\dotsc,m$}\}.
\end{equation}
It is clear that $(k \bmod{M}) \in \cZ_{\ell_1,\dotsc,\ell_m}$.
We wrote a short \texttt{Magma} script which,
for each pair $(q,n)$, computed $\cZ_{\ell_1,\dotsc,\ell_m}$
where
$\ell_1,\ell_2,\dotsc,\ell_m$ are the odd primes $\le 101$
distinct from $n$ and $q$. %In $190$ of the $191$ cases we found that
In all $191$ cases we found that
$\cZ_{\ell_1,\dotsc,\ell_m}=\emptyset$,  completing the desired contradiction. 
%The only case where this fails to be true is when 
%$(q,n)=(601,5)$, where we obtain a contradiction upon choosing  the $\ell_1,\ell_2,\dotsc,\ell_m$ to be the odd primes $\le 101$,
%coprime to  $5$.
\end{proof}

\section{Lucas Sequences and the Primitive Divisor Theorem}\label{sec:BHV}
The Primitive Divisor Theorem of Bilu, Hanrot and Voutier \cite{BHV}
shall be our main tool for treating equation \eqref{Main-eq}
when $\delta=0$ and $y$ is odd. In this section, we state this result and a related theorem of Carmichael that shall
be useful later.
%In this section we introduce
%Lucas sequences,
%following \cite{BHV}.
A pair of algebraic integers $(\gamma,\delta)$ is called a
\textbf{Lucas pair} 
if $\gamma+\delta$ and $\gamma\delta$ are
non-zero coprime rational integers, and $\gamma/\delta$
is not a root of unity. 
We say that two Lucas pairs $(\gamma_1,\delta_1)$
and $(\gamma_2,\delta_2)$ are \textbf{equivalent}
if $\gamma_1/\gamma_2=\pm 1$ and $\delta_1/\delta_2=\pm 1$.
%In particular, associated
%to the Lucas pair $(\alpha,\beta)$ is a
%\textbf{characteristic polynomial}
%\[
%X^2-(\alpha+\beta) X + \alpha \beta \; \in \; \Z[X].
%\]
%This polynomial has discriminant $D=(\alpha-\beta)^2 \in \Z \setminus \{0\}$.
Given a Lucas pair $(\gamma,\delta)$ we define
the corresponding \textbf{Lucas sequence} by
\[
L_m=\frac{\gamma^m-\delta^m}{\gamma-\delta}, \qquad m=0,1,2,\dotsc.
\]
A prime $\ell$ is said to be a
\textbf{primitive divisor} of the $m$-th term 
if $\ell$ divides $L_m$ but $\ell$ does not divide $(\gamma-\delta)^2 \cdot u_1 u_2 \dotsc u_{m-1}$.
\begin{thm}[Bilu, Hanrot and Voutier \cite{BHV}]\label{thm:BHV}
Let $(\gamma,\delta)$ be a Lucas pair and write $\{L_m\}$
for the corresponding Lucas sequence.
\begin{enumerate}
\item[(i)]
If $m \ge 30$, then $L_m$ has a primitive divisor.
\item[(ii)] If $m \ge 11$ is prime, then $L_m$ has a primitive
divisor. 
\item[(iii)] $L_7$ has
	a primitive divisor unless $(\alpha,\beta)$
		is equivalent to 
		$\left((a-\sqrt{b})/2 \, ,\, (a+\sqrt{b})/2\right)$
where
\begin{equation}\label{eqn:7list}
	(a,b) \; \in \; \{ (1,-7),  (1,-19) \}.
\end{equation}
\item[(iv)] $L_5$ has a primitive divisor unless
	$(\alpha,\beta)$ is equivalent to 
		$\left( (a-\sqrt{b})/2 \, ,\, 
		(a+\sqrt{b})/2 \right)$
where
\begin{equation}\label{eqn:5list}
        (a,b) \; \in \;  \{(1,5), (1,-7), (2,-40), (1,-11),  (1,-15), 
		(12,-76), (12,-1364) \}.
\end{equation}
\end{enumerate}
\end{thm}
Let $\ell$ be a prime. 
We define the \textbf{rank of apparition of $\ell$}
in the Lucas sequence $\{L_m\}$ to be the smallest
positive integer $m$ such that $\ell \mid L_m$. We denote the
rank of apparition of $\ell$ by $m_\ell$.
The following theorem of Carmichael \cite{Car} will be useful for us; a concise proof
may be found in \cite[Theorem 8]{BGPS}
\begin{thm}[Carmichael  \cite{Car}]\label{thm:Carmichael}
Let $(\gamma,\delta)$ be a Lucas pair, and $\{L_m\}$
the corresponding Lucas sequence. Let $\ell$ be a prime.
\begin{enumerate}[label=(\roman*)]
\item If $\ell \mid \gamma \delta$ then $\ell \nmid L_m$
for all positive integers $m$.
\item Suppose $\ell \nmid \gamma \delta$.
Write $D=(\gamma-\delta)^2 \in \Z$.
\begin{enumerate}[label=(\alph*)]
\item If $\ell \ne 2$ and $\ell \mid D$, then
$m_\ell=\ell$.
\item If $\ell \ne 2$ and $\left(\frac{D}{\ell}\right)=1$, then $m_\ell \mid (\ell-1)$.
\item If $\ell \ne 2$ and $\left(\frac{D}{\ell} \right)=-1$, then $m_\ell \mid (\ell+1)$.
\item If $\ell=2$, then $m_\ell=2$ or $3$.
\end{enumerate}
\item If $\ell \nmid \gamma \delta$ then
\[
\ell \mid L_m \iff m_\ell \mid m.
\]
\end{enumerate}
\end{thm}

%------------------------------------------------------------------------------------
\section{The equation $x^2+q^{2k}=y^n$ : the proof of Theorem \ref{thm:even+}} \label{Thm4}
%------------------------------------------------------------------------------------

In this section, we prove Theorem~\ref{thm:even+} with the help of the Primitive Divisor Theorem.
We are concerned with the equation
\begin{equation}\label{eqn:evenodd2}
x^2+q^{2k} = y^n, \qquad \text{$x$, $k$, $n$ positive integers}, \quad n \ge 3, \quad \gcd(x,y)=1. 
\end{equation}
Here $q \ge 3$ is a prime. Considering this equation modulo $8$ immediately tells us that $y$ is odd and $x$ is even.

\begin{lem}
Solutions to \eqref{eqn:evenodd2} with $n=4$ and odd prime $q$
satisfy $k=1$, $q^2=2y^2-1$ and $x=(q^2-1)/2$. In particular, the only solutions to
\eqref{eqn:evenodd2} with $n=4$ and prime $3 \le q <1000$
correspond to the identities
\[
		24^2+ 7^2= 5^4, \; \; 
		840^2+ 41^2= 29^4 \; \; \mbox{ and } \; \; 
		28560^2+ 239^2= 169^4.
\]
\end{lem}
\begin{proof}
Suppose $n=4$. Then $(y^2+x)(y^2-x)=q^{2k}$,
and so
\[
	y^2+x=q^{2k}  \; \; \mbox{ and } \; \;  y^2-x=1.
\]
Thus $2y^2=q^{2k}+1$. 
By Theorem 1 of \cite{BenS},
the only solutions to the equation
$A^n+B^n=2C^2$ with $n \ge 4$, $ABC \ne 0$ and $\gcd(A,B)=1$
are with $|AB|=1$ or $(n,A,B,C)=(5,3,-1,\pm 11)$. It follows that
the equation $2y^2=q^{2k}+1$ has no solutions
with $k \ge 2$. Therefore $k=1$, and hence
$q^{2}=2y^2-1$. 
The only primes in the range $3 \le q < 1000$
such that $q^2=2y^2-1$ are $q=7$, $41$ and $239$,
which lead to the solutions in the statement of the lemma.
\end{proof}
If we are interested in solutions with $4 \mid n$,
then we note also the solution
$28560^2+ 239^2= 13^8$.

Henceforth, we will suppose that $n$ is an odd prime.
Thus
$x+q^k i = \alpha^n$, where we can write $\alpha=a+bi$,
for $a$ and $b$ coprime integers with
$y=a^2+b^2$.
Subtracting this equation from its conjugate yields
\begin{equation}\label{eqn:preLucas}
q^k \; = \; b \cdot \frac{\alpha^n-\overline{\alpha}^n}{\alpha-\overline{\alpha}}.
\end{equation}
\begin{lem}\label{lem:n=3even+}
Solutions to \eqref{eqn:evenodd2} with $n=3$ and odd prime $q$
must satisfy
\begin{enumerate}
\item[(i)] either $q=3$ and $(k,x,y)=(2,46,13)$;
\item[(ii)] or $q=3a^2-1$ for some positive integer $a$ and $(k,x,y)=(1,a^3-3a,a^2+1)$;
\item[(iii)] or $q^2=3a^2+1$ for some positive integer $a$ and $(k,x,y)=(1,8a^3+3a,4a^2+1)$.
\end{enumerate}
In particular, the only solutions to \eqref{eqn:evenodd2} with $n=3$
and prime $3 \leq q <1000$ correspond to the identities
\begin{gather*}
46^2+3^4=13^3, \quad
524^2+ 7^2= 65^3, \quad
2^2+ 11^2= 5^3, \quad
52^2+ 47^2= 17^3, \\
1405096^2+ 97^2= 12545^3, \quad
198^2+ 107^2= 37^3, \quad
488^2+ 191^2= 65^3, \quad
1692^2+ 431^2= 145^3, \\
2702^2+ 587^2= 197^3 \; \; \mbox{ and } \; \; 
5778^2+ 971^2= 325^3.
\end{gather*}
\end{lem}
\begin{proof}
Let $n=3$.
Thanks to Table~\ref{table:p=3}, we know that the only solution
with $q=3$ is the one given in (i). We may thus suppose that
$q \ge 5$.
Equation~\eqref{eqn:preLucas} gives
\[
q^k=b(3a^2-b^2).
\]
By the coprimality of $a$ and $b$, we have $b=\pm 1$ or $b=\pm q^k$. We note that
$b=-1$ gives $q^k=1-3a^2$ which is impossible. Also if $b=q^k$ then
$3a^2-q^{2k}=1$ which is impossible modulo $3$. Thus either $b=1$ or $b=-q^k$.
If $b=1$, then
$$
q^k=3a^2-1,
$$
and if $b=-q^k$ then
$$
q^{2k}=3 a^2+1.
$$
From Theorem 1.1 of \cite{BenS}, these equations have no solutions in positive integers if $k \geq 4$ or $k \geq 2$, respectively. 
If $k=3$, the elliptic curve corresponding to the first equation has Mordell-Weil rank $0$ over $\mathbb{Q}$ and it is straightforward to show that the equation has no integer solutions. We therefore have that  $k=1$ in either case. 
%i.e. that
Thus $q=3a^2-1$ or $q^2=3a^2+1$, and these yield the parametric solutions in (ii) and (iii).
%
%\begin{equation} \label{k1}
%q=3a^2-1
%\end{equation}
%or
%\begin{equation} \label{k2}
%q^2=3a^2+1.
%\end{equation}
%In each case, we find a corresponding solution to equation (\ref{covid-plus})
%from the identities
%$$
%(a^3-3a)^2 + (3a^2-1)^2= (a^2+1)^3
%$$
%and
%$$
%(8a^3+3a)^2+(3a^2+1)=(4a^2+1)^3,
%$$
%respectively.
For $5 \le  q < 1000$, the primes $q$ of the  form $3a^2-1$ are
$$
11, 47, 107, 191, 431, 587 \mbox{ and } 971.
$$
For $5 \le q < 1000$, the primes $q$ satisfying $q^2=3a^2+1$ are
$q=7$ and $97$.
These yield the solutions given in the statement of the lemma.
%$$
%7, 97 \mbox{ and } 708158977.
%$$
\end{proof}
We expect that there are infinitely many primes $q$ of the form $3a^2-1$, but are very unsure about the number of primes $q$  satisfying $q^2=3a^2+1$ (the only ones known are $7, 97$ and $708158977$). Quantifying such results, in any case, is well beyond current technology.

In view of Lemma~\ref{lem:n=3even+}, we henceforth suppose that $n$ is $\ge 5$ and prime. The following lemma now completes the proof of Theorem~\ref{thm:even+}.

\begin{lem}
Let $(k,x,y,n)$
be a solution to \eqref{eqn:evenodd2} with prime $n \ge 5$ and odd prime $q$. Then $k$ is odd, 
\begin{equation}\label{eqn:nq}
	\begin{cases}
		n \mid (q-1) & \text{if $q \equiv 1 \mod{4}$}\\
		n \mid (q+1) & \text{if $q \equiv 3 \mod{4}$},
	\end{cases}	
\end{equation}
and there is an integer $a$ such that
\[
	y=a^2+1, \qquad x=\frac{(a+i)^n+(a-i)^2}{2},
	\qquad \frac{(a+i)^n-(a-i)^n}{2i} = \pm q^k.
\]
In particular, the only solutions to \eqref{eqn:evenodd2}
with prime $3 \le q <1000$ and prime $n \ge 5$
correspond to the identities
\[
38^2+41^2 = 5^5, \qquad
278^2+ 29^2 = 5^7.
\]
\end{lem}
\begin{proof}
Suppose $n$ is $\ge 5$ and prime in \eqref{eqn:evenodd2}.
By Theorem 1 of \cite{BeElNg}, the equation $A^4+B^2=C^m$ has no solutions satisfying $\gcd(A,B)=1$, $AB \ne 0$ and $m \ge 4$.
We conclude that $k$ is odd.
We note that $(\alpha,\overline{\alpha})$ is a Lucas pair and write
$\{L_m\}$ for the corresponding Lucas sequence. By Theorem~\ref{thm:BHV},
$L_n$ must have a primitive divisor, and from \eqref{eqn:preLucas} this
primitive divisor is $q$. In particular, $q$ does not divide $D=(\alpha-\overline{\alpha})^2=-4b^2$. Thus by \eqref{eqn:preLucas} we have $b= \pm 1$ and $D=-4$.
Moreover, the rank of apparition of $q$ in the sequence is $n$.
By Theorem~\ref{thm:Carmichael}, we have $n \mid (q-1)$ if $q \equiv 1 \mod{4}$ and
$n \mid (q+1)$ if $q \equiv 3 \mod{4}$. 

We now let $q$ be a prime in the range $3 \le q < 1000$.
There are $168$ pairs $(q,n)$ with $q$ in this range
and $n$ a prime $\ge 5$ satisfying \eqref{eqn:nq},
ranging from $(19,5)$ to $(997,83)$.
For each of these pairs $(q,n)$,
and each sign $\eta=\pm 1$, we need to consider the
equation
\begin{equation}\label{eqn:preTMa}
		\frac{(a+i)^n-(a-i)^n}{2i}\;  =\;  \eta \cdot q^k
\end{equation}
where $k$ is an odd integer. 
We shall follow the sieving
approach of Lemma~\ref{lem:sieve}
to eliminate all but two of the possible $2 \times 168=336$
triples $(q,n,\eta)$.
Fix such a triple $(q,n,\eta)$. 
Let $f_{n} \in \Z[X]$ be the polynomial
\[
	f_n(X) \; =\; \frac{(X+i)^n-(X-i)^n}{2i}.
\]
Let $\ell \nmid  nq$ be an odd prime, and let $\cZ_{\ell}$ be the set
$\beta \in \Z/(\ell-1)\Z$ such that $\beta$
is odd and $f_n(X)-\eta \cdot q^\beta$ has a root in $\F_\ell$.
It follows that $(k \mod \ell) \in \cZ_\ell$.
Now let $\ell_1,\ell_2,\cdots,\ell_m$
be a collection of odd primes $\nmid qn$. Define $M$
and $\cZ_{\ell_1,\dotsc,\ell_m}$ by \eqref{eqn:M}
and \eqref{eqn:cZ} respectively. 
It is clear that $(k \mod M) \in \cZ_{\ell_1,\dotsc,\ell_m}$.
We wrote a short \texttt{Magma} script which,
for each triple $(q,n,\eta)$, computed $\cZ_{\ell_1,\dotsc,\ell_m}$
where $\ell_1\,\dotsc,\ell_m$ are the odd primes $<150$
distinct from $n$ and $q$. In all but two of the $336$ cases
we found that $\cZ_{\ell_1,\dotsc,\ell_m}=\emptyset$.
The two exceptions are $(q,n,\eta)=(41,5,1)$ and $(29,7,-1)$,
and so these are the only two cases we need to consider.
Let
\[
	F_n(X,Y)=\frac{(X+iY)^n-(X-iY)^n}{2i Y}.
\]
This is a homogeneous degree $n-1$ polynomial
belonging to $\Z[X,Y]$. Now  \eqref{eqn:preTMa}
can be written as $F_n(a,1)=\eta \cdot q^k$.
Thus it is sufficient to solve the
Thue--Mahler equations $F_n(X,Y)=\eta \cdot q^k$
for $(q,n,\eta)=(41,5,1)$ and $(29,7,-1)$.
	%[ 41, 5, 1 ]
	%5*X^4 - 10*X^2*Y^2 + Y^4
	%[ 29, 7, -1 ]
	%7*X^6 - 35*X^4*Y^2 + 21*X^2*Y^4 - Y^6
Explicitly these equations are
\begin{equation}\label{eqn:41k}
5 X^4 - 10 X^2 Y^2 + Y^4=41^k
\end{equation}
and
\begin{equation}\label{eqn:29k}
	7X^6 - 35X^4 Y^2 + 21X^2Y^4 - Y^6=-29^k.
\end{equation}
Using the \texttt{Magma} implementation of the Thue--Mahler solver
described in \cite{GKMS},
we find that the solutions to \eqref{eqn:41k}
are $(X,Y,k)=(\pm 2,\pm 1,1)$ and $(0,\pm 1,0)$,
and that the solutions to \eqref{eqn:29k}
are also $(X,Y,k)=(\pm 2,\pm 1,1)$ and $(0,\pm 1,0)$.
%{
%    [ 2, 1, 1 ],
%    [ 0, 1, 0 ],
%    [ 2, -1, 1 ],
%    [ -2, 1, 1 ],
%    [ -2, -1, 1 ],
%    [ 0, -1, 0 ]
%}
%{
%    [ 2, 1, 1 ],
%    [ 0, 1, 0 ],
%    [ -2, 1, 1 ],
%    [ 2, -1, 1 ],
%    [ -2, -1, 1 ],
%    [ 0, -1, 0 ]
%}
These lead to the two solutions in the statement of the lemma.
\end{proof}

\section{The equation $x^2-q^{2k}=y^n$ with $y$ even: reduction to Thue--Mahler equations}\label{sec:reduction}
Section~\ref{sec:elementary} dealt with \eqref{covid}
in the case that $y$ is odd, using purely
elementary means. We now turn our attention
to \eqref{covid} with $y$ even, and consider
the equation
\begin{equation}\label{eqn:eveneven}
x^2-q^{2k} = y^n, \qquad \text{$x$, $k$, $n$ positive integers}, \quad n \ge 3, \quad \gcd(x,y)=1, \quad \text{$y$ an even integer}.
\end{equation}
Here $q \ge 3$ is a prime. 
\begin{lem}\label{lem:n=4even}
Write $\gamma=1+\sqrt{2}$.
Any solution to \eqref{eqn:eveneven}
with $n=4$ and $q$ an odd prime must satisfy $k=1$, 
\begin{equation}\label{eqn:para2m}
	q=\frac{\gamma^{2m}+\gamma^{-2m}}{2},
	\; \; 
	x=\frac{\gamma^{4m}+6+\gamma^{-4m}}{8}
	\; \; \mbox{ and } \;  \; 
	y=\frac{\gamma^{2m}-\gamma^{-2m}}{2\sqrt{2}},
\end{equation}
for some integer $m$.
In particular, the only solutions with
$3 \le q <1000$ correspond to the identities
\[
	5^2-3^2=(\pm 2)^4, \; \; 
	145^2-17^2=(\pm 12)^4 \; \; \mbox{ and } \; \; 
	166465^2-577^2=(\pm 408)^4.
\]
\end{lem}
\begin{proof}
Suppose $n=4$. Then $(x+y^2)(x-y^2)=q^{2k}$, and so, by the coprimality of $x$ and $y$,
\[
x+y^2=q^{2k} \; \; \mbox{ and } \; \; x-y^2=1,
\]
or equivalently
\begin{equation}\label{eqn:n=4factor}
	x=\frac{q^{2k}+1}{2} \; \; \mbox{ and } \; \;  q^{2k}-2y^2=1.
\end{equation}
First we show that $k=1$. From the second equation, we have $(q^k+1)(q^k-1)=2y^2$.
Since the greatest common divisor of the two factors on the left is $2$ we see that
one of the two factors must be a perfect square, i.e. $q^k+1=z^2$ or $q^k-1=z^2$
for some non-zero integer $z$, and it is easy to see that $k$ must be odd.
The impossiblity of these cases for $k \ge 3$ follows
from Mih\u{a}ilescu's theorem \cite{Mih} (i.e. Catalan's conjecture).
Hence $k=1$.

The second equation in \eqref{eqn:n=4factor} implies that
$q+y \sqrt{2}$ is a totally positive unit in $\Z[\sqrt{2}]$. Thus
\begin{equation} \label{mickey}
q+y \sqrt{2}=\gamma^{2m} \; \; \mbox{ and } \; \;  q-y\sqrt{2}=\gamma^{-2m},
\end{equation}
for some integer $m$.  The formulae for $q$ and $y$ in \eqref{eqn:para2m} follow
from this, and the formula for $x$ follows from the first relation in \eqref{eqn:n=4factor}.

We focus on primes $3 \le q <1000$. From the first relation in \eqref{mickey},
\[
	\lvert m\rvert \; < \; 
	\frac{\log(2q)}{2 \log{\gamma}} \; < \;
	\frac{\log(2000)}{2 \log(1+\sqrt{2})} \; < 5.
\]
Thus $-4 \le m \le 4$. The
values $m=\pm 1$, $\pm 2$, $\pm 4$, 
respectively, give
the three solutions in the statement of the Lemma.
If $m=0$ or $\pm 3$, then we obtain
$q=1$ or $99$ which are not prime.
\end{proof}

In view of Lemma~\ref{lem:n=4even}, we may henceforth suppose that $n \ge 3$ is odd.
Let $x^\prime$ be either $x$ or $-x$,
chosen so that $x^\prime \equiv q^k \mod{4}$.
From \eqref{eqn:eveneven}, we deduce the existence of relatively prime integers $y_1$ and $y_2$ for which
\begin{equation}\label{eqn:factoredeven}
x^\prime+q^k=2 y_1^n \; \; \mbox{ and } \; \;  x^\prime-q^k=2^{n-1} y_2^n,
\end{equation}
with $y=2y_1y_2$, so that we have
\begin{equation} \label{easy}
y_1^n-2^{n-2} y_2^n = q^k.
\end{equation}
We have thus reduced the resolution of \eqref{eqn:eveneven} for particular $q$ and $n$ to solving 
a degree $n$ Thue--Mahler equation.
\begin{lem}\label{lem:n=3or5}
The only solutions to \eqref{eqn:eveneven} with $n \in \{ 3, 5 \}$
and $3 \le q < 1000$ an odd prime correspond to the identities
\begin{gather*}
53^2-3^2=40^3,\; 
1^2-3^2=(-2)^3,\; 
7^2-3^4=(-2)^5,\; 
29^2-5^4=6^3,\; 
9^2-7^2=2^5,\; 
43^2-11^2=12^3,\\
499^2-13^2=12^5,\;
15^2-17^2=-4^3,\;
397^2-17^4=42^3,\;
25^2-29^2=(-6)^3,\;
11^2-43^2=(-12)^3,\\
1415^2-43^2=126^3,\;
30042907^2-43^8=96222^3,\;
500047^2-47^2=6300^3,\;
55^2-53^2=6^3,\\
179^2-71^2=30^3,\;
4785^2-71^4=(-136)^3,\;
775^2-89^2=84^3,\;
155^2-101^2=24^3,\\
13609^2-109^2=570^3,\;
141^2-109^2=20^3,\;
129^2-127^2=8^3,\;
123^2-127^2=(-10)^3,\\
2789^2-127^2=198^3,\;
71^2-179^2=(-30)^3,\;
4197^2-197^2=260^3,\;
277^2-223^2=30^3,\\
249^2-251^2=(-10)^3,\;
235^2-251^2=(-6)^5,\;
255^2-257^2=-4^5,\;
118901^2-263^2=2418^3,\\
223^2-277^2=(-30)^3,\;
2355^2-307^2=176^3,\;
143027^2-307^4=2262^3,\;
557^2-307^2=60^3,\\
327^2-359^2=(-28)^3,\;
4146689^2-383^4=25800^3,\;
289^2-397^2=(-42)^3,\;
433^2-431^2=12^3,\\
431^2-433^2=(-12)^3,\;
4308693^2-433^4=26462^3,\;
979^2-479^2=90^3,\;
13^2-499^2=(-12)^5,\\
23831^2-503^2=828^3,\;
307^2-557^2=(-60)^3,\;
93^2-593^2=(-70)^3,\;
857^2-601^2=72^3,\\
713^2-659^2=42^3,\;
85016415^2-683^2=193346^3,\quad
297053^2-701^2=4452^3,\quad
731^2-727^2=18^3,\\
3655^2-739^2=234^3,\; \; 
561^2-811^2=(-70)^3,\; \; 
601^2-857^2=(-72)^3 \; \; \mbox{ and } \; \; 
1051^2-997^2=48^3.
\end{gather*}
\end{lem}
\begin{proof}
%The equation \eqref{easy} is a Thue--Mahler equation
%of degree $n$. 
For $n \in \{ 3, 5 \}$ and primes $3 \le q < 1000$,
we solved the Thue--Mahler equation  \eqref{easy}
using the \texttt{Magma} implementation of the Thue--Mahler solver
described in \cite{GKMS}. The computation resulted in the solutions given in
the statement of the lemma.
\end{proof}

\section{The modular approach to Diophantine equations: some background}\label{sec:modular}
Let $F/\Q$ be an elliptic curve over the rationals of conductor $N_F$ and minimal discriminant $\Delta_F$.
Let $n \ge 5$ be a prime. The action of $\Gal(\overline{\Q}/\Q)$ on the $n$-torsion $E[n]$ gives
rise to a $2$-dimensional mod $n$ representation
\[
\overline{\rho}_{F,n} : \Gal(\overline{\Q}/\Q) \rightarrow \GL_2(\F_n).
\]
Suppose $\overline{\rho}_{F,n}$ is irreducible (i.e. $F$ does not have an $n$-isogeny);
this can often be established by appealing to Mazur's isogeny theorem \cite{Mazur-1978}.
A standard consequence of Ribet's lowering theorem \cite{Ribet-1990}, building on the 
modularity of elliptic curves over $\Q$
due to Wiles and others \cite{Wiles}, \cite{BreuilConradDiamondTaylor01}, 
is that $\overline{\rho}_{F,n}$ arises from a weight $2$ newform of level 
\[
N \; =\; \left. {N_F} \middle/ \prod_{\stackrel{\ell \mid\mid N_F}{n \mid \ord_\ell(\Delta_F)}} \ell  \right. .
\]
More precisely, there is a newform $f$ of weight $2$ and level $N$
with normalized $\fq$-expansion
\begin{equation}\label{eqn:qexp}
f\; = \; \fq+\sum_{m=2}^\infty c_n \fq^n
\end{equation}
such that
\begin{equation}\label{eqn:cong}
\overline{\rho}_{F,n} \sim \overline{\rho}_{f,\gn}
\end{equation}
where $\gn$ is a prime ideal above $n$ of the 
ring of integers $\OO_f$ of the Hecke eigenfield
$K_f=\Q(c_1,c_2,\dotsc)$.

The original motivation for the great theorems of Ribet and Wiles included
Fermat's Last Theorem. To motivate what is to come in later sections,
we quickly sketch the deduction of FLT from the above.
Let $x$, $y$ and $z$ be non-zero coprime rational 
integers satisfying $x^n+y^n+z^n=0$ where $n \ge 5$ is prime.
After appropriately permuting $x$, $y$ and $z$, we may suppose that
$2 \mid y$ and that $x^n \equiv -1 \mod{4}$.
Let $F$ be the Frey--Hellegouarch curve
\[
Y^2=X(X-x^n)(X+y^n).
\]
It follows from Mazur's isogeny theorem and related results
that $\overline{\rho}_{E,n}$ is irreducible. A short computation reveals that
\[
\Delta_F=2^{-8} (xyz)^{2n} \; \mbox{ and } \; N_F=2 \Rad(xyz),
\]
where $\Rad(m)$ denotes the product of the prime divisors of $m$.
We find that $N=2$. Thus $\overline{\rho}_{F,n}$ 
arises from a newform $f$ of weight $2$ and level $2$; the nonexistence of such newforms provides the desired contradiction.

It is possible to use a similar strategy to treat
various Diophantine problems
including generalised Fermat equations $A x^p+By^q=Cz^r$,
for certain signatures $(p,q,r)$.
This is done by Kraus \cite{Krausppp}
for signature $(n,n,n)$ and 
by the first author and Skinner \cite{BenS}
for signature $(n,n,2)$.
Fortunately, these papers provide
recipes for the Frey--Hellegouarch curves $F$
and for the levels $N$, and establish the
required irreducibility of $\overline{\rho}_{F,n}$.
We shall make frequent use of these recipes
in later sections. It is known (and easily checked
using standard dimension formulae) that there are no weight $2$
newforms at levels
\begin{equation}\label{eqn:zerolist}
1,2,3,4,5,6,7,8,9,10,12,13,16,18,22,25,28,60,
\end{equation}
but there are newforms at all other levels. Thus, if the level $N$
predicted by the recipes is not in the list \eqref{eqn:zerolist}
then we do not immediately obtain a contradiction. 
Instead, we may compute the possible newforms using
implementations (e.g.\ in \texttt{Magma} or \texttt{SAGE})
of modular symbols algorithms due to Cremona \cite{Cre}
and Stein \cite{Stein}.
We then use the relation \eqref{eqn:cong}
to help us extract information about
the solutions to our Diophantine equation.
In doing this, we shall often make use of the following standard result, e.g. \cite{KO}, \cite[Section 5]{Siksek}.
\begin{lem}\label{lem:cong}
	Let $F/\Q$ be an elliptic curve of conductor $N_F$. Let $f$ be a weight $2$ newform of level $N$ having $\fq$-expansion as in \eqref{eqn:qexp}.
Suppose \eqref{eqn:cong} holds for some prime $n \ge 5$.
Let $\ell \ne n$ be a rational prime. 
\begin{enumerate}
\item[(i)] If $\ell \nmid  N_F N$ then $a_\ell(F) \equiv c_\ell \mod{\gn}$.
\item[(ii)] If $\ell \nmid  N$ but $\ell \mid\mid N_F$ then $\ell+1 \equiv \pm c_\ell \mod{\gn}$.
\end{enumerate}
If $f$ is a rational newform (i.e. $K_f=\Q$) then (i), (ii) also hold
for $\ell=n$.
\end{lem}

We will also make frequent use of
the following theorem of Kraus \cite[Proposition 2]{Krausppp}.
\begin{thm}[Kraus]\label{thm:Kraus}
Let $f$ be a newform of weight $2$ and level $N$ with
$\mathfrak{q}$-expansion as in \eqref{eqn:qexp}, and Hecke eigenfield
$K_f$ with ring of integers $\OO_f$.
Write
\[
M=\lcm(4,N) \; \; \mbox{ and } \; \;  \mu(M)=M \cdot \prod_{\stackrel{r \mid M}{\text{$r$ prime}}}
\left(1+\frac{1}{r}\right).
\]
Let $\frak{n}$ be a prime ideal of $\OO_f$
and suppose the following two conditions hold.
\begin{enumerate}
\item[(i)] For all primes $\ell \le \mu(M)/6$, $\ell \nmid 2N$,
we have
\[
\ell+1 \equiv c_{\ell} \mod{\frak{n}}.
\]
\item[(ii)] For all primes $\ell \le \mu(M)/6$, $\ell \mid 2N$,
$\ell^2 \nmid 4N$, we have
\[
(\ell+1) (c_\ell-1) \equiv 0 \mod{\frak{n}}.
\]
Then $\ell+1 \equiv c_{\ell} \mod{\frak{n}}$ for all primes $\ell \nmid 2N$.
\end{enumerate}
\end{thm}

%--------------------------------------------------------------------------------------
\section{The equation $x^2-q^{2k}=y^n$ with $y$ even: an approach via Frey curves}\label{sec:via}
%-------------------------------------------------------------------------------------

We are still concerned with equation \eqref{eqn:eveneven}.
In view of the results of Section~\ref{sec:reduction}, 
we may suppose that $n \ge 7$ is prime.
To show that \eqref{eqn:eveneven} has no solutions for a particular pair $(q,n)$,
it is enough to show the same for \eqref{easy}. We shall think of \eqref{easy}
as a Fermat equation of signature $(n,n,n)$ by writing it as
$y_1^n-2^{n-2} y_2^n=q^k \cdot 1^n$. This enables us to
apply recipes of Kraus \cite{Krausppp} for Frey--Hellegouarch curves
and level lowering. The following lemma will eliminate
some cases when applying those recipes.
\begin{lem}
Suppose $n \ge 7$ is prime. Then $\gcd(k,2n)=1$.
\end{lem}
\begin{proof}
Theorem 1.2 of \cite{BenS} asserts that the 
equation $A^p+2^\alpha B^p=C^2$ with prime $p \ge 7$ has no solutions
in non-zero integers with $\gcd(A,B,C)=1$ and
$\alpha \ge 2$. It immediately follows from \eqref{easy}
that $k$ is odd. Moreover, Theorem 3 of \cite{Ribet97}
asserts that the equation $A^p+2^\alpha B^p+C^p=0$
has no solutions with $ABC \ne 0$ for prime $p \ge 7$ and
$2 \le \alpha \le p-1$.
It follows that $n \nmid k$.
\end{proof}
% Note, appealing to Theorems 1.1 and 1.2 of \cite{BenS}, if $n \geq 7$ is prime, we may suppose that $k \equiv 1 \mod{2}$. 

%We note that \eqref{easy} can written as $y_1^n - 2^{n-2} y_2^n = q^k \cdot 1^n$
%and so regarded as an equation of signature $(n,n,n)$.
%We follow the recipes of Kraus \cite{Krausppp} in attaching
Following Kraus, we attach to a solution of \eqref{easy}
a Frey--Hellegouarch curve $F$, where
\begin{equation} \label{4Frey1}
F \; \; : \; \; Y^2 = X (X+y_1^n)(X+2^{n-2} y_2^n)
\end{equation}
if $q \equiv 1 \mod{4}$, and
\begin{equation} \label{4Frey2}
F \; \; : \; \; Y^2 = X (X-q^k)(X+2^{n-2} y_2^n),
\end{equation}
if $q \equiv 3 \mod{4}$.   The Frey--Hellegouarch curve $F$ is semistable,
and has minimal discriminant and conductor respectively given by
\begin{equation}\label{eqn:disccond}
\Delta_F=2^{2n-12} q^{2k} (y_1 y_2)^n \; \; \mbox{ and } \; \;  N_F=2q \cdot \Rad_2(y_1 y_2),
\end{equation}
where $\Rad_2(m)$ denotes the product of the odd primes dividing $m$.
%Similarly, if we are in case (\ref{eqn:beginning}), we put
%$$
%F_{4,1} \; \; : \; \; Y^2 = X (X+y_1^n)(X+2 q^k)
%$$
%or
%$$
%F_{4,2} \; \; : \; \; Y^2 = X (X-y_1^n)(X-2 q^k),
%$$
%depending on whether $y_1^n \equiv 1 \mod{4}$ or $y_1^n \equiv -1 \mod{4}$, respectively. 
From Kraus \cite{Krausppp}, the mod $n$ representation
of $F$ arises from a newform $f$
of weight $2$ and level $N=2q$. 
%it is in invoking the recipes of Kraus for the level
%that we make use of the fact that $n \nmid k$.

%We return to considering the Frey curve $F$ defined in \eqref{4Frey1} if $q \equiv 1 \mod{4}$
%and \eqref{4Frey2} if $q \equiv 3 \mod{4}$.
Let $\ell \nmid 2q$ be a prime.
Write
\[
T \; =\; \{a \in \Z \cap [-2\sqrt{\ell},2\sqrt{\ell}] \; : \;
a \equiv \ell+1 \mod{4}\}.
\]
%To be precise, there is a prime ideal $\frak{n} \mid n$ of $\OO_{K_f}$ such that, defining
Let
\[
\mathcal{D}_{f,\ell}^\prime=
((\ell+1)^2-c_\ell^2) \cdot
 \prod_{a \in T} (a-c_\ell),
\]
and
\[
\mathcal{D}_{f,\ell}=
\begin{cases}
\ell \cdot \mathcal{D}_{f,\ell}^\prime & \text{if $K_f \ne \Q$}\\
\mathcal{D}_{f,\ell}^\prime & \text{if $K_f=\Q$},
\end{cases}
\]
\begin{lem}\label{lem:irratbound}
Let $f$ be a newform of weight $2$ and level $2q$,
and suppose that \eqref{eqn:cong} holds. 
Let $\ell \nmid 2q$ be a prime. 
%If
%$\overline{\rho}_{F,n} \sim \overline{\rho}_{f,\frak{n}}$ 
%then 
Then $\gn \mid \mathcal{D}_{f,\ell}$.
\end{lem}
\begin{proof}
If $\ell \nmid y_1 y_2$, then $\ell \nmid N_F$ and so
is a prime of good reduction for $F$. As $F$ has full $2$-torsion we deduce that
$4 \mid (\ell+1-a_\ell(F))$. 
By the Hasse--Weil bounds,  $a_\ell(F)$ belongs to the set $T$. 
If $\ell \mid y_1 y_2$, then $\ell \mid\mid N_F$. The lemma now
follows from Lemma~\ref{lem:cong}.
\end{proof}
It is straightforward from Lemma~\ref{lem:irratbound}
and the fact that $\gn \mid n$ that $n \mid \Norm(\mathcal{D}_{f,\ell})$.
Thus if $\mathcal{D}_{f,\ell} \ne 0$, we immediately obtain an upper bound
upon the exponent $n$.
This approach will result in a bound on the exponent $n$ in
\eqref{eqn:eveneven} unless $f$ corresponds to an
elliptic curve over $\Q$ with full $2$-torsion and conductor $N=2q$; for this
see \cite[Section 9]{Siksek}.
Mazur showed that such an elliptic curve exists if and only if $q \ge 31$ is a Fermat
or a Mersenne prime; see for example \cite[Theorem 8]{Siksek}.
We note that $31$, $127$ and $257$ are the only such primes in our range $3 \le q < 1000$.
We shall exploit this approach to prove the following.
\begin{prop}\label{prop:FM}
Let $n \ge 7$  and $3 \le q <1000$ be  primes.
\begin{enumerate}
\item[(i)] If $q \not\in \{ 31, 127, 257 \}$, then \eqref{eqn:eveneven} has no solutions. 
\item[(ii)] Suppose $q \in \{ 31, 127, 257 \}$, write $q=2^m+\eta$ where $\eta=\pm 1$, and let 
%If $q \geq 31$ is a Fermat or Mersenne prime, say
%$q=2^m+(-1)^\delta$ with $\delta \in \{ 0, 1 \}$,  define an elliptic curve
%$E_q$ via \begin{equation} \label{Mersenne-Fermat}
%E_q \; \; : \; \; y^2 + xy = x^3 + (-1)^\delta 2^{m-2} x^2 + (-1)^\delta  2^{m-4} x.
%\end{equation}
%%% WARNING: the attempted minimal model above is wrong.
%%% There is no reason to give E_q in minimal form.
%If $q \geq 31$ is a Fermat or Mersenne prime, say
%$q=2^m+\eta$ with $\eta \in \{ -1, 1 \}$,  define an elliptic curve
%$E_q$ via 
\begin{equation} \label{eqn:Mersenne-Fermat}
E_q \; \; : \; \; Y^2=X(X \, +\, 1)(X \, -\, \eta \cdot 2^m).
\end{equation}
\noindent Suppose $(k,x,y)$ is a solution to \eqref{eqn:eveneven}
and let $F$ be as above. Then
%Combining  Lemma \ref{lem:irratbound} and  Theorem \ref{thm:Kraus}, we obtain the following.
%\begin{prop} \label{spec}
%Suppose that $x, y, k, q$ and $n$ are  integers satisfying (\ref{covid}), with $k$ positive, $y$ even, $q$ and $n$ prime,  $n \geq 7$,  $q < 1000$ and $q \nmid x$. Then either $n=7$ and 
%\begin{equation} \label{lantz}
%q \in \{ 43, 101, 103, 139, 163, 379,  467, 509, 557, 569,  839, 937, 947, 977 \},
%\end{equation}
%or $q \in \{ 31, 127, 257 \}$ and
%there exist coprime integers $y_1$ and $y_2$ such that $y=2y_1y_2$ and, for the corresponding Frey--Hellegouarch curves (\ref{4Frey1}) or (\ref{4Frey2}), we have
\[
\overline{\rho}_{F,n} \sim \overline{\rho}_{E_q,n}
\]
\end{enumerate}
\end{prop}
In Cremona's notation, these $E_q$ are the elliptic curves \texttt{62a2}, \texttt{254d2} and \texttt{514a2}, for $q=31$, $127$ and $257$, respectively.
\begin{proof}
There are no newforms of weight $2$ and levels $6$, $10$ and $22$.
Therefore the proof is complete in the cases where $q \in \{ 3, 5, 11 \}$.
We may thus suppose that $7 \le q<1000$ is prime and that $q \ne 11$.

For a newform $f$ of weight $2$ and level $2q$,
and a collection of primes $\ell_1,\dotsc,\ell_m$
(all coprime to $2q$), we write
$\mathcal{D}_{f,\ell_1,\dotsc,\ell_m}$ for the ideal
of $\OO_f$ generated by $\mathcal{D}_{f,\ell_1},\dotsc,\mathcal{D}_{f,\ell_m}$.
Let $\mathcal{B}_{f,\ell_1,\dotsc,\ell_m} \in \Z$ be the norm
of the ideal $\mathcal{D}_{f,\ell_1,\dotsc,\ell_m}$.
If $\overline{\rho}_{F,n} \sim \overline{\rho}_{f,\mathbf{n}}$,
then $\gn \mid \mathcal{D}_{f,\ell_1,\dotsc,\ell_m}$
by Lemma~\ref{lem:irratbound}. As $\gn \mid n$,
we deduce that $n \mid \mathcal{B}_{f,\ell_1,\dotsc,\ell_m}$.
In our computations we will take $\ell_1,\dotsc,\ell_m$
to be all the primes $<200$ distinct from $2$ and $q$,
and write $\mathcal{B}_f$ for $\mathcal{B}_{f,\ell_1,\dotsc,\ell_m}$.

We wrote a short \texttt{Magma} script which
computed, for all newforms $f$ at all levels $2q$ under consideration,
the integer
 $\mathcal{B}_f$.
We found that $\mathcal{B}_f \ne 0$
for all newforms $f$ except for three rational
newforms of levels $62$, $254$ and $514$  (corresponding to $q=31, 127$ and $257$, respectively).
% which
%corresponds to the elliptic curve \texttt{62a1}.
Thus, for all other newforms, we at least
obtain a bound on $n$. In many cases this bound
is already sharp enough to contradict our assumption that $n \geq 7$.
We give a few examples.
%In particular, for primes $q < 1000$, all newforms of level $2q$
%have $\mathcal{B}_f$  divisible only by primes
%$\le 5$ for
%$$
%\begin{array}{c} 
%q \in \left\{ 7, 17, 19, 23, 29, 47, 53, 59, 71, 79, 89, 107, 149, 179, 191, 199,  239,   \right. \\
%\left. 269, 359, 383, 431, 449, 479,  499, 599, 647, 719, 809, 863,  971 \right\} \\
%\end{array}
%$$

%For every remaining prime $13 \leq q < 1000$, there exists at least one newform $f$ for which $\mathcal{B}_f$ has a prime divisor exceeding $5$.
Let $q=13$. Then there are two eigenforms $f_1$, $f_2$
of level $2q=26$, and
\[
\mathcal{B}_{f_1}=3 \times 5, \qquad
\mathcal{B}_{f_2}=3 \times 7.
\]
Thus we eliminate $f_1$ from consideration, and also conclude
that $n=7$. It is natural to wonder if $n=7$
can be eliminated by increasing the size of our set of
primes $\ell_1,\dotsc,\ell_m$, but this is not the case.
The newform $f_2$ is rational and corresponds to the
elliptic curve $\texttt{26b1}$ with Weierstrass model
$$
E^\prime \; : \; Y^2 + X Y + Y = X^3 - X^2 - 3 X + 3.
$$
The torsion subgroup of $E^\prime(\Q)$ is isomorphic
to $\Z/7\Z$, generated by the point $(1,0)$. In particular,
for any prime $\ell \nmid 26$, we have
$7 \mid (\ell+1-a_\ell(E^\prime))$. Since $a_\ell(E^\prime)=c_\ell(f_2)$,
we have $7 \mid \mathcal{B}_{f_2,\ell}$. Thus $7 \mid \mathcal{B}_{f,\ell_1,\dotsc,\ell_n}$ regardless of the set of primes $\ell_1,\dotsc,\ell_m$
that we choose. However we can still obtain a contradiction
for $n=7$ in this case. Indeed,
we have $\overline{\rho}_{F,7} \sim \overline{\rho}_{f_2,7} \sim
\overline{\rho}_{E^\prime,7}$. Since $E^\prime$ has non-trivial
$7$-torsion, the representation $\overline{\rho}_{E^\prime,7}$
is reducible. However, the representation of the Frey
curve $\overline{\rho}_{F,7}$ is irreducible
as shown by Kraus \cite[Lemme 4]{Krausppp}, contradicting the fact that $F$ has full rational $2$-torsion.

For $q=31$, there are two newforms, $g_1$ and $g_2$.
We find that $\mathcal{B}_{g_1}=0$ and $\mathcal{B}_{g_2}=2^3 \times 3^2$;
thus we may eliminate $g_2$ for consideration.
The eigenform $g_1$ is rational and corresponds
to the elliptic curve $E_{31}$ with Cremona label \texttt{62a2}.
Hence
 $\overline{\rho}_{F,p} \sim \overline{\rho}_{g_1,p} \sim
\overline{\rho}_{E_{31},p}$, whence the proof is complete for $q=31$.

For $q=37$, there are two newforms, $h_1$ and $h_2$.
We find that $\mathcal{B}_{h_1}=3^3$ and $\mathcal{B}_{h_2}=19$.
Thus $n=19$
and
\begin{equation}\label{eqn:h2}
\overline{\rho}_{F,19} \sim \overline{\rho}_{h_2,19}.
\end{equation}
The newform $h_2$ has $\fq$-expansion
\[
h_2=\fq + \fq^2 + \alpha \fq^3 + \fq^4 + (-3\alpha - 1) \fq^5 + \alpha \fq^6 + 2 \alpha \fq^7
+\cdots, \; \; \mbox{ where } \; \;  \alpha=\frac{-1+\sqrt{5}}{2},
\]
and Hecke eigenfield $K=\Q(\sqrt{5})$. Let $\gn$ 
be the prime ideal $\gn=(4-\alpha) \cdot \OO_K$
having norm $19$.
We checked, using Theorem~\ref{thm:Kraus},
that $\ell+1 \equiv c_\ell \mod{\gn}$
for all primes $\ell \nmid 2 \cdot 37$, where $c_\ell$
is the $\ell$-th coefficient of $h_2$. From relation \eqref{eqn:h2},
we know that
\[
a_\ell(F) \equiv c_\ell \mod{\gn}
\]
for all primes $\ell$ of good reduction for $F_{3,1}$.
Thus $19 \mid (\ell+1-a_\ell(F_{3,1}))$ for all primes $\ell$
of good reduction. As before, this now implies that $\overline{\rho}_{F,19}$
is reducible \cite[IV-6]{Serre}, giving a contradiction. The proof
is thus complete for $q=37$.

The above arguments allow us to prove (ii) in the statement of the proposition,
and to obtain a contradiction for all $3 \le q <1000$, $q \not\in \{ 31, 127, 257 \}$,
except when $n=7$ and $q$ belongs to the list
%\begin{equation} \label{lantz}
\[
43,\; 101,\; 103,\; 139,\; 163,\; 379,\;  467,\; 509,\; 557,\; 569,\;  839,\; 937,\; 947,\; 977.
\]
%\end{equation}
For $n=7$ and these values of $q$, we checked using
the aforementioned Thue--Mahler solver that the only solutions
to \eqref{easy} are $(y_1,y_2,k)=(1,0,0)$. Since $k \ne 0$ in
\eqref{eqn:eveneven},
the proof is complete.
%
%The above arguments in fact allow us to deal with the remaining primes
%$41 \le q < 1000$, with the exceptions of the primes $q$ listed in (\ref{lantz}), where they fail for at at least one form $f$ in case $n=7$, and also
%for $q \in\{ 127, 257 \}$, where we find that $\overline{\rho}_{F_{3,i},n}  \sim
%\overline{\rho}_{E_{q},n}$, for $n \geq 7$.
%In all other  cases,  if $n \mid \mathcal{B}_f$
%is a prime $\ge 7$, then we can use Theorem~\ref{thm:Kraus}
%to show that $\overline{\rho}_{F_{3,i},n}$ is reducible.
%This completes the proof.
\end{proof}

\subsection*{Symplectic criteria} 
When $q \ge 31$ is a Fermat or Mersenne prime, 
%the approach employed 
%above does not yield a sharp upper bound for $n$, for 
it does not seem
to be possible, working purely with Galois representations of elliptic curves,
to eliminate the possibility that  $\overline{\rho}_{F,n} \sim \overline{\rho}_{E_q,n}$.
However, the so-called \lq symplectic method\rq\ of Halberstadt and Kraus 
\cite{HK}
allows us to derive an additional restriction on the solutions to \eqref{eqn:eveneven}.
\begin{lem}\label{lem:symplectic}
Let $q=2^m+\eta$ be a Fermat or Mersenne prime. Let $n \ge 7$ be a prime 
$\ne q$.
Suppose $(x,y,k)$ is a solution to \eqref{eqn:eveneven}, and let
$F$ be the Frey--Hellegouarch curve constructed above, and $E_q$ be given
by \eqref{eqn:Mersenne-Fermat}. Suppose $\overline{\rho}_{F,n} \sim \overline{\rho}_{E_q,n}$.
Then either $n \mid (m-4)$ or
\begin{equation} \label{symplectic}
\left( \frac{(24-6m)k}{n} \right)=1.
\end{equation}
\end{lem}
\begin{proof}
We note that the curves $F$ and $E_q$ have multiplicative reduction
at  both $2$ and $q$. Write $\Delta_1$ and $\Delta_2$ for the minimal discriminants
of $F$ and $E_q$, respectively.
By \cite[Lemme 1.6]{HK}, the ratio
\[
\frac{\ord_2(\Delta_1) \cdot \ord_q(\Delta_1)}{\ord_2(\Delta_2) \cdot \ord_q(\Delta_2)}
\]
is a square modulo $n$, provided $n \nmid \ord_2(\Delta_i)$, $n \nmid \ord_q(\Delta_i)$. It is in invoking this result of Halberstadt
and Kraus that we require the assumption that $n \ne q$.
%
%
%Suppose that $q = 2^m+ (-1)^\delta \geq 31$ ($\delta \in \{ 0, 1 \}$) is a Fermat or Mersenne prime and that we have a solution to equation (\ref{easy}) in integers $k, y_1$ and $y_2$ with $\gcd (y_1,y_2)=1$, and prime $n \geq 7$. Suppose further that we have
%$$
%\overline{\rho}_{F_{3,1+\delta},n} \sim \overline{\rho}_{E_q,n}.
%$$
%We note that the curves $F_{3,1+\delta}$ and $E_q$ have minimal discriminants
We find that
\[
\Delta_{1} = 2^{2n-12} q^{2k} (y_1y_2)^{2n}
\; \; \mbox{ and } \; \; 
\Delta_{2} = 2^{2m-8} q^2.
\]
We have previously noted that $n \nmid k$ by appealing to a result of Ribet.
Suppose $n \nmid (m-4)$. Then the valuations
$\ord_2(\Delta_i)$ and $\ord_q(\Delta_i)$ are all
indivisible by $n$. 
The result follows.
\begin{comment}
Appealing to Proposition 2 of Kraus and Oesterl\'e \cite{KO} with, in the notation of that result,  $\ell=2$ and $p=n$, we find that the representations $\overline{\rho}_{F_{3,1+\delta},n}$ and $\overline{\rho}_{E_q,n}$ are symplectically isomorphic precisely when 
$$
\left( \frac{24-6m}{n} \right)=1.
$$
Assuming further that $q \neq n$, we can again apply Proposition 2 of Kraus and Oesterl\'e \cite{KO}  (now with $\ell=q$ and $p=n$) to conclude that $\overline{\rho}_{F_{3,1+\delta},n}$ and $\overline{\rho}_{E_q,n}$ are symplectically isomorphic iff 
$$
\left( \frac{k}{n} \right)=1.
$$
It follows that, in all cases, necessarily
\end{comment}
%\begin{equation} \label{symplectic}
%\left( \frac{(24-6m)k}{n} \right)=1,
%\end{equation}
%at least provided $n \geq 7$ and $n \neq q$.
\end{proof}

%------------------------------------------------------------------------------------
\section{The equation $x^2-q^{2k}=y^n$: an upper bound for the exponent $n$ } 
\label{sec:upper}
%------------------------------------------------------------------------------------

%Consider equation (\ref{covid}),
%where $q$ and $n$ are odd primes, $k$ is a positive integer, and $x$ and $y$
%are coprime integers. We would like to solve
%this for all $3 \le q <1000$ (recall that
%we dealt already with the case $n=3$ in Section~\ref{subsec:peq3}). 
To help us complete the proof of Theorem~\ref{thm:even}, we begin by deriving an upper bound for $n$. Our approach
is essentially a minor sharpening of Theorem 3 of Bugeaud \cite{Bu-Acta} in a slightly special case. Since this result is valid for an arbitrary prime $q$, it may be of independent interest. 
\begin{thm} \label{thm:upper}
Let $x$, $y$, $q$, $k \ge 1$ and $n \geq 3$ be integers satisfying equation \eqref{covid}, with $n$ and $q$ prime, and $q \nmid x$. Then 
\[
n \; < \; 1000 \, q \log q.
\]
\end{thm}

\begin{proof}
%\subsection{$y$ odd}
If $q=2$, then we have that $n \leq 5$ from Theorem 1.2 of \cite{BenS}. We may thus suppose that $q$ is odd and, additionally, 
 that $y$ is even, or, via Proposition \ref{prop:even-odd}, we immediately obtain the much stronger result that $n \mid (q-1)$.
We are therefore in case (\ref{eqn:factoredeven}). 
By Proposition~\ref{prop:FM}, we may suppose that $q =31$ or that $q \geq 127$. 
Set $Y=\max \{ |y_1|, |2y_2| \}$ and suppose first  that
\begin{equation} \label{split}
q^k \geq Y^{n/2},
\end{equation}
or equivalently
\begin{equation} \label{split2}
2 k \log q \geq n \log Y.
\end{equation}
We set
$$
\Lambda = \frac{q^k}{(2y_2)^n}= \left( \frac{y_1}{2y_2} \right)^n - \frac{1}{4};
$$
we wish to apply an upper bound for linear forms in $q$-adic logarithms to $\Lambda$, in order to bound $k$. To do this, we must first treat the case where $y_1/2y_2$ and $1/4$ are multiplicatively dependent, i.e.\
 where $y_1y_2$ has no odd prime divisors. Under this assumption, since $y_1$ is odd, we find from \eqref{easy} that
$$
2^j \pm 1 = q^{k},
$$
for an integer $j$ with $j \equiv -2 \mod{n}$.
Via Mihailescu's theorem \cite{Mih}, if $n \geq 7$, necessarily $k=1$, $y_1=\pm 1$, $y_2=-2^\kappa$ for some integer $\kappa$ and
$$
q = 2^{(\kappa +1) n -2} \pm 1.
$$
In this case, we find a solution to equation (\ref{covid}) corresponding to the identity
$$
(-q \pm 2 )^2 - q^2 = 4 \mp 4q = \left( \mp 2^{\kappa+1} \right)^n,
$$
whereby, certainly $n < 1000 \, q \log q$.

Otherwise, we may suppose that $y_1/2y_2$ and $1/4$ are multiplicatively independent and that $Y \geq 3$. We will appeal to 
Th\'eor\`eme 4 of Bugeaud and Laurent \cite{BuLa}, with, in the notation of that result, $(\mu,\nu)=(10,5)$ (see also Proposition 1 of Bugeaud \cite{Bu-Acta}). Before we state this result, we require some notation. Let $\overline{\mathbb{Q}_q}$ denote an algebraic closure of the $q$-adic field $\mathbb{Q}_q$, and define $\nu_q$ to be the unique extension to $\overline{\mathbb{Q}_q}$ of the standard $q$-adic valuation over $\mathbb{Q}_q$, normalized so that $\nu_q(q)=1$. For any algebraic number $\alpha$ of degree $d$ over $\mathbb{Q}$,  define the {\it absolute logarithmic height} of $\alpha$ via the formula
\begin{equation}\label{eqn:htdef}
h(\alpha)= \dfrac{1}{d} \left( \log \vert a_{0} \vert + \sum\limits_{i=1}^{d} \log \max \left( 1, \vert \alpha^{(i)}\vert \right) \right), 
\end{equation}
where $a_{0}$ is the leading coefficient of the minimal polynomial of $\alpha$ over $\mathbb{Z}$ and the $\alpha^{(i)}$ are the conjugates of $\alpha$ in $\mathbb{C}$.

\begin{thm}[Bugeaud-Laurent] \label{qlog}
Let $q$ be a prime number and let $\alpha_1, \alpha_2$ denote algebraic numbers which are $q$-adic units. Let $f$ be the residual degree of the extension $\mathbb{Q}_q(\alpha_1,\alpha_2)/\mathbb{Q}_q$ and put $D=[\mathbb{Q}_q(\alpha_1,\alpha_2) : \mathbb{Q}_q]/f$. Let $b_1$ and $b_2$ be positive integers and put
$$
\Lambda_1 = \alpha_1^{b_1}-\alpha_2^{b_2}.
$$
Denote by $A_1 > 1$ and $A_2 > 1$ real numbers such that
$$
\log A_i \geq \max \left\{ h(\alpha_i), \frac{\log q}{D} \right\}, \; \; i \in \{ 1, 2 \},
$$
and put
$$
b^\prime = \frac{b_1}{D \log A_2} + \frac{b_2}{D \log A_1}.
$$
If $\alpha_1$ and $\alpha_2$ are multiplicatively independent, then we have the bound
$$
\nu_q(\Lambda_1) \leq \frac{24 q (q^f-1)}{(q-1)\log^4 (q)} \, D^4 \left( \max \left\{ \log b^\prime + \log \log q + 0.4, \frac{10 \log q}{D}, 5 \right\} \right)^2 \cdot \log A_1  \cdot \log A_2.
$$
\end{thm}

We apply this with
$$
f=1, \; D=1, \; \alpha_1=y_1/2y_2, \; \alpha_2=1/4, \; b_1=n, \; b_2=1, \;
$$
so that we may choose
$$
\log A_1 = \max \{ \log Y, \log q \}, \; \log A_2 = \max \{ 2 \log 2, \log q \}, \;
$$
and
$$
b' = \frac{n}{\log A_2} + \frac{1}{\log A_1}.
$$

Let us assume now that %$q \geq 31$ and that
\begin{equation} \label{starter}
n \geq 1000 \, q \log q,
\end{equation}
whilst recalling that either $q  = 31$ or $q \geq 127$.
We therefore have
$$
b' < 1.001  \frac{n}{\log q}
$$
and hence find that
\begin{equation} \label{igby}
k \leq 24 \,  \frac{q}{\log^3q}
 \left( \max \left\{ \log n + 0.401, 10 \log q  \right\} \right)^2 \log A_1,
\end{equation}
whence, from (\ref{split2}),
\begin{equation} \label{full}
n \log Y \leq  48 \, \frac{q}{\log^2 q}
 \left( \max \left\{  \log n + 0.401, 10 \log q  \right\} \right)^2 \log A_1.
\end{equation}

Let us suppose first that
$$
 \log n + 0.401 \geq  10 \log q.
 $$
If $q \geq Y$, we have that $\log A_1=\log q$ and hence
$$
\frac{n \log Y}{\left( \log n + 0.401 \right)^2} \leq  48 \, \frac{q}{\log q}.
$$
From (\ref{starter}), we thus have
$$
\frac{\log^2 q}{\left( \log (1000 q \log q) + 0.401 \right)^2}
\leq  \frac{0.048}{\log Y} \leq  \frac{0.048}{\log 3},
$$
contradicting $q \geq 31$. If, on the other hand, $q < Y$, then  $\log A_1=\log Y$ and so
\begin{equation} \label{pasta}
\frac{n}{\left( \log n + 0.401 \right)^2} \leq  48 \, \frac{q}{\log^2 q}.
\end{equation}
With  (\ref{starter}), this implies that
$$
\log^3 q < 0.048 \left( \log (1000 q \log q) + 0.401 \right)^2,
$$
again contradicting $q \geq 31$.

We may therefore assume that
$$
 \log n + 0.401 <   10 \log q,
 $$
 so that
 $$
 n \log Y \leq 4800 \, q \log A_1.
 $$
 If $q \geq Y$, then, from (\ref{starter}),
 $$
 \log Y < 4.8,
 $$
 whereby $3 \leq Y \leq 121$. If $|y_1| \geq 2|y_2|$, it follows from equation (\ref{easy})  that
 \begin{equation} \label{fowl}
 q^k \geq |y_1|^n - \frac{1}{4} |y_1|^n = \frac{3}{4} Y^n.
 \end{equation}
Suppose, conversely,  that $|y_1| \leq 2 |y_2|-1$ (so  that $1 \leq |y_2| \leq 60$). If $y_1>0$ and $y_2 < 0$, it follows  from (\ref{easy})  that 
 \begin{equation} \label{fowl2}
 q^k > \frac{1}{4} Y^n.
 \end{equation}
We may thus suppose that $y_1$ and $y_2$ have the same sign, whence, from (\ref{easy}),  (\ref{starter}) and $|y_2| \leq 60$, 
   \begin{equation} \label{fowl3}
 q^k = 2^{n-2} |y_2|^n - |y_1|^n > 0.24 \cdot |2y_2|^n = 0.24 \cdot Y^n.
 \end{equation}
 Combining (\ref{fowl}), (\ref{fowl2}) and (\ref{fowl3}), we thus have from (\ref{igby})  that
 $$
 n \log Y + \log (0.24) < k \log q \leq 2400 \, q \log q,
 $$
 contradicting (\ref{starter}) and $q \geq 31$. If $q < Y$, then, via (\ref{starter}),
 \begin{equation} \label{gap}
 1000 \, q \log q \leq n \leq 4800 \, q,
 \end{equation}
 a contradiction for $q \geq 127$. We may thus suppose that $q=31$, $Y > 31$ and, from (\ref{pasta}), $n \leq 12119$, which contradicts (\ref{starter}).

Next suppose that inequality (\ref{split}) (and hence  also inequality (\ref{split2}))  fails to hold.
In this case, we will apply lower bounds for linear forms in two complex logarithms.
Following Bugeaud, we take
$$
\Lambda_1= 4 \Lambda = \frac{4q^k}{(2y_2)^n}= 4  \left( \frac{y_1}{2y_2} \right)^n - 1,
$$
so that
\begin{equation} \label{janet}
\log \left| \Lambda_1 \right| =2 \log 2 + k \log q -n \log (|2y_2|).
\end{equation}

If  $Y=\max \{ |y_1|, |2y_2| \} = |y_1|$, then, from (\ref{eqn:factoredeven}), it follows that
$$
q^k \geq \frac{3}{4} |y_1|^n =  \frac{3}{4} Y^n,
$$
contradicting $q^k < Y^{n/2}$. It follows that $Y=|2y_2|$ and so, from (\ref{janet}),
\begin{equation} \label{gumball}
\log \left| \Lambda_1 \right| =2 \log 2 + k \log q -n \log Y \leq 2 \log 2 - \frac{n}{2} \log Y.
\end{equation}

From (\ref{starter}), we have that $\left| \Lambda_1 \right| \leq 1/2000$, so that
\begin{equation} \label{awful}
\left| n \log \left| \frac{2y_2}{y_1} \right| - 2 \log 2 \right| \leq \left| \log \left( 1-\Lambda_1 \right) \right| \leq 1.001 \left| \Lambda_1 \right|.
\end{equation}

We will appeal to Corollary 1 of Laurent \cite{Lau} :
\begin{thm}[Laurent] \label{LFL2}
Consider the linear form 
$$
\Lambda = c_2 \log  \beta_2 - c_1 \log \beta_1, 
$$
where $c_1$ and $c_2$ are positive integers, and $\beta_1$ and $\beta_2$ are multiplicatively independent algebraic numbers. Define
$\,D=[\mathbb{Q}(\beta_1,\beta_2) : \mathbb{Q}]\bigm/[\mathbb{R}(\beta_1,\beta_2) : \mathbb{R}]$ and set
$$
b^\prime =\frac{c_1}{D \log B_2}+ \frac{c_2}{D \log B_1},
$$
where $B_1, B_2 > 1$ are real numbers such that
$$
\log B_i \geq \max \{ h (\beta_i), |\log \beta_i|/D, 1/D \}, \; \; i \in \{ 1, 2 \}.
$$
Then
$$
\log \left| \Lambda \right| \geq -C D^4 \left( \max \{ \log b^\prime + 0.21, m/D, 1 \} \right)^2 \log B_1 \log B_2,
$$
for each pair $(m,C)$ in the following set
$$
\begin{array}{c}
 \left\{ (10,32.3), (12,29.9), (14,28.2), (16,26.9), (18,26.0), (20,25.2), \right. \\
\left.   (22,24.5), (24,24.0), (26,23.5), (28, 23.1), (30, 22.8)  \right\}.\\
\end{array}
$$
\end{thm}
Applying this result to the left hand side of (\ref{awful}), with $(m,C)=(10,32.3)$,
$$
\beta_2=|2y_2/y_1|, \; \beta_1 = 4, \; c_2=n, \; c_1=1, \; D=1,
$$
$$
\log B_2 = \log Y, \; \log B_1 = 2 \log 2 \; \mbox{ and } \; b'= \frac{n}{2 \log 2} + \frac{1}{\log Y} < \frac{1.001n}{2\log 2},
$$
we may conclude that
$$
\log \left| \Lambda_1 \right| \geq -0.001-44.8 \left( \max \left\{ \log n -0.11, 10 \right\} \right)^2  \log Y.
$$
Combining this with (\ref{gumball}), 
we thus have
$$
n  \leq 89.6  \left( \max \left\{ \log n -0.11, 10 \right\} \right)^2
+ \frac{1.4}{\log Y}.
$$
After a little work we find that
$$
n \leq 8961,
$$
contradicting (\ref{starter}) and $q \geq 31$.
%
%It remains, therefore, to treat equation (\ref{covid}) for
%$$
%q \in \{ 3, 5, 7, 11, 13, 17, 19, 23, 29 \}
%$$
%and $y$ is even. In each case, a stronger version of the desired bound is a direct consequence of applying Lemma \ref{lem:irratbound} to the Frey--Hellegouarch curve $F_{3,i}$.
\end{proof}

\begin{comment}
It remains, then, to handle the case where $y$ is even in equation \eqref{covid}, i.e. equation (\ref{easy}).
For a fixed exponent $n$, this  is a Thue--Mahler equation.
We deal with $n=5$ separately.
\begin{lem}
The only solutions to $x^2-q^{2k}=y^5$ where $3 \le q < 1000$
is prime, $\gcd(x,y)=1$ and $x$, $k$ are both positive correspond to the identities
\begin{equation}\label{eqn:solsp5}
3^4-7^2=2^5, \; \; 
499^2-13^2=12^5, \; \; 
235^2-251^2 = (-6)^5, \; \; 
255^2 - 257^2=(-4)^5. \; \; 
\end{equation}
\end{lem}
\begin{proof}
We used the Magma implementation accompanying \cite{GKMS}
to solve \eqref{easy} for all primes $3 \le q < 1000$
with $n=5$. The only solutions we obtained with $k>0$ are for $(q,k,y_1,y_2)$ one of
$$
(3,2,1,-1), \, (7,1,-1,-1), \, (13,1,-3,-2), \, (251,1,3,-1), \, (257,1,-1,-2) \mbox{ or } (499,1,3,-2).
$$
These yield the solutions in \eqref{eqn:solsp5}.
\end{proof}

\end{comment}

%---------------------------------------------------------------------------------------------------------------------------
\section{The equation $x^2-q^{2k}=y^n$: proof of Theorem~\ref{thm:even}}\label{sec:sieve}
%---------------------------------------------------------------------------------------------------------------------------

In this section, we complete the proof of Theorem~\ref{thm:even}.
Let $3 \le q <1000$ be a prime and let $(k,x,y,n)$
be a solution to \eqref{covid}
where $x$, $k \ge 1$ and $n \ge 3$ are positive integers satisfying $q \nmid x$.
Thanks to Lemmata~\ref{lem:sieve},~\ref{lem:n=4even} and~\ref{lem:n=3or5}, we may suppose that $y$ is even
and that $n \ge 7$ is prime. It follows from Proposition~\ref{prop:FM}
that $q=31$, $127$ or $257$ and $\overline{\rho}_{F,n} \sim \overline{\rho}_{E_q,n}$,
where $E_q$ is given in \eqref{eqn:Mersenne-Fermat}, and $F$ is the Frey--Hellegouarch curve
given in \eqref{4Frey1} or \eqref{4Frey2} according to whether $q \equiv 1$ or $3 \mod{4}$.
From Theorem \ref{thm:upper}, we have
\[
n \, < \, 1000 \times 257 \times \log (257) \, < \, 1.5 \times 10^6.
\]

%From Section~\ref{sec:elementary} we may suppose that $y$ is even.

%To treat the larger prime exponents in equation (\ref{easy}), we begin by noting that an initial appeal to Proposition \ref{spec} enables us to 
%restrict attention to $q \in \{ 31, 127, 257 \}$. From Theorem \ref{thm:upper}, we may suppose in each case that $n < 1.5 \times 10^6$. 
We now give a method, which for given exponent $n$ and prime $q \in \{ 31, 127, 257 \}$,
is capable of showing that \eqref{easy} 
has no solutions. This is an adaptation of the
method called \lq predicting the exponents of constants\rq\ in
\cite[Section 13]{Siksek}.
%Suppose that $q \in \{ 31, 127, 257 \}$, write $q = 2^m+ (-1)^\delta$ and 
%let $E_{q}$ be defined as in Proposition \ref{spec}.
 Let $n \ge 7$ be prime and 
choose $\ell \ne q$ to be a prime satisfying the following:
\begin{enumerate}
\item[(i)] $\ell=tn+1$ for some positive integer $t$;
\item[(ii)] $n \nmid ((\ell+1)^2-a_\ell(E_{q})^2)$.
\end{enumerate}
For $\kappa \in \F_{\ell}$, $\kappa \not\in \{ 0, 1 \}$ set
\[
E({\kappa}) \; : \; Y^2=X(X-1)(X-\kappa).
\]
Let $g$ be a primitive root for $\ell$ (i.e. a generator
for $\F_\ell^*$) and let $h=g^n$. Define $\mathcal{X}_\ell \subset \F_\ell^*$ via
\[
\mathcal{X}_\ell=\left\{ h^r/4 \; : \; 0 \le r \le t-1 \;
\text{and} \; h^r \not\equiv 4 \mod{\ell} \right\}
\]
and
\[
\mathcal{Y}_\ell=\left\{ (\kappa-1) \cdot (\F_\ell^*)^n \; : \;
 \kappa \in \mathcal{X}_\ell \quad
\text{and} \quad
a_\ell(E(\kappa))^2 \equiv a_\ell(E_{q})^2 \mod{n}
\right\} \; \subset \; \F_{\ell}^*/(\F_\ell^*)^n.
\]
Define further
\[
\phi : \Z/n\Z \rightarrow \F_{\ell}^*/(\F_\ell^*)^n \; \; \mbox{ via } \; \; 
\phi(s) =q^{s} \cdot (\F_\ell^*)^n.
\]
Finally, let 
\[
\mathcal{Z}_\ell=\left\{ s \in \phi^{-1} \left( \mathcal{Y}_\ell \right) \; : \;
\left(\frac{(24-6m)s}{n} \right) =1\right\},
\]
where $q=2^m\pm 1$; thus $m=5$, $7$ and $8$ for $q=31$, $127$ and $257$,
respectively. We note that $n \nmid (m-4)$ in all cases,
so that  \eqref{symplectic} holds.

\begin{lem}\label{lem:231}
Let $q \in \{ 31, 127, 257 \}$ and $n \geq 7$, $n \neq q$ be prime. Let $\ell_1,\dotsc,\ell_t$ be
primes $\ne q$ satisfying (i) and (ii) above, and also
\begin{equation}\label{eqn:intersection}
\bigcap_{i=1}^{t} \mathcal{Z}_{\ell_i} \; = \; \emptyset.
\end{equation}
Then equation \eqref{covid} has no solutions 
with $k \ge 1$ and $q \nmid x$.
\end{lem}
\begin{proof}
%Suppose that we have a solution to equation (\ref{easy}) in integers
%$y_1, y_2$ and $k$ with $\gcd(y_1, y_2)=1$. As noted previously, we 
%may suppose that $n \nmid  k$.
From our preceding work,  $\overline{\rho}_{F,n} \sim \overline{\rho}_{E_q,n}$. 
The minimal discriminant and conductor of $F$ are given
	in \eqref{eqn:disccond}.
	Thus a prime $\ell \nmid 2q$ satisfies
	$\ell \mid\mid N_F$ if and only if $\ell \mid y_1 y_2$,
	otherwise $\ell \nmid N_F$.
Let $\ell \ne q$ be a prime satisfying (i) and  (ii).
By (ii) we know, thanks to Lemma~\ref{lem:cong}, that $\ell \nmid y_1 y_2$, 
and so
$a_{\ell}(F) \equiv a_\ell(E_q) \mod{n}$. Let $\kappa \in \F_\ell$
satisfy 
$$
\kappa \equiv 2^{n-2} y_2^n/y_1^n \mod{\ell}. 
$$
Then $E(\kappa)/\F_\ell$
is a quadratic twist of $F/\F_\ell$ and so
$a_\ell(E(\kappa))=\pm a_\ell(F)$. We conclude
that $a_\ell(E(\kappa))^2 \equiv a_\ell(E_q)^2 \mod{n}$.

Recall that $\ell=tn+1$ and $h=g^n$, where $g$
is a primitive root of $\F_\ell$.
Observe that
\[
4 \kappa \equiv 2^n y_2^n/y_1^n
\equiv h^r \mod{\ell},
\]
for some $0 \le r \le t-1$. Moreover,
\[
\kappa -1 \; \equiv \;  \frac{2^{n-2} y_2^n}{y_1^n} - 1 \; \equiv \; - \frac{q^{k}}{y_1^n} \not\equiv 0
\mod{\ell}.
\]
In particular, $\kappa \ne 1$ and so $\kappa \in \mathcal{X}_\ell$
and $q^{k} \cdot (\F_\ell^*)^n=(\kappa-1) \cdot (\F_\ell^*)^n
\in \mathcal{Y}_\ell$. Hence $s \in \phi^{-1} (\mathcal{Y}_{\ell})$,
where $s=\overline{k} \in \Z/n\Z$.
Since  $k$ also satisfies \eqref{symplectic}, we conclude that
$s \in \mathcal{Z}_\ell$. As this is true for $\ell=\ell_1,\dotsc,\ell_t$,
the element $s$ belongs to the intersection \eqref{eqn:intersection}
giving a contradiction.
\end{proof}

\begin{cor}
For $q \in \{ 31, 127, 257 \}$ and prime $n$ with  $7 \le n < 1.5 \times 10^6$, equation
\eqref{covid} has no solutions 
with $k \ge 1$ and $q \nmid x$.
%\eqref{easy}  has no solutions in integers $k, y_1$ and $y_2$ with $\gcd(y_1,y_2)=1$, unless possibly $(q,n)=(31,7)$.
\end{cor}
\begin{proof}
For $n \neq q$, we ran a short \texttt{Magma} script that searches for
suitable primes $\ell_i$ and verifies the criterion of Lemma~\ref{lem:231}.
This succeeded for all the primes $7 \le n < 1.5 \times 10^6$ in a few minutes, except for $(q,n)=(31,7)$. In this case, we found that
$\cap \mathcal{Z}_{\ell_i}=\{\overline{1}\}$ no matter how many primes
$\ell_i$ we chose. The reason for this is that there is
a solution to equation (\ref{easy}) with $n=7$ and $k=1$, namely
$(-1)^7 - 2^5 \cdot (-1)^7=31^{1}$.

In case $n=q$, we are unable to appeal directly to Lemma~\ref{lem:231} as we no longer necessarily have (\ref{symplectic}).
We can however, still derive a slightly weaker analog of Lemma~\ref{lem:231} with the $\mathcal{Z}_\ell$ replaced by
the (typically) larger sets
\[
\mathcal{Z}_\ell^\prime=\phi^{-1} \left( \mathcal{Y}_\ell \right).
\]
For $n=q$, we find that
\[
\mathcal{Z}_{311}^\prime \cap \mathcal{Z}_{373}^\prime = \emptyset, \; \; 
\mathcal{Z}_{509}^\prime \cap \mathcal{Z}_{2287}^\prime = \emptyset \; \; 
	\mbox{ and } \; \;
\mathcal{Z}_{1543}^\prime = \emptyset, \; \; 
\]
for $q=31, 127$ and $257$, respectively.
\end{proof}

To complete the proof of Theorem \ref{thm:even}, 
it remains only to solve the Thue--Mahler equation
\[
y_1^7-32 y_2^7 = 31^k.
\]
Using the Magma implementation of \cite{GKMS}, we find  that the only solution with $k$ positive is with $k=1$ and  $y_1=y_2=-1$, corresponding to the solution $(q,k,y,n)=(31,1,2,7)$ to equation (\ref{covid}).

%----------------------------------------------------------------------------------
\section{The equation $x^2+q^{2k+1}=y^n$ with $y$ odd}\label{sec:oddodd+}
%----------------------------------------------------------------------------------
In previous sections, we have completed the proofs of Theorems~\ref{thm:even+}
and~\ref{thm:even}, therefore solving equation \eqref{Main-eq}
with $3 \le q<1000$ prime, for even exponents $\alpha$. The remainder of the 
paper is devoted to solving \eqref{Main-eq}
for odd exponents $\alpha$, and for the more modest range $3 \le q<100$.
In this section, we focus on the equation
\begin{equation}\label{eqn:odd+}
x^2+q^{2k+1}=y^n, \qquad \text{$x$, $y$, $k$ integers}, \quad k \ge 0, \quad \gcd(x,y)=1, \quad \text{$y$ odd,}
\end{equation}
with exponent $n \ge 5$ prime; here $q \ge 3$ is prime.
\begin{thm}[Arif and Abu Muriefah]\label{thm:AM}
Suppose $q \ge 3$ and $n \ge 5$ are prime, and that $n$ does not divide the class number of $\Q(\sqrt{-q})$.
Then the only solution to \eqref{eqn:odd+} corresponds to the identity
\begin{equation}\label{eqn:solodd+}
 22434^2+19 \, = \, 55^5.
\end{equation}
\end{thm}
\begin{proof}
The proof given by Arif and Abu Muriefah \cite{AM} is somewhat lengthy and slightly incorrect.
For the convenience of the reader we give a corrected and simplified proof.
Let $M=\Q(\sqrt{-q})$ and suppose that $n$ does not divide the class number of $M$.
This and the assumptions in \eqref{eqn:odd+} quickly lead us to conclude
that
\[
x+q^k \sqrt{-q}=\alpha^n
\]
for some $\alpha \in \OO_M$ with $\Norm(\alpha)=y$. Thus
\begin{equation}\label{eqn:preLucas2}
\alpha^n-\overline{\alpha}^n = 2 q^k \sqrt{-q}.
\end{equation}
If $\alpha/\overline{\alpha}$ is a root of unity, then by the coprimality
of $\alpha$ and $\overline{\alpha}$, we can conclude that $\alpha$ is a unit
and so $y=1$ giving a contradiction. Thus $\alpha/\overline{\alpha}$ is not a root of unity.
Therefore 
\[
u_m=\frac{\alpha^m-\overline{\alpha}^m}{\alpha-\overline{\alpha}}
\]
is a Lucas sequence. 
Since $\alpha \overline{\alpha}=y$, we note that $\alpha \overline{\alpha}$
is coprime to $2q$.
Suppose that the term $u_n$ has a primitive divisor $\ell$.
By definition, this is a prime $\ell$ dividing $u_n$ that does not
divide $(\alpha-\overline{\alpha})^2 \cdot u_1 u_2 \cdots u_{n-1}$.
However
$\alpha = u+v \sqrt{-q}$ or $\alpha=(u+v\sqrt{-q})/2$
where $u$, $v \in \Z$. Thus $(\alpha-\overline{\alpha})^2=-4q$
or $-q$ respectively. In particular $\ell \ne q$. It follows from \eqref{eqn:preLucas}
that $\ell=2$. By Theorem~\ref{thm:Carmichael} and the primality of $n$,
we have $n=m_2$, the rank of apparition of $\ell=2$ in the sequence $u_n$.
Again by Theorem~\ref{thm:Carmichael}, $n=m_2=2$ or $3$ contradicting
our assumption that $n \ge 5$.
It follows that $u_n$ does not have a primitive divisor.

We now invoke the Primitive Divisor Theorem (Theorem~\ref{thm:BHV})
to conclude that $n=5$ or $7$  and
that  $(\alpha,\overline{\alpha})$ is equivalent to $((a-\sqrt{b})/2,(a+\sqrt{b})/2)$
where possibilities for $(a,b)$ are given by \eqref{eqn:5list} if $n=5$, and 
by \eqref{eqn:7list} if $n=7$. For illustration, we take $n=5$
and $(a,b)= (12,-76)$. Thus $\alpha=(\pm 12 \pm \sqrt{-76})/2=\pm 6 \pm \sqrt{-19}$, whence
 $q=19$ and $y=\Norm(\alpha)=55$, quickly giving the solution in \eqref{eqn:solodd+}.
The other possibilities for $(a,b)$ in \eqref{eqn:7list} and \eqref{eqn:5list}
do not yield solutions to \eqref{eqn:odd+}.
\end{proof}

\begin{cor}\label{cor:odd+}
The only solutions to \eqref{eqn:odd+} with $3 \leq q < 100$ and $n \geq 5$ prime correspond to the identities
\[
 22434^2+19 \, = \, 55^5, \; \;  14^2+47 \, =\, 3^5 \; \; \mbox{ and } \; \;  46^2+71 \, =\, 3^7.
\]
\end{cor}
\begin{proof}
Write $h_q$ for the class number of $M=\Q(\sqrt{-q})$.
Thanks to Theorem~\ref{thm:AM},
if $n \nmid h_q$ then the only corresponding solution is $22434^2+19 \, = \, 55^5$.
Thus we may suppose that $n \mid h_q$.
The only values of $q$ in our range with $h_q$ divisible by a prime $\ge 5$ are
$q=47$, $71$ and $79$, where $h_q=5$, $7$ and $5$, respectively.
We therefore reduce to considering the three cases
$(q,n)=(47,5)$, $(71,7)$ and $(79,5)$, with $h_q=n$ in all three cases. From \eqref{eqn:odd+}, we have
\[
(x+q^k \sqrt{q} ) \cdot \OO_M \; = \; \fA^n.
\]
If $\fA$ is principal, then we are in the situation of the proof of Theorem~\ref{thm:AM}
and we obtain a contradiction. Thus $\cA$ is not principal. Now for the three quadratic
fields under consideration the class group is generated by the class $[\fP]$
where
\[
\fP \; = \; 2 \cdot \OO_M+ \frac{(1+\sqrt{-q})}{2} \cdot \OO_M
\]
is one of the two prime ideals dividing $2$. We conclude that $[\fA]=[\fP]^{-r}$
for some $1 \le r \le n-1$. Observe that $\fC \overline{\fC}$ is principal
for any ideal $\fC$ of $\OO_M$, so $[\overline{\fC}]=[\fC]^{-1}$.
We choose $\fB=\fA$ or $\overline{\fA}$
so that $[\fB]=[\fP]^{-r}$ where $1 \le r \le (n-1)/2$. We note that
\[
(x \pm q^k \sqrt{-q}) \cdot \OO_M \; = \;  \fB^n \; = \; (\fP^{-n})^r \cdot (\fP^r \fB)^n
\]
where the $\pm$  sign is $+$ if $\fB=\fA$ and $-$ is $\fB=\overline{\fA}$. We note that 
both $\fP^{-n}$ and $\fP^r \fB$ are principal. We find that
$\fP^{-n}=2^{-n-1} (u+v \sqrt{-q}) \cdot \OO_M$ where $u$, $v$ are 
\[
(u,v)=\begin{cases}
                (-9, 1) & \text{if $q=47$,}\\
                (-21, 1) & \text{if $q=71$,}\\
                (7, 1) & \text{if $q=79$}.
\end{cases}
\]
The ideal $\fP^r \fB$ is integral as well as principal, and so 
has the form $(X^\prime+Y^\prime \sqrt{-q}) \cdot \OO_M$
where $X^\prime$ and $Y^\prime$ are either both integers,
or both halves of odd integers. 
We conclude that
\[
2^{s+rn+r} (x \pm q^k \sqrt{-q} )\; = \; (u+v \sqrt{-q})^r \cdot (X+Y \sqrt{-q})^n
\]
where $X$, $Y \in \Z$ and $s=0$ or $n$. Equating imaginary parts gives
\[
G_r(X,Y) = \pm 2^{s+rn+r} q^k
\]
where $G_r \in \Z[X,Y]$ is a homogeneous polynomial of degree $n$.
We solved this Thue--Mahler equation using the Thue--Mahler solver
associated to the paper \cite{GKMS},
for each of our three pairs $(q,n)$ and each $0 \le r \le (n-1)/2$.
For illustration, we consider the case $q=47$, $n=5$, $r=2$. Thus
$(u,v)=(-9,1)$. We find 
\[
	G_2(X,Y)= 2 
	(-9 X^5 + 85 X^4 Y + 4230 X^3 Y^2 - 7990 X^2 Y^3 - 99405 X Y^4 +
	37553 Y^5)
\]
and are therefore led to  solve the Thue--Mahler equation
\[
	-9 X^5 + 85 X^4 Y + 4230 X^3 Y^2 - 7990 X^2 Y^3 - 99405 X Y^4 +
	37553 Y^5 \; = \; \pm 2^j q^k.
\]
We find that the solutions are 
\[
	(X,Y,j,k)\; =\; (1, 1, 16, 0 ) \; \; \mbox{ and } \; \; (-1, -1, 16, 0),
\]
and compute $G_2(1,1)=-2^{17}$, $G_2(-1,-1)=2^{17}$. We note that
$17=n+rn+r$, therefore $s=n=5$. 
We deduce that
\[
	x\pm 47^{k} \sqrt{-47} = \pm (-9+\sqrt{-47})^2 \cdot (1+\sqrt{-47})^5=
	\pm (14-\sqrt{-47}).
\]
Thus $x=\pm 14$ and $k=0$, giving the solution 
$14^2+47 \, =\, 3^5$. 
The other cases are similar.
\end{proof}

%----------------------------------------------------------------------------------------------------
\section{The equation $x^2-q^{2k+1}=y^5$}\label{sec:n=5}
%---------------------------------------------------------------------------------------------------

We will soon apply Frey--Hellegouarch curves to study the
equation $x^2 - q^{2k+1}=y^n$ for prime exponents $n \ge 7$,
and for $q$ a prime in the range $3 \le q < 100$. In Section~\ref{sec:tiny},
we have solved this equation for $n \in \{ 3, 4 \}$. This leaves
only exponent $n=5$ which we now treat through reduction to Thue--Mahler equations.
\begin{lem}\label{lem:odd5+}
Let $3 \le q < 100$ be a prime. The only solutions to the equation
\begin{equation}\label{eqn:odd5-}
x^2-q^{2k+1}=y^5, \qquad \text{$x$, $y$, $k$ integers}, \quad k \ge 0, \quad \gcd(x,y)=1,
\end{equation}
correspond to the identities
\begin{gather*}
2^2-3=1^5, \; 
2^2-5=(-1)^5, \;
10^2-7^3=(-3)^5, \;
56^2-11=5^5,\;
16^2-13=3^5, \\
4^2-17=(-1)^5, \;
7^2-17=2^5, \;
3^2-41=(-2)^5, \;
411^2-41^3=10^5 \; \mbox{ and } \;
11^2-89=2^5.
\end{gather*}
\end{lem}
%[ 3, 0, 2, 1 ]
%[ 5, 0, 2, -1 ]
%[ 7, 1, 10, -3 ]
%[ 11, 0, 56, 5 ]
%[ 13, 0, 16, 3 ]
%[ 17, 0, 4, -1 ]
%[ 17, 0, 7, 2 ]
%[ 41, 1, 411, 10 ]
%[ 41, 0, 3, -2 ]
%[ 89, 0, 11, 2 ]
\begin{proof}
Let $M=\Q(\sqrt{q})$. 
For $q$ in our range, the class number of $M$ is $1$, unless $q=79$
in which case the class number is $3$. Suppose first that $y$ is odd. Then
\[
(x+q^k\sqrt{q}) \OO_M= \fA^5
\]
where $\fA$ is an ideal of $\OO_M$. Since the class number is not divisible by $5$,
we see that $\fA$ is principal and conclude that
\begin{equation}\label{eqn:preThue}
x+q^{k} \sqrt{q} = \epsilon^r \cdot \alpha^5,
\end{equation}
where $\epsilon$ is some fixed choice of a fundamental unit for $M$,
$-2 \le r \le 2$, and $\alpha \in \OO_M$. Note that 
\[
-x+q^{k} \sqrt{q} = \epsilon^{-r} \cdot \beta^5,
\]
where $\beta$ is one of $\pm \overline{\alpha}$. Thus we may without loss
of generality suppose that $0 \le r \le 2$. The case $r=0$ is easily shown
not to lead to any solutions by following the approach in the proof
of Theorem~\ref{thm:AM}. Thus we suppose $r=1$ or $2$. 

Let 
\[
\theta=\begin{cases}
\sqrt{q} & \text{if $q \equiv 3 \mod{4}$}\\
(1+\sqrt{q})/2 & \text{if $q \equiv 1 \mod{4}$}.
\end{cases}
\]
Then $\{ 1, \theta \}$ is a $\Z$-basis for $\OO_M$ and so we may write $\alpha=X+Y \theta$ where $X$, $Y \in \Z$. 
It follows that
\[
\epsilon^r \cdot \alpha^5=F_r(X,Y)+G_r(X,Y) \theta,
\]
where $F_r$, $G_r$ are homogeneous degree $5$ polynomials in $\Z[X,Y]$.
Equating the coefficients of $\theta$ in \eqref{eqn:preThue} yields
the Thue--Mahler equations
\[
G_r(X,Y)=\begin{cases}
q^k & \text{if $q \equiv 3 \mod{4}$}\\
2q^k & \text{if $q \equiv 1 \mod{4}$}.
\end{cases}
\]
Solving these equations for prime $3 \le q <100$
and for $r \in \{ 1, 2 \}$ leads to the solutions given
in the statement of the theorem with $y$ odd.

Next we consider the case when $y$ is even, so that $q \equiv 1 \mod{8}$.
The possible values of $q$ in our range
are $17$, $41$, $73$, $89$ and $97$ (whence, in each case, $M$ has class number $1$).
We can rewrite the equation $x^2-q^{2k+1}=y^5$ as
\[
\left(\frac{x+q^k \sqrt{q}}{2} \right) \left(\frac{x-q^k\sqrt{q}}{2} \right) = 2^3 y_1^5
\]
where $y_1=y/2$. The two factors on the left-hand side are coprime. Let $\beta$
be a generator of 
\[
\fP=2\OO_M + \left(\frac{1+\sqrt{q}}{2}\right) \cdot \OO_M
\]
which is one of the two prime ideals above $2$. After possibly replacing $x$ by $-x$ we obtain
\[
\frac{x-q^k}{2}+q^k \theta \; = \; \frac{x+q^k \sqrt{q}}{2} \; = \; \epsilon^r \beta \alpha^5
\]
where $-2 \le r \le 2$. Writing $\alpha=X+Y \theta$ and equating the coefficients of $\theta$
on both sides gives, for each choice of $q$ and $r$, a Thue--Mahler equation. Solving these 
leads to the solutions in the statement of the theorem with $y$ even.
\end{proof}

%-----------------------------------------------------------------------
\section{Frey--Hellegouarch curves for a ternary equation of signature $(n,n,2)$} \label{sec:FreyH}
%-----------------------------------------------------------------------

In studying equation \eqref{covid}, we employed a factorisation
argument which reduced to \eqref{easy} (which in turn we treated
as a special case of a Fermat equation having signature $(n,n,n)$).
In the remainder of the paper, we are primarily interested
in the equation $x^2+(-1)^\delta q^{2k+1}=y^n$,
where $q$ is a prime. We shall treat this, for prime $n \ge 7$,
as a Fermat equation of signature $(n,n,2)$
by rewriting this as $y^n+q^{2k+1} (-1)^{(\delta+1)n} =x^2$,
a special case of 
\begin{equation} \label{eqn:major}
y^n + q^\alpha z^n = x^2, \qquad \qquad \gcd(x,y)=1.
\end{equation}
Equation \eqref{eqn:major} has previously been studied
by Ivorra and Kraus \cite{IvKr}, and by the first author and Skinner \cite{BenS}.
In this section, we recall some of these results and strengthen
them slightly before specialising them to the case $z=\pm 1$
in forthcoming sections.
%As noted earlier, 
%a starting point for our work is to view a solution to equation (\ref{LebNag}) as corresponding to
%a solution to a ternary equation of the shape $y^n - D \cdot 1^n = x^2$. In case $D$ is divisible by a single prime $q$, machinery exists in the literature for treating such equations somewhat more generally.
%In particular, the main Th\'eor\`eme of Ivorra and Kraus \cite{IvKr} is as follows.
\begin{thm}[Ivorra and Kraus] \label{general}
Suppose that $q$ is a prime with the property that $q$ cannot be written in the form
$$
q = |t^2 \pm 2^k|,
$$
where $t$ and $k$ are integers, with  $k=0$, $k=3$ or $k \geq 7$. Then there are no solutions to the Diophantine equation
\eqref{eqn:major}
in integers $x, y, z, n$ and $\alpha$ with %$\gcd(x,y)=1$ and 
$n$ prime satisfying
\begin{equation} \label{pess}
n > \left( \sqrt{8(q+1)}+1 \right)^{2(q-1)}.
\end{equation}
\end{thm}

To verify whether or not a given prime $q$ can be written as $|t^2-2^k|$, an old result of Bauer and the first author \cite{BaBe} can be helpful.
We have, from Corollary 1.7 of \cite{BaBe}, if $t$ and $k$ are positive integers with $k \geq 3$ odd,
$$
\left| t^2 - 2^k \right| > 2^{13k/50},
$$
unless
$$
(t,k) \in \{ (3,3), (181,15) \}.
$$
In particular, a short computation reveals that Theorem \ref{general} is applicable to the following primes $q < 100$ :
\begin{equation} \label{good-boys}
q \in \{ 11, 13,  19, 29, 43, 53, 59, 61, 67, 83 \}.
\end{equation}
We shall make Theorem~\ref{general}
more precise for these particular values of $q$.
To this end we attach to a solution of \eqref{eqn:major} a certain
Frey--Hellegouarch curve, following the recipes of
the first author and Skinner.
If $yz$ is even in \eqref{eqn:major}, then we define, assuming, without loss of generality,  that $x \equiv 1 \mod{4}$,
\begin{equation}\label{eqn:Fyeven}
F \; \; : \; \; Y^2+XY = X^3+\left( \frac{x-1}{4} \right) X^2 + \frac{y^n}{64} X, \; \; \mbox{ if $y$ is even},
\end{equation}
and
\[
F \; \; : \; \; Y^2+XY = X^3+\left( \frac{x-1}{4} \right) X^2 + \frac{q^\alpha z^n}{64} X, \; \; \mbox{ if $z$ is even}.
\]
If, on the other hand, $yz$ is odd, we define
\begin{equation}\label{eqn:Fy1mod4}
F \; \; : \; \; Y^2 = X^3+2x X^2 + q^\alpha z^n X
\end{equation}
or
\begin{equation}\label{eqn:Fy3mod4}
F \; \; : \; \; Y^2 = X^3+2x X^2 + y^n X,
\end{equation}
depending on whether $y \equiv 1 \mod{4}$ or $y \equiv -1 \mod{4}$, 
respectively.
Let
\begin{equation}\label{eqn:kappa}
	\kappa=\begin{cases}
		1 & \text{if $yz$ is even}\\
		5 & \text{if $yz$ is odd}.
	\end{cases}
\end{equation}
By the results of \cite{BenS}, in each case, we may suppose that $n \nmid \alpha$ and that the mod $n$ representation
of $F$ arises from a newform $f$
of weight $2$ and level $N=2^\kappa \cdot q$.
%\[
%N_i=\begin{cases}
%2q & \mbox{ if } i=1\\
%32q & \mbox{ if } i=2.
%\end{cases}
%\]
Let the $\fq$-expansion of $f$ be given by \eqref{eqn:qexp}.
As before, we denote the Hecke eigenfield by $K_f=\Q(c_1,c_2,\dotsc)$
and its ring of integers by $\OO_f$. In particular, there
is a prime ideal $\gn$ of $\OO_f$ 
such that \eqref{eqn:cong} holds.
%\begin{equation}\label{eqn:qexp-first}
%f=\mathfrak{q}+\sum_{m=1}^\infty c_m \mathfrak{q}^m
%\end{equation}
%for the $\fq$-expansion of $f$, and let $K=\Q(c_1,c_2,\ldots)$ be its Hecke eigenfield.
%The coefficients $c_i$ belong to the ring of integers $\OO_K$ and there is a prime ideal $\gn \mid n$ of $\OO_K$
%such that
%\begin{equation}\label{eqn:compare-first}
%\overline{\rho}_{F_{i,j},n} \sim \overline{\rho}_{f,\gn}.
%\end{equation}
Let $\ell \nmid 2q$ be prime and
\[
T=\{a \in \Z \cap [-2\sqrt{\ell},2\sqrt{\ell}] \; : \;
a \equiv 0 \mod{2}\}.
\]
We write
\[
\mathcal{D}_{f,\ell}^\prime=
((\ell+1)^2-c_\ell^2) \cdot
 \prod_{a \in T} (a-c_\ell),
\]
and
\[
\mathcal{D}_{f,\ell}=
\begin{cases}
\ell \cdot \mathcal{D}_{f,\ell}^\prime & \text{if $K_f \ne \Q$}\\
\mathcal{D}_{f,\ell}^\prime & \text{if $K_f=\Q$}.
\end{cases}
\]
\begin{lem}\label{lem:irratbound-first}
Let $f$ be a newform of weight $2$ and level $N=2^\kappa \cdot q$.
Let $\ell \nmid 2q$ be a prime. If
$\overline{\rho}_{F,n} \sim \overline{\rho}_{f,\gn}$
then $\gn \mid \mathcal{D}_{f,\ell}$.
\end{lem}
\begin{proof}
%This is a standard sort of result
%and follows in a straightforward fashion from \eqref{eqn:compare-first},
%using Propositions 5.1 and 5.2 of \cite{Siksek}.
The proof is almost identical to the proof of Lemma~\ref{lem:irratbound}.
The only difference is the definition of $T$
which takes into account the fact $F$ has
a single rational point of order $2$ instead of full $2$-torsion.
\end{proof}

The following is a slight refinement of  Theorem 1.3 of \cite{BenS}. 
\begin{prop} \label{general100}
Suppose that $q$ belongs to \eqref{good-boys}.
Then there are no solutions to  equation \eqref{eqn:major}
in integers $x, y, z, n$ and $\alpha$ with $\gcd(x,y)=1$ and $n \geq 7$ prime, unless, possibly, $n=7$ and 
$q \in \{ 29, 43, 53, 59, 61 \}$, or one of the following holds
$$
q=11, \, n=7 \mbox{ and } yz \equiv 1 \mod{2}, \mbox{ or }
$$
$$
q=19, \, n=7 \mbox{ and } yz \equiv 1 \mod{2}, \mbox{ or }
$$
$$
q=43, \, n =11 \mbox{ and } yz \equiv 1 \mod{2}, \mbox{ or }
$$
$$
q=53, \, n = 17 \mbox{ and } yz \equiv 1 \mod{2}, \mbox{ or }
$$
$$
q=59, \, n = 11  \mbox{ and } yz \equiv 0 \mod{2}, \mbox{ or }
$$
$$
q=61, \, n = 13  \mbox{ and } yz \equiv 1 \mod{2}, \mbox{ or }
$$
$$
q=67, \, n \in \{ 7, 11, 13, 17 \} \mbox{ and } yz \equiv 1 \mod{2}, \mbox{ or }
$$
$$
q=83, \, n=7 \mbox{ and } yz \equiv 1 \mod{2}.
$$
\end{prop}
\begin{proof}
For a weight $2$ newform $f$ of level $N$ and primes $\ell_1,\dots,\ell_m$
(all coprime to $2q$), 
write
$\mathcal{D}_{f,\ell_1,\dotsc,\ell_m}$ for the ideal
of $\OO_f$ generated by $\mathcal{D}_{f,\ell_1},\dotsc,\mathcal{D}_{f,\ell_m}$.
Let $\mathcal{B}_{f,\ell_1,\dotsc,\ell_m} \in \Z$ be the norm
of the ideal $\mathcal{D}_{f,\ell_1,\dotsc,\ell_m}$.
If $\overline{\rho}_{F,n} \sim \overline{\rho}_{f,\mathbf{n}}$
then $\gn \mid \mathcal{D}_{f,\ell_1,\dotsc,\ell_m}$
by Lemma~\ref{lem:irratbound-first}.
Write $\mathcal{B}_{f,\ell_1,\dotsc,\ell_m}=\Norm(\mathcal{D}_{f,\ell_1,\dotsc,\ell_m})$.
Thus $n \mid \cB_{f,\ell_1,\dotsc,\ell_m}$.
In our computations, we take
$\ell_1,\dotsc,\ell_m$ to be the primes
$<100$ coprime to $2q$, and we let $\cB_f=\cB_{f,\ell_1,\dotsc,\ell_m}$.
If $\cB_f \ne 0$, then we certainly have a bound on $n$.
If $\cB_f$ is divisible only by primes $\le 5$, then 
we know that \eqref{eqn:cong} does not hold
for that particular $f$, and we can eliminate
it from further consideration.

For primes $q$ in \eqref{good-boys},
we apply this with  newforms $f$ of levels $N =2^\kappa q$, $\kappa \in \{ 1, 5 \}$. 
We obtain the desired conclusion that equation \eqref{eqn:major} has no
solutions provided $n \geq 7$ is prime,  unless $q
\in \{ 29, 43, 53, 59, 61 \}$ and $n=7$, or $(q,n,\kappa)$ is one of 
$$
\begin{array}{l}
(11,7,5),\; (13,7,1),\; (19,7,5),\; (43,11,1),\;  (43,11,5),\;  (53,17,5),\;
(59,11,1),\; (61,31,1),  \\
 (61,13,5),\;   (67,17,1),\; (67,7,5),\; (67,11,5),\; (67,13,5),\; 
(67,17,5),\; (83,7,1),\; (83,7,5). 
\end{array}
$$

We show that the triples $(13,7,1)$, $(43,11,1)$, $(61,31,1)$, $(67,17,1)$ and
$(83,7,1)$ do not have corresponding solutions; the remaining triples
lead to the noted possible exceptions. 
For illustration, take $q=83$ and $\kappa=1$, so that $N=2 \times 83=166$.
There are three conjugacy classes of weight $2$ newforms of level $N$,
which we denote by $f_1$, $f_2$, $f_3$, which respectively have
Hecke eigenfields $\Q$, $\Q(\sqrt{5})$ and $\Q(\theta)$
where $\theta^3-\theta^2-6 \theta+4=0$.
We find
\[
	\cB_{f_1}=3^2 \times 5, \qquad \cB_{f_2}=5, \qquad \cB_{f_3}=7.
\]	
We therefore reduce that $f=f_3$ and $n=7$.
%To apply this result to treat the case $(q,n,\kappa)=(83,7,1)$, for example, we take $N=2 \cdot 83$ and (using Magma's labelling for Galois conjugacy classes of newforms), $f=f_3$. We note that the forms $f_1$ and $f_2$ are eliminated from consideration through application of Lemma \ref{lem:irratbound-first}.
In fact, $\cD_f=(7,3+\theta)$ is a prime ideal above $7$, so we take
$\gn=(7,3+\theta)$.
%We let $\gn$  be the prime ideal $\gn=(3+\theta) \cdot \OO_K$, where $\theta^3-\theta^2-6 \theta+4=0$, so that $\mu(M)/6 =84$. 
A short calculation 
verifies the congruences in hypotheses (i) and (ii) of Theorem~\ref{thm:Kraus},
whence
%	reveals that, for each prime $\ell$, $3 \leq \ell < 83$, we have
%$\ell+1 \equiv c_\ell \mod{\gn}$, whence, from Theorem~\ref{thm:Kraus}, in fact 
$\ell+1 \equiv c_\ell \mod{\gn}$ for all $\ell$ with $\ell \nmid 2 \cdot 83$. It follows from Lemma~\ref{lem:cong}
that
\[
a_\ell(F) \equiv c_\ell \mod{\gn}
\]
for all primes $\ell$ of good reduction for $F$ and hence $7 \mid (\ell+1-a_\ell(F))$ for all such  primes $\ell$
of good reduction.
 This now implies that $\overline{\rho}_{F,7}$
is reducible \cite[IV-6]{Serre}, giving a contradiction.

We argue similarly for
$$
(q,n,\kappa) = (13,7,1), \; (43,11,1), \; (61,31,1) \mbox{ and } (67,17,1).
$$
In each case, Lemma \ref{lem:irratbound-first} eliminates all but one class of newforms %($f_2, f_2, f_3$ and $f_2$, respectively) 
which are in turn treated via Theorem \ref{thm:Kraus}.
\end{proof}

For other odd primes $q < 100$, outside the set (\ref{good-boys}), we can, in
certain cases,  still show that equation \eqref{eqn:major} has no nontrivial
solutions for suitably large $n$, under the additional assumption that $yz
\equiv 0 \mod{2}$ or, for other $q$, under the assumption that $yz \equiv 1
\mod{2}$. To be precise, we have \begin{prop} \label{general100-even}
Suppose that
$q \in \{ 3, 5, 37, 73  \}$.
Then there are no solutions to  equation \eqref{eqn:major}
in integers $x, y, z, n$ and $\alpha$ with $yz \equiv 0 \mod{2}$, $\gcd(x,y)=1$ and $n \geq 7$ prime, unless, possibly, $(q,n)=(73,7)$.
\end{prop}
\noindent and
\begin{prop} \label{general100-odd}
Suppose that
$q \in \{ 23, 31, 47, 71, 79, 97 \}$.
Then there are no solutions to  equation \eqref{eqn:major}
in integers $x, y, z, n$ and $\alpha$ with $yz \equiv 1 \mod{2}$, $\gcd(x,y)=1$ and $n \geq 7$ prime, unless, possibly,
$n=7$ and $q \in \{ 23, 31, 47, 71, 97 \}$, or $(q,n)=(79,11)$, or $(q,n)=(97,29)$.
\end{prop}

As in the case of Proposition \ref{general100}, these results follow after a small amount of computation, by applying Lemma \ref{lem:irratbound-first} and Theorem \ref{thm:Kraus}.

%---------------------------------------------
\section{The equation $x^2\pm q^{2k+1}=y^n$ and proofs
of Theorems~\ref{thm:odd} and~\ref{thm:modest}} \label{Sec14}
%--------------------------------------

We now specialize and improve on the results of Section~\ref{sec:FreyH}, proving
the following.
\begin{prop}\label{prop:primary}
Let $(x,y,k)$ be a solution to the equation
\begin{equation}\label{eqn:primary}
	x^2+(-1)^\delta q^{2k+1} = y^n, \qquad
	\delta \in \{0,1\}, \quad k \ge 0, \quad \gcd(x,y)=1
\end{equation}
where $q$ is a prime in the range
$3 \le q <100$, and $n \ge 7$ is prime.
Suppose moreover that
\begin{enumerate}
\item[(a)] if $y$ is odd then $\delta=1$;
\item[(b)] if $\delta=1$ then $q \not\in \{ 3, 5, 17, 37 \}$.
\end{enumerate}
If $y$ is even, suppose without loss of generality that $x \equiv 1 \mod{4}$.
Write
\begin{equation}\label{eqn:kappaNew}
	\kappa=\begin{cases}
		1 & \text{if $y$ is even}\\
		5 & \text{if $y$ is odd}.
	\end{cases}
\end{equation}
	Let $v \in\{0,1\}$ satisfy $k \equiv v \mod{2}$.
Attach to the solution $(x,y,k)$ the Frey--Hellegouarch curve
\[
	G=G_{x,k} \; : \; 
	\begin{cases}
		Y^2=X^3+4x X^2+ 4(x^2+(-1)^\delta q^{2k+1}) X 
		& \text{if $\kappa=1$},\\
		Y^2=X^3-4x X^2+ 4(x^2+(-1)^\delta q^{2k+1}) X 
		& \text{if $\kappa=5$ and $q \equiv (-1)^\delta \bmod{4}$},\\
		Y^2=X^3+2xX^2+(x^2+(-1)^\delta q^{2k+1})X 
		& \text{if $\kappa=5$ and $q \equiv (-1)^{\delta+1} \bmod{4}$}.
	\end{cases}
\]
Then either $n >1000$ and $\overline{\rho}_{G,n} \sim \overline{\rho}_{E,n}$
where $E/\Q$ is an elliptic curve of conductor $2^\kappa q$ given in Table~\ref{table:ells}
or the solution $(x,y,k)$ corresponds to one of the identities
\begin{gather*}
11^2+7=2^7, \;
45^2+23=2^{11}, \;
13^2-41=2^7, \;
9^2+47=2^7, \\
7^2+79=2^7, \;
91^2-89=2^{13} \; \mbox{ or } \; 
15^2-97=2^7.
%[ 23, 0, 0, 45, 2, 11 ]
%[ 41, 1, 0, 13, 2, 7 ]
%[ 47, 0, 0, 9, 2, 7 ]
%[ 79, 0, 0, 7, 2, 7 ]
%[ 89, 1, 0, 91, 2, 13 ]
%[ 97, 1, 0, 15, 2, 7 ]
\end{gather*}
\end{prop}

\begin{table}
\caption{Data for Proposition~\ref{prop:primary}.
Here the elliptic curves $E$ are
given by their Cremona labels.} \label{table:ells}
\centering
\begin{minipage}[t]{0.49\textwidth}%\centering
{
\tabulinesep=1.2mm
\begin{tabu}{|l|l|l|l|l|}
\hline
$q$ & $\delta$ & $\kappa$ & $v$ & $E$ \\
\hline\hline
$7$ &  $0$ & $1$ & $0$ & \texttt{14a1}\\
\hline
$7$ &  $0$ & $1$ & $1$ & \texttt{14a1}\\
\hline
$23$ & $0$ & $1$ & $0$ & \texttt{46a1}\\
\hline
$31$ & $0$ & $1$ & $0$ & \texttt{62a1}\\
\hline
$31$ & $0$ & $1$ & $1$ & \texttt{62a1}\\
\hline
$41$ & $1$ & $1$ & $0$ & \texttt{82a1} \\ 
\hline
$41$ & $1$ & $5$ & $0$ & \texttt{1312a1}, \texttt{1312b1}\\
\hline
$41$ & $1$ & $5$ & $1$ & \texttt{1312a1}, \texttt{1312b1}\\
\hline
\end{tabu}
}
\end{minipage}
\hskip -20mm
\begin{minipage}[t]{0.49\textwidth}%\centering
{
\tabulinesep=1.2mm
\begin{tabu}{|l|l|l|l|l|}
\hline
$q$ & $\delta$ & $\kappa$ & $v$ & $E$ \\
\hline\hline
$47$ & $0$ & $1$ & $0$ & \texttt{94a1}\\
\hline
$71$ & $0$ & $1$ & $0$ &\texttt{142c1}\\
\hline
$71$ & $0$ & $1$ & $1$ &\texttt{142c1}\\
\hline
$73$ & $1$ & $5$ & $0$ & \texttt{2336a1}, \texttt{2336b1}\\
\hline
$73$ & $1$ & $5$ & $1$ & \texttt{2336a1}, \texttt{2336b1}\\
\hline
$79$ & $0$ & $1$ & $0$ & \texttt{158e1}\\
\hline
$89$ & $1$ & $1$ & $0$ & \texttt{178b1}\\
\hline
$97$ & $1$ & $1$ & $0$ & \texttt{194a1}\\
\hline
\end{tabu}
}
\end{minipage}
\end{table}

% [* [* 7, 0, 1, 0, 14a1 *], [* 7, 0, 1, 1, 14a1 *], 
%[* 23, 0, 1, 0, 46a1 *], 
%[* 31, 0, 1, 0, 62a1 *], 
%[* 31, 0, 1, 1, 62a1 *], 
%[* 41, 1, 1, 0, 82a1 *], 
%[* 41, 1, 5, 0, 1312a1 *], [* 41, 1, 5, 0, 1312b1 *], [* 41, 1, 5, 1, 1312a1 *], [* 41, 1, 5, 1, 1312b1 *], 
%[* 47, 0, 1, 0, 94a1 *], 
%[* 71, 0, 1, 0, 142c1 *], [* 71, 0, 1, 1, 142c1 *], 
%[* 73, 1, 5, 0, 2336a1 *], [* 73, 1, 5, 0, 2336b1 *], 
%[* 73, 1, 5, 1, 2336a1 *], [* 73, 1, 5, 1, 2336b1 *], 
%[* 79, 0, 1, 0, 158e1 *], [* 89, 1, 1, 0, 178b1 *], [* 97, 1, 1, 0, 194a1 *] *]

Before proceeding to the proof of this result, we make 
a few remarks on the assumptions in Proposition~\ref{prop:primary}.
Our eventual goal is to prove Theorems~\ref{thm:everything},
\ref{thm:odd} and \ref{thm:modest},
and thus we are interested in the equation $x^2+(-1)^\delta q^\alpha=y^n$
where $3 \le q < 100$. Theorems~\ref{thm:even+} and~\ref{thm:even}
(proved in Sections~\ref{Thm4} and~\ref{sec:sieve}, respectively)
treat the case where $\alpha$ is even, so we 
are reduced to $\alpha=2k+1$. 
The results of Sections~\ref{sec:tiny},
Corollary~\ref{cor:odd+} and Lemma~\ref{lem:odd5+} allow us
to restrict the exponent $n$ to be a prime $\ge 7$.
Thanks to Theorem~\ref{thm:AM},
we need not consider the case where $\delta=0$ and $y$ is odd,
which explains the reason for assumption (a).
With a view to proving the proposition, we will soon provide a method
which is usually capable, for a fixed $q$, $\delta$ and $n$, of showing
that \eqref{eqn:primary} does not have a solution. 
If $\delta=1$, and $q$ is one of the values $3$, $5$, $17$ or $37$,
then there is a solution to \eqref{eqn:primary}
for all odd values of the exponent $n$:
\[
2^2-3=1^n, \; \;  2^2-5=(-1)^n, \; \; 4^2-17=(-1)^n \; \; \mbox{ and } \; \; 6^2-37=(-1)^n,
\]
and so our method fails if $\delta=1$ and $q$ is one of these
four values; this explains assumption (b) in the statement of the proposition.

We note that \eqref{eqn:primary} is a special case of
\eqref{eqn:major} with $z$ specialised to the value $(-1)^{\delta+1}$,
and with $\alpha=2k+1$. 
The value $\kappa$ in the statement of the proposition
agrees with value for $\kappa$ in \eqref{eqn:kappa}
given in the previous section.
We note that if $y$ is odd, then
$y \equiv (-1)^\delta \cdot q \mod{4}$.
The Frey--Hellegouarch curve $G$ is, up to isogeny, the same as the Frey--Hellegouarch curve
$F$ in the previous section, but is more convenient for our
purposes. More precisely, the model $G$ is isomorphic to $F$ given in \eqref{eqn:Fyeven}
if $y$ even (i.e. $\kappa=1$), and to $F$ given in \eqref{eqn:Fy3mod4}
if $y \equiv 3 \mod{4}$ (i.e. $\kappa=5$ and $q \equiv (-1)^{\delta+1} \bmod{4}$). It is $2$-isogenous to $F$ in
\eqref{eqn:Fy1mod4} if $y \equiv 1 \mod{4}$ (i.e. $\kappa=5$ and $q \equiv (-1)^\delta \bmod{4}$).
Thus $\overline{\rho}_{F,n} \sim \overline{\rho}_{G,n}$ in all three cases. 
%The model $G$ is simply more convenient
%for the sieving arguments we will apply later. 
We conclude from the previous section that $\overline{\rho}_{G,n} \sim \overline{\rho}_{f,\gn}$
where $f$ is a weight $2$ newform of level $N=2^\kappa q$.

Note that if $\kappa=1$ (i.e. $y$ is even) then $1+(-1)^\delta q \equiv 0 \mod{8}$. This together with
the assumptions of Proposition~\ref{prop:primary} 
shows that we are concerned with $30$ possibilities for  
the triple $(q,\delta,\kappa)$, namely
%\begin{gather*}
\begin{equation}\label{eqn:triples}
\begin{cases}
    ( 7, 0, 1 ), \;
    ( 7, 1, 5 ), \;
    ( 11, 1, 5 ), \;
    ( 13, 1, 5 ),\;
    ( 19, 1, 5 ),\;
    ( 23, 0, 1 ),\;
    ( 23, 1, 5 ),\;
    ( 29, 1, 5 ),\\
    ( 31, 0, 1 ),\;
    ( 31, 1, 5 ),\;
    ( 41, 1, 1 ),\;
    ( 41, 1, 5 ),\;
    ( 43, 1, 5 ),\;
    ( 47, 0, 1 ),\;
    ( 47, 1, 5 ),\;
    ( 53, 1, 5 ),\\
    ( 59, 1, 5 ),\;
    ( 61, 1, 5 ),\;
    ( 67, 1, 5 ),\;
    ( 71, 0, 1 ),\;
    ( 71, 1, 5 ),\;
    ( 73, 1, 1 ),\;
    ( 73, 1, 5 ),\;
    ( 79, 0, 1 ),\\
    ( 79, 1, 5 ),\;
    ( 83, 1, 5 ),\;
    ( 89, 1, 1 ),\;
    ( 89, 1, 5 ),\;
    ( 97, 1, 1 ),\;
    ( 97, 1, 5 ).
\end{cases}
\end{equation}
%\end{gather*}

\subsection*{Bounding the exponent $n$}
In the previous section we defined 
an ideal $\cD_{f,\ell_1,\dotsc,\ell_r}$ which if non-zero
allows us to bound the exponent $n$ in \eqref{eqn:major}.
That bound will also be valid for \eqref{eqn:primary}
since it is a special case of \eqref{eqn:major}.
We now offer a refinement that is often capable of yielding
a better bound for \eqref{eqn:primary}.

Fix a triple $(q,\delta,\kappa)$ from the above list.
%and one of the three cases $(\delta,\kappa)=(0,1)$, $(1,1)$, $(1,5)$.
We also fix $v \in \{0,1\}$ and suppose that $k \equiv v\mod{2}$. 
Let $f$ be a weight $2$ newform of level $N=2^\kappa q$
with $\fq$-expansion as in \eqref{eqn:qexp}. Write $K_f$ for the Hecke eigenfield of $f$, and $\OO_f$
for the ring of integers of $K_f$.
For a prime $\ell\ne 2$, $q$, define 
\[
\cS_{\ell} \; =\; \{\; a_\ell(G_{w,v}) \quad : \quad w \in \F_\ell, \quad w^2+(-1)^\delta q^{2v+1} \not \equiv 0 \mod{\ell} \;\}.
\]
Let
\[
\cT=\cT_{\ell} = \begin{cases}
\cS_{\ell} \cup \{\ell+1,-\ell-1\} & \text{if $(-1)^{\delta+1} q$ is a square modulo $\ell$}\\
\cS_{\ell} & \text{otherwise}.\\
\end{cases}
\]
Let
\[
\cE^\prime_{\ell}=\prod_{a \in \cT} (a-c_\ell) \; \; \mbox{ and } \; \; 
\cE_{\ell} = 
\begin{cases}
\ell \cdot \cE^\prime_{\ell} & \text{if $K_f \ne \Q$},\\
\cE^\prime_{\ell} & \text{if $K_f=\Q$},
\end{cases}
\]
where, as before, $c_\ell$ is the $\ell$-th coefficient in the $\fq$-expansion of $f$.
\begin{lem}\label{lem:cE}
Let $\gn$ be a prime ideal of $\OO_f$ above $n$.
If $\overline{\rho}_{G,n}  \sim \overline{\rho}_{f,\gn}$ then
$\gn \mid \cE_{\ell}$. 
\end{lem}
\begin{proof}
Write $k=2u+v$ with $u \in \Z$. Let $w \in \F_\ell$ satisfy
$w \equiv x/q^{2u} \mod{\ell}$. Hence
\[
y^n=x^2+(-1)^\delta q^{2k+1} \equiv q^{4u} \cdot (w^2+(-1)^\delta q^{2v+1}) \mod{\ell}.
\]
It follows that $\ell \mid y$ if and only if $w^2+(-1)^\delta q^{2v+1} \mod{\ell}$.
Suppose first that $\ell \nmid y$. The elliptic curves $G_{x,k}/\F_\ell$ and $G_{w,v}/\F_\ell$
are isomorphic, and so $a_\ell(G_{x,k})=a_\ell(G_{w,v})$. In particular,
$a_\ell(G_{x,k}) \in \cT_{\ell}$ and so $a_\ell(G_{x,k})-c_\ell$ divides
$\cE_{\ell}$. Likewise, if $\ell \mid y$ (which can only happen
if $(-1)^{\delta+1} q$ is a square modulo $\ell$) then
$(\ell+1)^2 - c_\ell^2$ divides $\cE_{\ell}$.
The lemma follows from Lemma~\ref{lem:cong}.
\end{proof}

\subsection*{A sieve}
Lemma~\ref{lem:cE} will soon allow us to eliminate most possibilities for the newform 
$f$ in a manner similar to Propositions~\ref{general100},~\ref{general100-even} and~\ref{general100-odd}.
We will still need to treat some cases for fixed exponent $n$.
To this end, we will employ a sieving technique similar to the one
in Section~\ref{sec:sieve}.

\medskip

Fix a prime $n \ge 7$, and let $\gn$ be a prime ideal of $\OO_f$ above $n$. Let $\ell \ne q$ be a prime. Suppose
\begin{enumerate}
\item[(i)] $\ell=tn+1$ for some positive integer $t$;
\item[(ii)] either $\gn \nmid (4-c_\ell^2)$, or $(-1)^{\delta+1}q$ is not a square modulo $\ell$.
\end{enumerate}
Let 
\begin{gather*}
A=\{m \in \{0,1,\dotsc,2n-1\} \quad : \quad m \equiv v \mod{2}, \quad n \nmid (2m+1) \},\\
\cX_{\ell}=\{(z,m) \in \F_\ell \times A \; : \; (z^2+(-1)^\delta q^{2m+1})^t \equiv 1 \mod{\ell} \},\\
\cY_\ell=\{(z,m) \in \cX_{\ell} \; : \;   a_\ell(G_{z,m}) \equiv c_\ell \mod{\gn}\}, \\
\cZ_\ell=\{ m \; : \; \text{there exists $z$ such that $(z,m) \in \cY_\ell$}\}.
\end{gather*}
\begin{lem}\label{lem:crudeKraus}
Let $\ell_1,\dotsc,\ell_r$ be primes $\ne q$ satisfying (i), (ii). 
Let
\[
\cZ_{\ell_1,\dotsc,\ell_r}=\bigcap_{i=1}^r \cZ_{\ell_i}.
\]
If $\overline{\rho}_{G,n} \sim \overline{\rho}_{f,\gn}$
then
\[
(k \bmod {2n}) \in \cZ_{\ell_1,\dotsc,\ell_r}.
\]
\end{lem}
\begin{proof}
Let $m$ be the unique element of  $\{0,1,\dotsc,2n-1\}$ satisfying $k \equiv m \mod{2n}$.
Let $\ell \ne q$ be a prime satisfying (i) and (ii). 
It is sufficient
to show that $m \in \cZ_\ell$.
First we will demonstrate that $\ell \nmid y$.
If $(-1)^{\delta+1} q$ is not a square modulo $\ell$
then $\ell \nmid y$ from \eqref{eqn:primary}. Otherwise,
by (ii), $\gn \nmid (4-c_\ell^2)$. However,
from (i) and the fact that $\gn \mid n$ we have
$\ell+1 \equiv 2 \mod{\gn}$ and so $\gn \nmid ((\ell+1)^2-c_\ell^2)$.
It follows from Lemma~\ref{lem:cong} that $\ell$ is a prime of good reduction for 
$G_{x,k}$ and so $\ell \nmid y$. We deduce from Lemma~\ref{lem:cong}
that $a_\ell(G_{x,k}) \equiv c_\ell \mod{\gn}$.

In the previous section, we observed that $n \nmid \alpha$
in equation~\ref{eqn:major} thanks to the results of \cite{BenS}, whence
$n \nmid (2k+1)$.
Since $k \equiv v \mod{2}$, we know that $m \in A$.
Write $k=2nb+m$ with $b$ a non-negative integer and let
$z \in \F_\ell$ satisfy $z \equiv x/q^{2nb} \mod{\ell}$.
Then
\[
z^2+(-1)^\delta q^{2m+1} \; \equiv \; \frac{1}{q^{4nb}} (x^2+(-1)^\delta q^{2k+1}) \; \equiv \; (y/q^{4b})^n \mod{\ell}.
\]
From (i), we deduce that
\[
(z^2+(-1)^\delta q^{2m+1})^t \; \equiv \; (y/q^{4b})^{\ell-1} \; \equiv \; 1 \mod{\ell}.
\]
Thus $(z,m) \in \cX_{\ell}$. Moreover, we have that $G_{x,k}/\F_\ell$ and $G_{z,m}/\F_\ell$
are isomorphic elliptic curves, whence $a_\ell(G_{z,m}) =a_\ell(G_{x,k}) \equiv c_\ell \mod{\gn}$. 
Thus $(z,m) \in \cY_\ell$ and so $m \in \cZ_\ell$ as required.
\end{proof}
\noindent \textbf{Remarks.}
We would like to explain how to compute $\cZ_\ell$ efficiently, given $n$ and $\ell$.
\begin{enumerate}
\item In our computations, the value $t$ will be relatively small compared to $n$ and to $\ell=tn+1$.
Let $g$ be a primitive root modulo $\ell$ (i.e. a cyclic generator for $\F_\ell^\times$),
and let $h=g^n$.
The set $\cX_\ell$ consists of pairs $(z,m) \in \F_\ell \times A$
such that $(z^2+(-1)^\delta q^m)^t \equiv 1 \mod{\ell}$.
Hence $z^2+(-1)^\delta q^m$ is one of the values $1,h,h^{2},\dotsc,h^{(t-1)}$.
Thus, to compute $\cX_\ell$, we run through $i=0,1,\dotsc,t-1$ and $m \in A$ and solve
$z^2=h^i - (-1)^\delta q^m$. We note that the expected cardinality
of $\cX_\ell$ should be roughly $t \times \#A \approx t \times n \approx \ell$.
\item 
It seems at first that, in order to compute $\cY_\ell$ and $\cZ_\ell$,
we need to compute $a_\ell(G_{z,m})$ for all $(z,m) \in \cX_\ell$,
and this might be an issue for large $\ell$.
There is in fact a shortcut that often means
that we need to perform few of these computations. 
In fact we will need to compute $\cZ_\ell$ for large values
of $\ell$ 
only for rational newforms $f$
that correspond to elliptic curves $E/\Q$
with non-trivial $2$-torsion.
In this case, 
we note that $a_\ell(G_{z,m}) \equiv a_\ell(E) \mod{2}$,
as both elliptic curves have non-trivial $2$-torsion.
If $(z,m) \in \cY_\ell$, then
$a_\ell(G_{z,m}) \equiv a_\ell(E) \mod{2n}$.
However, by the Hasse--Weil bounds,
\[
\lvert a_\ell(G_{z,m}) - a_\ell(E) \rvert \;  \le \; 4 \sqrt{\ell}.
\]
Suppose  moreover that $n^2>4\ell$ (which will be usually satisfied
as $t$ is typically small). 
Then, 
the congruence
$a_\ell(G_{z,m}) \equiv c_\ell=a_\ell(E) \mod{2n}$
is equivalent to the equality $a_\ell(G_{z,m})=a_\ell(E)$,
and so to 
$\# G_{z,m}(\F_\ell)=\# E(\F_\ell)$. 
To check whether the equality
$\# G_{z,m}(\F_\ell)=\# E(\F_\ell)$ holds
for a particular pair $(z,m) \in \cX_\ell$,
we first 
choose a random point $Q \in G_{z,m}(\F_\ell)$ and
check whether $\#E(\F_\ell) \cdot Q=0$.
Only for pairs $(z,m) \in \cX_\ell$ that pass this test do
we need to compute $a_\ell(G_{z,m})$ and check
the congruence $a_\ell(G_{z,m}) \equiv a_\ell(E) \mod{n}$. 
\end{enumerate}

\subsection*{A refined sieve}
We note that if $\cZ_{\ell_1,\dotsc,\ell_r}=\emptyset$
then $\overline{\rho}_{G,n} \nsim \overline{\rho}_{f,\gn}$.
In our computations, described later,
we are always able to find suitable
primes $\ell_1,\dotsc,\ell_r$ satisfying (i), (ii),
so that $\cZ_{\ell_1,\dotsc,\ell_r}=\emptyset$, at least
for $n$ suitably large. For smaller values
of $n$ (say less than $50$), we occasionally failed.
We now describe a refined sieving method
that, whilst being somewhat slow, has a better chance
of succeeding for those smaller values of the exponent $n$.

Let $(q,\delta,\kappa)$ be one of our $30$ triples given in \eqref{eqn:triples}, and let $n\ge 7$ be
a prime.
Suppose that $(x,y,k)$ is a solution
to \eqref{eqn:primary} where 
$y$ is even if and only if $\kappa=1$.
Let $\phi=\sqrt{(-1)^{\delta+1} q}$ and set $M=\Q(\phi)$.
Let $\fP$ be one of the prime ideals of $\OO_M$ above $2$.

Our first goal is to produce a finite set $\cS \subset M^*$,
such that
\begin{equation}\label{eqn:descent}
x+q^k \phi = \gamma \cdot \alpha^n
\end{equation}
for some $\gamma \in \cS$ and $\alpha \in \OO_M$.
This is the objective of Lemmata~\ref{lem:cSkappa=5}
and~\ref{lem:cSkappa=1}. Both of these
make an additional assumption on the class group,
but this assumption will in fact be satisfied
in all cases where we need to apply our 
refined sieving idea.
\begin{lem}\label{lem:cSkappa=5}.
Let $\kappa=5$.
Suppose that the class group $\Cl(\OO_M)$ of $\OO_M$ 
is cyclic and generated by the class $[\fP]$.
Let $h=\# \Cl(\OO_M)$ and set 
\[
\cI = \{ 0 \le i \le h-1 \; : \; \text{$\fP^{-ni}$ is principal}\}.
\]
Choose for each $i \in \cI$ a generator $\beta_i$
for $\fP^{-ni}$. Let $\epsilon$ be a fundamental
unit for $M$ (recall that if $\kappa=5$
then $\delta=1$ and so $M$ is real). 
Let 
\[
\cS=\{ \epsilon^j \beta_i \; : \; -(n-1)/2 \le j \le (n-1)/2, \;  i\in \cI\}.
\]
Then there is some $\gamma \in \cS$ and $\alpha \in \OO_M$
such that \eqref{eqn:descent} holds.
Moreover, $\Norm(\alpha)=2^\mu y$ for some $\mu \ge 0$.
\end{lem}
\begin{proof}
As $\kappa=5$, we have that $y$ is odd. Then
\[
(x+q^k \phi ) \OO_M = \fA^n,
\]
where $\fA$ is an ideal of $\OO_M$ with norm $y$.
Since $[\fP]$ generates the class group, the same is true of 
$[\fP]^{-1}$. 
Hence $[\fA]=[\fP]^{-i}$ for some $i \in \{0,1,\dotsc,h-1\}$.
Now
\[
(x+q^k \theta) \OO_M = \fP^{-ni} \cdot ( \fP^i \cdot \fA)^n.
\]
Since $\fP^i \cdot \fA$ is principal, it follows that  $\fP^{-ni}$ is also principal. 
The lemma obtains.
\end{proof}
\begin{lem}\label{lem:cSkappa=1}
Let $\kappa=1$.
Suppose that the class group $\Cl(\OO_M)$ of $\OO_M$ 
is cyclic and generated by the class $[\fP]$.
Let $h=\# \Cl(\OO_M)$
and set
\[
\cI = \{ 0 \le i \le h-1 \; : \; \text{$\fP^{n(1-i)-2}$ is principal}\}.
\]
Choose for each $i \in \cI$ a generator $\beta_i$
for $\fP^{n(1-i)-2}$. 
Let
\[
\cS^\prime=\{ \beta_i : i \in \cI   \} \cup 
\{\overline{\beta_i} : i \in \cI\},
\]
where $\overline{\beta_i}$ denotes
the Galois conjugate of $\beta_i$.
Let
\[
\cS = \begin{cases}
\{ 2 \cdot \beta \; : \;   \beta \in \cS^\prime\} & \text{if $\delta=0$}\\
\{ 2 \cdot \epsilon^j \cdot \beta \; : \; -(n-1)/2 \le j \le (n-1)/2, \;  \beta \in \cS^\prime\} & \text{if $\delta=1$},
\end{cases}
\]
where $\epsilon$ is a fundamental
unit for $M$.
Then there is some $\gamma \in \cS$ and $\alpha \in \OO_M$
such that \eqref{eqn:descent} holds.
Moreover, $\Norm(\alpha)=2^\mu y$ for some $\mu \in \Z$.
\end{lem}
\begin{proof}
As $\kappa=1$, we have that $y$ is even. 
Then
\[
\left(\frac{x+q^k \phi }{2}\right) \OO_M = \fC^{n-2}\fA^n
\]
where $\fA$ is an ideal of $\OO_M$ with norm $y/2$
and $\fC$ is one of $\fP$, $\overline{\fP}$.
Since $[\fP]$ generates the class group so does
$[\fC]^{-1}$. 
Hence $[\fA]=[\fC]^{-i}$ for some $i \in \{0,1,\dotsc,h-1\}$.
Now
\[
\left(\frac{x+q^k \phi }{2}\right) \OO_M = \fC^{n(1-i)-2} \cdot (\fC^i \fA)^n.
\]
But $\fC^i \cdot \fA$ is principal, whence $\fC^{n(1-i)-2}$ is principal,
and so $i \in \cI$ and $\fC^{n(1-i)-2}$ is generated by 
either $\beta_i$ or $\overline{\beta_i}$.
The lemma follows.
\end{proof}

We will now describe our refined sieve.
Fix $m \in \{0,1,\dotsc,2n\}$ and suppose $k \equiv m \mod{2n}$.
Let $\gn$ be a prime ideal of $\OO_f$ above $n$. Let $\ell \ne q$ be a prime. Suppose
\begin{enumerate}
\item[(a)] $\ell=tn+1$ for some positive integer $t$;
\item[(b)] $\gn \nmid (4-c_\ell^2)$;
\item[(c)] $(-1)^{\delta+1} q$ is a square modulo $\ell$.
\end{enumerate}
We choose an integer $s$ such that $s^2 \equiv (-1)^{\delta+1} q \mod{\ell}$.
Let
\[
\fL = \ell \OO_M + (s-\phi) \OO_M.
\]
By the Dedekind--Kummer theorem $\ell$ splits in $\OO_M$ and $\fL$ is one of the two
prime ideals above $\ell$. In particular, $\OO_M/\fL \cong \F_\ell$
and $\phi \equiv s \mod{\fL}$.
Let
\begin{gather*}
\cX_{\ell,m}=\{z \in \F_\ell  \; : \; (z^2+(-1)^\delta q^{2m+1})^t \equiv 1 \mod{\ell} \},\\
\cY_{\ell,m}=\{z \in \cX_{\ell,m} \; : \;   a_\ell(G_{z,m}) \equiv c_\ell \mod{\gn}\}, \\
\cU_{\ell,m} = \{ (z,\gamma) \; : \; z \in \cY_{\ell,m}, \quad
\gamma \in \cS, \quad (z+q^m \phi)^t \equiv \gamma^t \mod{\fL} \},\\
\cW_{\ell,m}=\{ \gamma \; : \; \text{there exists $z$ such that $(z,\gamma) \in \cU_{\ell,m}$}\}.
\end{gather*}
\begin{lem}\label{lem:refinedKraus}
Let $\ell_1,\dotsc,\ell_r$ be primes $\ne q$ satisfying (a), (b) and (c) above. 
Let
\[
\cW=\cW_{\ell_1,\dotsc,\ell_r}=\bigcap_{i=1}^r \cW_{\ell_i}.
\]
If $\overline{\rho}_{G,n} \sim \overline{\rho}_{f,\gn}$,
then there is some $\gamma \in \cW$ and 
some $\alpha \in \OO_M$ such that
\eqref{eqn:descent} holds.
\end{lem}
\begin{proof}
Suppose $\ell$ satisfies conditions (a), (b) and (c).
As $\ell$ satisfies (a) and (b), it also satisfies hypotheses (i) and (ii)
preceding the statement of Lemma~\ref{lem:crudeKraus}.
Write $k=2nb+m$ where $b$ is a non-negative integer, and let 
$z \equiv x/q^{2nb} \mod{\ell}$.
It follows from the proof of Lemma~\ref{lem:crudeKraus}
that $\ell \nmid y$ and that $z \in \cY_{\ell,m}$.
We know from Lemmata~\ref{lem:cSkappa=5} and~\ref{lem:cSkappa=1} that there is some $\gamma \in \cS$
such that
$x+q^k \phi = \gamma \alpha^n$
where $\alpha \in \OO_M$ satisfies $\Norm(\alpha)=2^\mu  y$ for some $\mu \in \Z$.
Note that $\gamma$ is supported only on the prime ideals above $2$. 
Since $\fL \mid \ell$, we have $\ord_\fL(\alpha)=\ord_\fL(\gamma)=0$.
Hence
\[
z+q^{m} \phi \equiv \frac{1}{q^{2nb}} (x+q^k \phi) \equiv \gamma \cdot \left(\frac{\alpha}{q^{2b}}\right)^n \mod{\cL}.
\]
Since $(\OO_M/\cL)^* \cong \F_\ell^*$ is cyclic of order $\ell-1=tn$, we have
\[
(z+q^m \phi)^t \equiv \gamma^t \mod{\cL}.
\]
Thus $(z,\gamma) \in \cU_{\ell,m}$ and hence $\gamma \in \cW_{\ell,m}$.
The lemma follows.
\end{proof}

\subsection*{Proof of Proposition~\ref{prop:primary}}
Our proof of Proposition~\ref{prop:primary}
is the result of applying \texttt{Magma}
scripts based on Lemmata~\ref{lem:cE},~\ref{lem:crudeKraus} and ~\ref{lem:refinedKraus},
as well as solving a few Thue--Mahler equations. 
Our approach subdivides the proof 
into $60$ cases
corresponding to $60$ quadruples $(q,\delta,\kappa,v)$:
here
$(q,\delta,\kappa)$ is one of the $30$ triples
in \eqref{eqn:triples},
and $v \in \{0,1\}$. 
Let $(x,y,k)$ be a solution to
\eqref{eqn:primary} with prime exponent $n \ge 7$. Suppose that $y$ is even
if $\kappa=1$ and $y$ is odd if $\kappa=5$. Suppose
moreover that $k \equiv v \mod{2}$.
Our first step is to compute the newforms $f$ of weight $2$
and level $N=2^\kappa q$. We know that for one these newforms $f$,
we have $\overline{\rho}_{G,n} \sim \overline{\rho}_{f,\gn}$
where $G=G_{x,k}$ is the Frey--Hellegouarch curve given in Proposition~\ref{prop:primary},
and $\gn \mid n$ is a prime ideal of $\OO_f$, the ring of integers of the
Hecke eigenfield $K_f$. 
Let $p_1,\dotsc,p_s$ be the primes $\le 200$
distinct from $2$, $q$ and let 
\[
\cE_{f}=\sum_{i=1}^s \cE_{p_i},
\]
where $\cE_{p_i}$ is as in Lemma~\ref{lem:cE}.
It follows from Lemma~\ref{lem:cE} that if 
$\overline{\rho}_{G,n} \sim \overline{\rho}_{f,\gn}$
then $\gn \mid \cE_{f}$, and so $n \mid \Norm(\cE_{f})$.

We illustrate this by taking $(q,\delta,\kappa,v)=(31,1,5,0)$.
There are $8$ newforms $f_1,\dotsc,f_8$
of weight $2$ and level $2^\kappa q=992$, which all happen to be
irrational. We find that
\[
\Norm(\cE_{f_j})=7,~7,~2^{10},~2^{10},~2^3,~2^3,~2^6 \times 3^2,~2^6 \times 3^2,
\]
respectively for $j=1,2,\dotsc,8$.
Thus $n=7$ and $f=f_1$ or $f_2$. We consider first 
\[
f=f_1=
\fq + \sqrt{2}\fq^3 - \fq^5 - (1+\sqrt{2}) \fq^7 - \fq^9 + 2(1-\sqrt{2}) \fq^{11} + \cdots, 
\]
with Hecke eigenfield $K_f=\Q(\sqrt{2})$ having ring of integers 
$\OO_f=\Z[\sqrt{2}]$. We found that $\cE_{f}=(1+2\sqrt{2})$
which is one of the two prime ideals above $7$.
Hence $\gn=(1+2\sqrt{2})$. Next we compute
$\cZ=\cZ_{\ell_1,\dotsc,\ell_{30}}$ as in Lemma~\ref{lem:crudeKraus}
where $\ell_1,\dotsc,\ell_{30} \ne 31$ are the $30$ primes 
satisfying (i) and (ii) with $t \le 200$.
We find that $\cZ=\{0,8\}$.
Thus, by Lemma~\ref{lem:crudeKraus},
we have $k \equiv 0$ or $8 \mod{14}$. 
Now for $m=0$ and $m=8$, we compute $\cW=\cW_{\ell_1,\dotsc,\ell_{36}}$
as in Lemma~\ref{lem:refinedKraus}, where
$\ell_1,\dotsc,\ell_{36} \ne 31$ are the $36$ primes
satisfying (a), (b) and (c) with $t \le 800$.
We found that $\cW=\emptyset$ for $m=0$
and that $\cW=\{\epsilon^3\}$ for $m=8$ where $\epsilon=1520+273\sqrt{31}$
is the fundamental unit of $M=\Q(\sqrt{31})$.
Hence we conclude, by Lemma~\ref{lem:refinedKraus},
that $k \equiv 8 \mod{14}$ and that
\[
x+31^k \sqrt{31} = (1520+273\sqrt{31})^3(X+Y\sqrt{31})^7,
\] 
for some integers $X$, $Y$. Equating the coefficients of $\sqrt{31}$ on both sides
results in a degree $7$ Thue--Mahler equation with huge coefficients.
However, using an algorithm of Stoll and Cremona \cite{StollCremona}
for reducing binary forms we discover that
this Thue--Mahler equation can be rewritten as
\[
31^k=
-56 U^7 + 112 U^6 V - 84 U^5 V^2 + 140 U^4 V^3 + 490 U^3 V^4 + 1596 U^2 V^5 + 
    2807 U V^6 + 2119 V^7,  
\]
where $U$, $V \in \Z$ are related to $X$, $Y$ via the unimodular substition
\[
U=2 X + 11 Y \; \; \mbox{ and } \; \; 
 V=   7 X + 39 Y.
\]
We applied the Thue--Mahler solver to this and found that it has
no solutions. Next we take $f=f_2$ which also has
Hecke eigenfield $K_f=\Q(\sqrt{2})$. 
We apply Lemmata~\ref{lem:cE},~\ref{lem:crudeKraus} and~\ref{lem:refinedKraus}
using the same sets of primes $p_j$ and $\ell_i$ as for $f_1$.
We
find $\cE_f=(1-2\sqrt{2})$, and so $\gn=(1-2\sqrt{2})$ and $n=7$.
Again we obtain $\cZ=\{0,8\}$ on applying Lemma~\ref{lem:crudeKraus}.
We find that
$\cW=\emptyset$ for $m=0$ and $\cW=\{\overline{\epsilon}^3\}$
for $m=8$. Again the corresponding Thue--Mahler
equation has no solutions. Thus equation \eqref{eqn:primary}
has no solutions with $n \ge 7$ prime
for $q=31$, $\delta=1$ and with $y$ odd (i.e. $\kappa=5$)
and $k \equiv 0 \mod{2}$.
We used the above approach to deal with
all the cases where $\cE_f$ is non-zero.
In all the cases where $\cE_f=0$,
the newform $f$ is rational, and in fact
corresponds to an elliptic curve $E/\Q$
with non-trivial $2$-torsion. These
elliptic curves are listed in Table~\ref{table:ells}.
Thus $\overline{\rho}_{G,n} \sim \overline{\rho}_{E,n}$.
What is required for Proposition~\ref{prop:primary}
is to show in these cases that
there are no solutions with prime 
$7 \le n < 1000$ apart from the ones listed
in the statement of the proposition.
We illustrate how this works by taking
$(q,\delta,\kappa,v)=(7,0,1,0)$. 
There is a unique newform $f$ of weight $2$
and level $N=2^\kappa q=14$
which corresponds to the elliptic curve
\[
	Y^2 + X Y + Y = X^3 + 4 X - 6
\]
with Cremona label \texttt{14a1}. 
For each prime $7 \le n <1000$ we computed
$\cZ=\cZ_{\ell_1,\dots,\ell_r}$ 
with $\ell_1,\dotsc,\ell_r$ being
the primes $\ne 7$ satisfying conditions (i), (ii) with
$t \le 200$. The results of this computation
is summarized in Table~\ref{table:cZ}. Note
that by Lemma~\ref{lem:crudeKraus},
$(k \mod{2n}) \in \cZ$.
We deduce that there
are no solutions for prime $n$
satisfying $17 \le n < 1000$, $n \ne 41$.
For $n=7$, $11$, $13$ and $41$,
and for each $m$ in the corresponding $\cZ$,
we compute
$\cW=\cW_{\ell_1,\dotsc,\ell_r}$
as in  Lemma~\ref{lem:refinedKraus}
where $\ell_1,\dotsc,\ell_r$
are now the primes $\ne q$
satisfying (a), (b) and (c)
with $t \le 800$.
We found that $\cW=\emptyset$
in all cases
except for $n=7$, $m=0$,
when $\cW={11-\sqrt{-7}}$.
It follows from Lemma~\ref{lem:refinedKraus}
that 
	$x+7^k \sqrt{-7} = (11-\sqrt{-7}) \cdot \alpha^7$
where $\alpha \in \Z[\theta]$ where $\theta=(1+\sqrt{-7})/2$.
Write $\alpha=(X+Y\theta)$ with $X$, $Y \in \Z$.
Thus
\[
	\frac{x-7^k}{2}+7^k \cdot \theta= 
	(6-\theta)\cdot (X+Y\theta)^7.
\]
Equating the coefficients of $\theta$ on either side
yields the Thue--Mahler equation
\[
-X^7 + 35 X^6 Y + 147 X^5 Y^2 - 105 X^4 Y^3 - 595 X^3 Y^4 - 231 X^2 Y^5 +
161 X Y^6 + 45 Y^7 \; = \; 7^k.
\]
We find that the only solution is $(X,Y,k)=(-1,0,0)$.
Hence $x=-11$, and the corresponding solution to \eqref{eqn:primary}
is $11^2+7=2^7$. We observe that $-11 \equiv 1 \mod{4}$
which is consistent with our assumption $x \equiv 1 \mod{4}$
if $\kappa=1$, made in the statement of Proposition~\ref{prop:primary}.
The other cases are similar. 
\begin{table}
	\begin{tabular}{|c|c|}
	\hline
	$n$ & $\cZ$ \\
	\hline
	\hline
	$7$ & $\{0,8,12\}$\\
	\hline
	$11$ & $\{8\}$\\
	\hline
	$13$ & $\{4\}$\\
	\hline
	$41$ & $\{44\}$\\
	\hline
	other values & $\emptyset$\\
	\hline
\end{tabular}
	\caption{
	For the quadruple $(q,\delta,\kappa,v)=
	(7,0,1,0)$ and for 
	prime $7 \le n <1000$
	we computed
	$\cZ=\cZ_{\ell_1,\dotsc,\ell_r}$
	as given by Lemma~\ref{lem:crudeKraus}.
	Here we chose $\ell_1,\dotsc,\ell_r$
	to be the primes $\ne q$
	satisfying (i) and (ii) with $t \le 200$.}
	\label{table:cZ}
\end{table}

\subsection*{Proofs of Theorems~\ref{thm:odd} and~\ref{thm:modest}}
We now deduce Theorems~\ref{thm:odd} and \ref{thm:modest}
from Proposition~\ref{prop:primary}.
These two theorems concern the equation
$x^2-q^{2k+1}=y^n$ with $n \ge 3$ and $q \nmid x$.
Thus we are in the $\delta=1$ case of the proposition.
By the remarks following the statement of 
the proposition we are reduced to the case $n \ge 7$
is prime. Theorem~\ref{thm:odd}
is concerned with the primes $q$ appearing in \eqref{eqn:goodq},
whilst Theorem~\ref{thm:modest}
deals with $q=41$, $73$, $89$ and $97$.
A glance at Table~\ref{table:ells}
reveals that all the elliptic curves $E$ appearing
in Proposition~\ref{prop:primary} for the case $\delta=1$
in fact correspond to the values $q=41$, $73$, $89$ and $97$.
Theorems~\ref{thm:odd} and~\ref{thm:modest}
now follow immediately from the proposition.

\subsection*{A Remark}
It is well-known that the exponent $n$
can be explicitly bounded in \eqref{eqn:primary}
in terms of the prime $q$. For example, if $\delta=1$
and $\kappa=5$ (i.e. $y$ is odd) then
Bugeaud \cite{Bu-Acta} shows that
% the case where $y$ is odd in  equation (\ref{covid1}), the situation is
%particularly amenable to application of bounds for linear forms in logarithms,
%and we find, via a result of Bugeaud \cite{Bu-Acta}, that 
\begin{equation} \label{Bu-upper}
n \leq 4.5 \times 10^6 q^2 \log^2 q.
\end{equation}
Let $(q,\delta,\kappa,v)=(73,1,5,1)$
and $E$ be the elliptic curve with Cremona label \texttt{2336a1};
this is one of the two outstanding cases 
from Table~\ref{table:ells}
for which the bound \eqref{Bu-upper} is applicable.
We are in fact able to substantially improve this bound for the case in
consideration
through a specialization and minor refinement (we omit the details) of Bugeaud's approach and deduce that
\[
	n < 6 \times 10^6.
\]
Theorem~\ref{thm:modest} 
only resolves $x^2-73^{2k+1}=y^n$ for $3 \le n \le 1000$.
It is natural to ask whether we can apply the same
technique, namely Lemma~\ref{lem:crudeKraus},
to show that there are no solutions for prime exponents
$n$ in the range $1000< n < 6 \times 10^6$.
Write $n_u$ for the smallest prime $>2^u$.
For $10 \le u \le 22$ the prime $n=n_u$
belongs to the range $1000<n<6 \times 10^6$.
For each of these $13$ primes we computed
primes $\ell_1,\dotsc,\ell_r$ satisfying
conditions (i) and (ii) such that $\cZ_{\ell_1,\dotsc,\ell_r}=\emptyset$,
whence by Lemma~\ref{lem:crudeKraus} there are no solutions for that
particular exponent $n$. Table~\ref{table:timings}
records the values of $\ell_1,\dotsc,\ell_r$ as 
well as the time taken to perform the corresponding computation
in \texttt{Magma} on a single processor.
There are $412681$ primes in the range  $1000 < n < 6 \times 10^6$.
On the basis of the timing in the table
we crudely estimate that it would take
around $60$ years to carry out the computation (on a single processor)
for all $412681$ primes.
\begin{table}
\begin{tabular}{|c|c|c|}
\hline
$n$ & $\{\ell_1,\dotsc,\ell_r \}$ & Time\\
\hline\hline
$2^{10}+7=1031$ & $\{  2063, 12373, 30931 \}$ & 0.18 seconds\\
\hline
$2^{11}+5=2053$ &
$\{ 94439, 110863, 143711, 168347, 197089 \}$ & 7.75 seconds\\
\hline
$2^{12}+3=4099$ & $\{ 73783, 98377, 114773 \}$ & 4.390 seconds\\
\hline
$2^{13}+17=8209$ & $\{ 246271, 525377, 574631 \}$ & 15.5 seconds\\
\hline
$2^{14}+27=16411$ & $\{ 98467, 459509, 590797 \}$ & 6.19 seconds\\
\hline
$2^{15}+3=32771$ & $\{ 65543, 983131, 1179757 \}$ & 3.91 seconds\\
\hline
$2^{16}+1=65537$ &
	$\{ 917519, 1310741, 1703963, 2359333 \}$ & 57.51 seconds\\
\hline
$2^{17}+29=131101$ & $\{ 2097617, 9439273, 11799091, 12585697 \}$ & 142.59 seconds\\
\hline
	$2^{18}+3=262147$ & 
	$\{ 1048589, 4194353, 6291529 \}$ & 65.89 seconds \\
\hline
	$2^{19}+21=524309$ &
	$\{ 6291709, 10486181, 23069597 \}$ & 402.12 seconds\\
\hline
	$2^{20}+7=1048583$ & 
	$\{ 20971661, 25165993, 44040487 \}$ & 1319.57 seconds\\
\hline
	$2^{21}+17=2097169$ & 
	$\{ 37749043, 176162197, 188745211 \}$ &
2468.46 seconds\\
\hline
	$2^{22}+15=4194319$ &  
	$\{75497743, 92275019, 100663657\}$ &  4983.07 seconds\\
\hline
%$2^{23}+9=8388617$ & $\{ 201326809, 218104043, 318767447 \}$ & 13740.49 seconds\\
%\hline
\end{tabular}
	\vskip2ex \caption{
Write $n_u$ for the smallest prime $>2^u$.
For $10 \le u \le 22$ the prime $n=n_u$
belongs to the range $1000<n<6 \times 10^6$.
The table lists the primes $n=n_u$
in this range
and, for each, a set of  
primes $\ell_1,\dotsc,\ell_r$ satisfying
conditions (i), (ii) such that $\cZ_{\ell_1,\dotsc,\ell_r}=\emptyset$.
It also records the time the computation took for each of these values of $n$, on a single processor.	
	}\label{table:timings}
\end{table}

We shall shortly give a substantially faster method for treating the case $\delta=0$.
Alas this method is not available for $\delta=1$, as we explain in due course.

%---------------------------------------------------------------------------------
\section{The proof of Theorem \ref{thm:everything} : large exponents} \label{Thm2}
%-----------------------------------------------------------------------------------

We now complete the proof of Theorem \ref{thm:everything} 
which is concerned, for prime $3 \le q <100$, with the equation $x^2+q^\alpha=y^n$, subject  to the assumptions that $q \nmid x$
and $n\ge 3$.
The exponents $n=3$ and $n=4$ were treated in Section~\ref{sec:tiny}, so
we may suppose that $n \ge 5$ is prime. The case $\alpha=2k$
was handled in Section~\ref{Thm4},
so we suppose further that $\alpha=2k+1$. The case with $y$  odd
was the topic of  Section~\ref{sec:oddodd+}, so we may assume that $y$ is even.
Finally,  the case with exponent $n=5$ was resolved in Section~\ref{sec:n=5}, whence we may suppose that $n \ge 7$ is prime. 
To summarize, we are reduced to treating the equation
\begin{equation}\label{eqn:summary}
x^2+q^{2k+1}=y^n, \qquad k \ge 0, \quad q \nmid x, \quad \text{$y$ even,} \; \; n \geq 7 \mbox{ prime.}. 
\end{equation}
By Proposition~\ref{prop:primary},
we may in fact suppose that $n >1000$ and that
\begin{equation}\label{qval}
q \in \{ 7, \, 23,\, 31,\, 47, \, 71,\, 79 \}.
\end{equation}
For convenience, we restate Proposition~\ref{prop:primary} specialized
to our current situation.
\begin{lem}\label{bounders}
Let $q$ be one of the values in \eqref{qval}.
Let $(x,y,k)$ satisfy \eqref{eqn:summary},
where $n>1000$ is prime. Suppose without loss of generality that $x \equiv 1 \mod{4}$.
Attach to this solution the Frey--Hellegouarch elliptic curve
\[
	G=G_{x,k} \; : \; Y^2=X^3+4x X^2+ 4(x^2+q^{2k+1}) X .
\]
Then $\overline{\rho}_{G,n} \sim \overline{\rho}_{E,n}$ where $E$ is an elliptic curve of conductor $2q$
and non-trivial $2$-torsion given in Table~\ref{table:qF}.
\end{lem}
\begin{table}[h]
\begin{centering}
\begin{tabular}{|c|c|c|}
\hline
$q$ & Cremona Label for $E$ & A minimal model for $E$ \\
\hline\hline
$7$ & \texttt{14a1} & $Y^2 + XY + Y = X^3 + 4X - 6$\\
\hline
$23$ & \texttt{46a1} & $Y^2 + X Y = X^3 - X^2 - 10 X - 12$\\
\hline
$31$ & \texttt{62a1} & $Y^2 + X Y + Y = X^3 - X^2 - X + 1$\\
\hline
$47$ & \texttt{94a1} & $Y^2 + X Y + Y = X^3 - X^2 - 1$\\
\hline
$71$ & \texttt{142c1} & $Y^2 + X Y = X^3 - X^2 - X - 3$\\
\hline
$79$ & \texttt{158e1} & $Y^2 + X Y + Y = X^3 + X^2 + X + 1$\\
\hline
\end{tabular}
\caption{}
\label{table:qF}
\end{centering}
\end{table}

%It remains to handle equation (\ref{Main-eq}) under the conditions (\ref{del0}). 

\subsection*{Upper bounds for $n$ : linear forms in logarithms, complex and $q$-adic}
We will appeal to bounds for linear forms in logarithms to deduce  an upper bound for the prime exponent $n$ in \eqref{eqn:summary}
where $q$ belongs to  \eqref{qval}.  
\begin{prop} \label{prop:upper}
Let $q$ belong to the list \eqref{qval}. Let $(x,y,k)$ satisfy \eqref{eqn:summary}
with prime exponent $n>1000$. Then $n < U_q$ where
\begin{equation}\label{eqn:U_q}
U_q=
\begin{cases}
1.43 \times 10^8 & \mbox{if } q=7, \\
5.51 \times 10^8 & \mbox{if } q=23, \\
2.49 \times 10^8 & \mbox{if } q=31, \\
1.10 \times 10^9 & \mbox{if } q=47, \\
1.18 \times 10^9 & \mbox{if } q=71, \\
1.13 \times 10^9 & \mbox{if } q=79. \\
\end{cases}
\end{equation}
\end{prop}

To obtain this result, our first order of business will be to produce a lower bound upon $y$.

\begin{lem} \label{ybound}
%If $q$ is in (\ref{qval}) and $x, y, \alpha, n$ is a solution to equation (\ref{Main-eq}) with $\delta=0$ and $n \geq 13$ prime, then 
If there exists a solution to equation (\ref{eqn:summary}), then $y > 4n -4 \sqrt{2n} + 2$.
\end{lem}
\begin{proof}
We suppose without loss of generality that $x \equiv 1 \mod{4}$, so that
we can apply Lemma~\ref{bounders}.
We first show that $y$ is divisible by an odd prime. Suppose otherwise and write $y=2^\mu$ with $\mu \ge 1$.
Then the Frey--Hellegouarch curve $G_{x,k}$ has conductor $2q$ and minimal discriminant $-2^{2n\mu-12} q^{2k+1}$.
A short search of Cremona's tables \cite{Cre} reveals that there are no such elliptic curves
for the values $q$ in \eqref{qval} (recall that $n>1000$).
Thus, there necessarily exists an odd prime $p \mid y$; since $q \nmid y$, we observe that $q \ne p$. By Lemma~\ref{lem:cong},
$$
a_p(E) \equiv \pm (p+1) \mod{n},
$$
where $E$ is given by Table~\ref{table:qF}. 
As $E$ has non-trivial $2$-torsion, we conclude
that $2n \mid (p+1 \mp a_p(E))$. However,
from the Hasse--Weil bounds
\[
0< p+1 \mp a_p(E) < (\sqrt{p}+1)^2 \le (\sqrt{y/2}+1)^2,
\]
and therefore $2n < (\sqrt{y/2}+1)^2$. 
The desired inequality follows.
\end{proof}

Now let $q$ be any of the values in (\ref{qval}), write $M=\Q(\sqrt{-q})$,
and let $\OO_M$ be its ring of integers. Note that the units
of $\OO_M$ are $\pm 1$.  Fix $\fP$ to be one
of the two prime ideals of $\OO_M$ above $2$. 
%Fix
%a prime exponent $n \ge 5$ and let $(x,y,k)$
%be a solution to (\ref{Main-eq}) with $\delta=0$, $\alpha=2k+1$ and $y$ even.
After possibly replacing $x$ by $-x$ we have
\begin{equation}\label{eqn:pretm}
\frac{x+q^k\sqrt{-q}}{2} \cdot \OO_M = \fP^{n-2} \cdot \fA^n,
\end{equation}
where $\fA$ is an ideal of $\OO_M$ with norm $y/2$.
Hence
\[
\frac{x-q^k \sqrt{-q}}{x+q^k \sqrt{-q}} \; = \; \left(\frac{\fP}{\overline{\fP}}\right)^2 \cdot \left(\frac{\overline{\fP} \cdot \overline{\fA}}{\fP \cdot \fA} \right)^{n}.
\]
For all six values of $q$ under consideration,
the class group is cyclic and generated by
the class $[\fP]$.
Let $h_q$ be the class number of $M$; this value is respectively  $1$, $3$, $3$, $5$, $7$ and $5$
for $q$ in \eqref{qval} (see Table~\ref{table:hq}). 
%If $n=h$ (i.e. $(q,n)=(47,5)$, $(71,7)$
%or $(79,5)$) then $[\fP^{n-2} \cdot \fA^n]=[\fP^{-2}]$
%which is a non-trivial element of the class group,
%contradicting \eqref{eqn:pretm}.
As $n >1000$ is prime,  $\gcd(n,h_q)=1$. 
Since $\OO_M$ has class number $h_q$, it follows that 
$\mathfrak{P}^{h_q}$ is principal, say 
$\mathfrak{P}^{h_q}= (\alpha_q) \cdot \OO_M$.
We fix our choice of $\fP$ so that $\alpha_q$ is given by
Table~\ref{table:hq}. 
\begin{table}
\begin{tabular}{|c|c|c||c|c|c|} 
\hline
$q$ & $h_q$ & $\alpha_q$ & $q$ & $h_q$ & $\alpha_q$ \\ 
\hline\hline
$7$ & $1$ &  $\frac{1+\sqrt{-7}}{2}$ & $47$ & $5$ &  $\frac{9+\sqrt{-47}}{2}$  \\
\hline
$23$ & $3$ &  $ \frac{3+\sqrt{-23}}{2}$ & $71$ & $7$ &  $\frac{21+\sqrt{-71}}{2}$ \\
\hline
$31$ & $3$ &  $\frac{1+\sqrt{-31}}{2}$ & $79$ & $5$ &  $\frac{7+\sqrt{-79}}{2}$ \\
\hline
\end{tabular}
\vskip2ex \caption{Here, $h_q$ denotes the class number
of $M=\Q(\sqrt{-q})$,
and $\alpha_q$ is a generator for the principal ideal $\fP^{h_q}$,
where $\fP$ is one of the two prime ideals of $\OO_M$
above $2$.}\label{table:hq}
\end{table}
Write $\beta_q = \alpha_q/\overline{\alpha_q}$. 
%\margnote{This originally said $\beta_q=\overline{\alpha_q}/\alpha$,
%but I changed it to be consistent with the rest of the choices
%on this page.}
Thus 
\begin{equation} \label{fig-newton}
\left( \frac{x-q^k \sqrt{-q}}{x+q^k \sqrt{-q}} \right)^{h_q} = \beta_q^2 \gamma^n,
\end{equation}
where $\gamma \in M$ is some generator for the 
principal ideal $((\overline{\fP} \cdot \overline{\fA})/(\fP \cdot \fA) )^{h_q}$.

To derive an upper bound upon $n$, we will begin by using (\ref{fig-newton}) to find a ``small'' linear form in logarithms.
Write
$$
\Lambda = \log \left(  \frac{x-q^k \sqrt{-q}}{x+q^k \sqrt{-q}} \right).
$$
\begin{lem} \label{smallform}
If  there exists a solution to equation (\ref{eqn:summary}) with  $y^n > 100 \, q^{2k+1}$, then
$$
\log \left|  \Lambda \right| < 0.75 + (k+1/2) \log q - \frac{n}{2} \log y.
$$
\end{lem}
\begin{proof}
The assumption that $y^n > 100 \, q^{2k+1}$, together with, say, Lemma B.2 of Smart \cite{Smart}, implies that 
$$
\left| \Lambda \right| \leq -10 \log (9/10) \left| \frac{x-q^k \sqrt{-q}}{x+q^k \sqrt{-q}}  -1 \right| = -20 \log (9/10)  \frac{q^k \sqrt{q}}{y^{n/2}},
$$
whence the lemma follows.
\end{proof}

To show that $\log |\Lambda|$ here is indeed small, we first require an upper bound upon $k$. From (\ref{fig-newton}), we have that 
$$
\left( \frac{x-q^k \sqrt{-q}}{x+q^k \sqrt{-q}} \right)^{h_q}-1 = \beta_q^2 \gamma^n-1
$$
and so
\begin{equation} \label{ake}
\frac{-2q^k \sqrt{-q}}{x+q^k \sqrt{-q}}  \sum_{i=0}^{h_q-1} \left( \frac{x-q^k \sqrt{-q}}{x+q^k \sqrt{-q}} \right)^i = \beta_q^2 \gamma^n-1.
\end{equation}

Since $\gcd (x,q)=1$, it follows from (\ref{ake}) that, if we set
$$
\Lambda_1 =\gamma^n - \overline{\beta_q}^2,
$$
then $\nu_q(\Lambda_1) \geq k$. To complement this with an upper bound for linear forms in $q$-adic logarithms, we will appeal to Theorem \ref{qlog},
with
$$
q \in \{ 7, 23, 31, 47, 79 \}, \; \; f=1, \;  \; D=2, \; \alpha_1= \gamma, \;  \; \alpha_2 = \overline{\beta_q}, \;
b_1=n, \; b_2=2,
$$
$$
\log A_1= \frac{h_q}{2} \log y, \; \; \log A_2 = \, \frac{1}{2} \log q \; \; \mbox{ and } \; \; b'= \frac{n}{ \log q} + \frac{2}{h_q \log y}.
$$
Here, we use Lemma 13.2 of Bugeaud, Mignotte and Siksek \cite{BMS2} which implies that
$$
h(\alpha_1)= \frac{h_q}{2} \log (y) \; \mbox{ and } \; h(\alpha_2) = \frac{h_q}{2} \log 2.
$$
In case $q=71$, we make identical choices except to take $\log A_2 = \, \frac{7}{2} \log 2$, whence
$$
b'= \frac{n}{ 7 \log 2} + \frac{2}{7 \log y}.
$$
Theorem \ref{qlog} thus yields the inequality
$$
\nu_q(\Lambda_1) \leq  \frac{96 \, q \,  h_q}{\log^3 q} \cdot  \left( \max \left\{ \log \left( b' \right) + \log \log
q + 0.4, 5 \log q \right\} \right)^2 \log y,
$$
for $q \in \{ 7, 23, 31, 47, 79 \}$, and 
%\margnote{again changed $73$ to $79$} 
$$
\nu_{71}(\Lambda_1) \leq  701.2 \cdot  \left( \max \left\{ \log \left( b' \right) + \log \log
71 + 0.4, 5 \log 71 \right\} \right)^2 \log y,
$$
if $q=71$.

Let us now suppose that
\begin{equation} \label{starter!}
n > 10^8,
%n > 
%\begin{cases} 10^{8} & \mbox{if } q \in \{ 7, 23, 31 \},  \\
%10^9 & \mbox{if } q =47,  \\
%10^{10} & \mbox{if } q \in \{ 71, 79 \}. \\
%\end{cases}
\end{equation}
which will certainly be the case if $n \geq U_q$, for $U_q$ as defined in (\ref{eqn:U_q}).
Then, from Lemma \ref{ybound},  in all cases we have that 
$$
b' < 1.001 \frac{n}{\log q}
$$
and hence obtain the inequalities
\begin {equation} \label{corona}
k < \frac{96 \, q \, h_q}{\log^3 q} \cdot  \left( \max \{ \log \left( n \right) + 0.4001, 5 \log q \} \right)^2 \log y, \; \; \mbox{ if } \; \; q \in \{ 7, 23, 31, 47, 79 \}
\end{equation}
and
\begin {equation} \label{corona2}
k < 701.2 \cdot  \left( \max \{ \log \left( n \right) + 0.4001, 5 \log 71 \} \right)^2 \log y, \; \; \mbox{ if } \; \; q =71.
\end{equation}

Now consider
\begin{equation} \label{vch}
\Lambda_2 =h_q \log \left(  \frac{x-q^k \sqrt{-q}}{x+q^k \sqrt{-q}} \right) = n \log \left( \epsilon_1 \gamma \right) + 2 \log \left( \epsilon_2 \beta_q \right) + j \pi i,
\end{equation}
where we take the principal branches of the logarithms and the integers $\epsilon_i \in \{ -1, 1 \}$ and $j$ are chosen so that 
$\mbox{Im}(\log \left( \epsilon_1 \gamma \right))$ and $\mbox{Im}(\log \left( \epsilon_2 \beta_q\right))$ have opposite signs, and we have both
$$
 \left| \log \left( \epsilon_2 \beta_q\right) \right|  < \frac{\pi}{2}
$$
and $|\Lambda_2|$ minimal. Explicitly,
$$
\begin{array}{|c|c|c|c|r|c|} \hline
q & \epsilon_2 & |\log \left( \epsilon_2 \beta_q \right)| & q & \epsilon_2 & |\log \left( \epsilon_2 \beta_q \right)| \\ \hline
7 & -1 &  \arccos(3/4) & 47 & 1 &  \arccos(17/64)  \\
23 & -1 &  \arccos(7/16) & 71 & 1 &  \arccos(185/256) \\
31 & -1 &  \arccos(15/16)  & 79 & -1 &  \arccos(15/64) \\\hline
\end{array}
$$

Assume first that
$$
y^n \leq  100 \, q^{2k+1}.
$$
If $q \in \{ 7, 23, 31, 47, 79 \}$, it follows from (\ref{corona}) that
$$
n  < \frac{2 \log 10}{\log y} + \frac{\log q}{\log y}  +  \frac{192 \, q \, h_q}{\log^2 q} \cdot  \left( \max \{  \log \left( n \right) + 0.4001, 5 \log q \}  \right)^2,
$$
in each case contradicting Lemma \ref{ybound} and (\ref{starter!}). We obtain a like contradiction in case $q=71$ upon considering (\ref{corona2}).

It follows, then that we may assume $y^n > 100 \, q^{2k+1}$ and hence conclude, from Lemma \ref{smallform},  that
$$
\log \left|  \Lambda_2 \right| < \log h_q + 0.75 + (k+1/2) \log q - \frac{n}{2} \log y.
$$
If $q \in \{ 7, 23, 31, 47, 79 \}$,  (\ref{corona}) thus implies that
$$
\log \left|  \Lambda_2 \right| < \log h_q + 0.75 + \frac{1}{2}  \log q  + \frac{96 \, q \, h_q}{\log^2 q} \cdot  \left( \max \{  \log \left( n \right) + 0.4001, 5 \log q \} \right)^2 \log y - \frac{n}{2} \log y.
$$
An analogous inequality holds for $q=71$, upon appealing to (\ref{corona2}).
From Lemma \ref{ybound} and (\ref{starter!}), we find that
\begin{equation} \label{upsies}
\log \left|  \Lambda_2 \right| < -\kappa_q  \,  n \,  \log y,
\end{equation}
where
\begin{equation} \label{kappaq}
\kappa_q =
\begin{cases} 0.499 & \mbox{if } q =7,  \\
0.497 & \mbox{if } q  \in \{ 23, 31 \},  \\
0.494 & \mbox{if } q  =47,  \\
0.486 & \mbox{if } q =71,  \\
0.490 & \mbox{if } q = 79. \\
\end{cases}
\end{equation}
It therefore follows  from the definition of $\Lambda_2$ that 
$$
|j| \pi < \pi  n + 2 \arccos(15/64)+y^{-0.486 n} < \pi  n+ \pi,
$$
and so
\begin{equation} \label{j-bound}
|j| \leq n.
\end{equation}

\subsubsection*{Linear forms in three logarithms}
To deduce an initial lower bound upon the linear form in logarithms $\left| \Lambda_2 \right| $, we will use 
the following,  the main result (Theorem 2.1) of Matveev \cite{Mat}.
\begin{thm}[Matveev] \label{Matveev} 
Let $\mathbb{K}$ be an algebraic number field of degree $D$ over $\mathbb{Q}$ and put $\chi=1$ if $\mathbb{K}$ is real, $\chi=2$ otherwise. Suppose that $\alpha_1, \alpha_2, \ldots, \alpha_{n_0} \in \mathbb{K}^*$ with absolute logarithmic heights $h(\alpha_i)$ for $1 \leq i \leq n_0$, and suppose that
$$
A_i \geq \max \{ D \, h (\alpha_i), \left| \log \alpha_i \right| \}, \; 1 \leq i \leq n_0,
$$
for some fixed choice of the logarithm. Define
$$
\Lambda = b_1 \log \alpha_1 + \cdots + b_{n_0} \log \alpha_{n_0},
$$
where the $b_i$ are integers and  set
$$
B = \max \{ 1, \max \{ |b_i| A_i/A_{n_0} \; : \; 1 \leq i \leq n_0 \} \}.
$$
Define, with $e := \exp(1)$, further, 
$$
\Omega =A_1 \cdots A_{n_0}, 
$$
$$
C(n_0) = C(n_0,\chi) = \frac{16}{n_0! \chi} e^{n_0} (2n_0+1+2 \chi) (n_0+2)(4n_0+4)^{n_0+1} \left( en_0/2 \right)^{\chi},
$$
$$
C_0 = \log \left( e^{4.4 n_0+7} n_0^{5.5} D^2 \log ( e D) \right) \; \mbox{ and } \; W_0 = \log \left(
1.5 e B D \log (eD) \right).
$$
Then, if $\log \alpha_1, \ldots, \log \alpha_{n_0}$ are linearly independent over $\mathbb{Z}$ and $b_{n_0} \neq 0$, we have
$$
\log \left| \Lambda \right| > - C(n_0) \, C_0 \, W_0 \, D^2 \, \Omega.
$$
\end{thm}

We apply Theorem \ref{Matveev} to $\Lambda=\Lambda_2$ with
$$
D=2, \; \chi = 2, \; n_0=3, \; b_3=n, \; \alpha_3=\epsilon_1 \gamma, \; b_2 = -2, \; \alpha_2 = \epsilon_2 \beta_q, \; b_1=j \;  \mbox{ and } \; \alpha_1=-1.
$$
We may thus take 
$$
A_3 =  \log y, \;  A_2 = \max \{ h_q \log 2, |\log \left( \epsilon_2 \beta_q \right)| \}, \; A_1 = \pi \; \mbox{ and } \;  B =n.
$$
Since
$$
4 \, C(3) \, C_0 = 2^{18} \cdot 3 \cdot 5 \cdot 11 \cdot e^5  \cdot \log \left( e^{20.2}  \cdot 3^{5.5} \cdot 4 \log (2e) \right) < 
1.80741 \times 10^{11},
$$
and 
$$
W_0 =  \log \left(  3 e n \log (2e) \right) < 2.63+\log n,
$$
we may  therefore conclude that
$$
\log  \left| \Lambda_2 \right| > - 5.68 \times 10^{11}  \max \{ h_q \log 2, |\log \left( \epsilon_2 \beta_q \right)| \} \left(  2.63+\log n \right) \log y.
$$
It thus follows from (\ref{upsies}) that 
$$
n < \kappa_q^{-1} \, 5.68 \times 10^{11} \max \{ h_q \log 2, |\log \left( \epsilon_2 \beta_q \right)| \} \left(  2.63+\log n \right)
$$
and hence
\begin{equation} \label{first-bound}
n < 
\begin{cases}
2.77 \times 10^{13} & \mbox{if } q=7, \\
8.24 \times 10^{13} & \mbox{if } q \in \{ 23, 31 \}, \\
1.42 \times 10^{14} & \mbox{if } q \in \{ 47, 79 \} \\
2.02 \times 10^{14} & \mbox{if } q=71. \\
\end{cases}
\end{equation}

To improve these inequalities, we appeal to a sharper, rather complicated  lower bound for linear forms in three complex logarithms, due to Mignotte  (Theorem 2 of \cite{Mig2}).

\begin{thm}[Mignotte] \label{miggy}
Consider three non-zero  algebraic numbers $\alpha_1$, $\alpha_2$
and $\alpha_3$, which are either all real and ${}>1,$ or all complex of modulus
one and all ${}\not=1$. Further, assume that the three numbers $\alpha_1, \alpha_2$ and $\alpha_3$ are either all multiplicatively independent, or that two of the numbers are multiplicatively independent and the third is a root of unity.
We also consider three positive
rational integers $b_1$, $b_2$, $b_3$ with $\gcd(b_1,b_2,b_3)=1$, and the linear form
$$
   \Lambda = b_2\log \alpha_2-b_1 \log \alpha_1-b_3\log \alpha_3 ,
$$
where the logarithms of the $\alpha_i$ are arbitrary determinations of the logarithm,
but which are all real or all purely imaginary.
We assume that
$$
0 < |\Lambda| < 2\pi/w,
$$
where $w$ is the maximal order of a root of unity in $\mathbb{Q}(\alpha_1,\alpha_2,\alpha_3)$. Suppose further that
\begin{equation} \label{queer}
   b_2 |\log \alpha_2| =
 b_1\,\vert \log \alpha_1 \vert+  b_3 \,\vert\log \alpha_3\vert \pm  \vert\Lambda\vert 
\end{equation}
and put
$$
d_1 = \gcd(b_1,b_2), \; \;   d_3 = \gcd(b_3,b_2) \; \mbox{ and } \; b_2=d_1 b_2^\prime = d_3 b_2^{\prime\prime} 
$$
Let
$K, L, R, R_1, R_2, R_3, S, S_1, S_2, S_3, T, T_1, T_2, T_3$
be positive rational integers with
$$
K \geq 3, \; L \geq 5, \; R > R_1+R_2+R_3, \; S > S_1+S_2+S_3 \; \mbox{ and } \; T > T_1+T_2+T_3
$$

Let $\rho\ge 2$  be a real number. Let $a_1, a_2$ and $a_3$ be real numbers such that
$$
   a_i \ge  \rho \vert \log \alpha_i \vert
   - \log  \vert \alpha_i\vert +2 D \,{\rm h}\kern .5pt(\alpha_i), \qquad
   i \in \{1, 2, 3 \},
$$
where
$\,D=[\mathbb{Q}(\alpha_1,\alpha_2,\alpha_3) : \mathbb{Q}]\bigm/[\mathbb{R}(\alpha_1,\alpha_2,\alpha_3) : \mathbb{R}]$, and set
$$
U = \left( \frac{KL}{2}+\frac{L}{4}-1- \frac{2K}{3L} \right) \log \rho.
$$
Assume further that
\begin{equation} \label{needed}
U \geq (D+1) \log (K^2L) + gL(a_1R+a_2S+a_3T) + D(K-1) \log b - 2 \log (e/2),
\end{equation}
where 
$$
g=\frac{1}{4}-\frac{K^2L}{12RST} \; \mbox{ and } \;  b=(b_2^\prime \eta_0)(b_2^{\prime\prime} \zeta_0) \left( \prod_{k=1}^{K-1} k! \right)^{-\frac{4}{K(K-1)}},
$$ 
with
$$
\eta_0=\frac{R-1}{2}+\frac{(S-1)b_1}{2b_2} \; \mbox{ and } \; \zeta_0=\frac{T-1}{2} + \frac{(S-1)b_3}{2b_2}.
$$
Put 
$$
\mathcal{V} = \sqrt{(R_1+1)(S_1+1)(T_1+1)}.
$$
If, for some positive real number $\chi$, we have
\begin{enumerate}[label=(\roman*)]
\item $(R_1+1)(S_1+1)(T_1+1) > K \mathcal{M}$, \\
\item $\mbox{Card} \{ \alpha_1^r \alpha_2^s\alpha_3^t \; : \; 0 \leq r \leq R_1, \; 0 \leq s \leq S_1, \; 0 \leq t \leq T_1 \} > L$, \\
\item $(R_2+1)(S_2+1)(T_2+1) > 2 K^2$, \\
\item $\mbox{Card} \{ \alpha_1^r \alpha_2^s\alpha_3^t \; : \; 0 \leq r \leq R_2, \; 0 \leq s \leq S_2, \; 0 \leq t \leq T_2 \} > 2K L$,  and \\
\item $(R_3+1)(S_3+1)(T_3+1) > 6 K^2 L$, \\
\end{enumerate}
where
$$
\mathcal{M} =  \max\Bigl\{R_1+S_1+1,\,S_1+T_1+1,\,R_1+T_1+1,\,\chi \; \mathcal{V} \Bigr\}, \; \; 
$$
then either
\begin{equation} \label{smosh}
\left| \Lambda \right| \cdot \frac{LSe^{LS|\Lambda|/(2b_2)}}{2|b_2|} > \rho^{-KL},
\end{equation}
or at least one of the following conditions ({\bf{C1}}), ({\bf{C2}}), ({\bf{C3}}) holds :

\vskip2ex
\noindent ({\bf{C1}})  \hskip4ex $|b_1| \leq R_1$  and $|b_2| \leq S_1$ and  $|b_3| \leq T_1$,

\vskip1.5ex
\noindent ({\bf{C2}})  \hskip4ex $|b_1| \leq R_2$  and $|b_2| \leq S_2$ and  $|b_3| \leq T_2$,

\vskip1.5ex
\noindent \mbox{({\bf{C3}})}   {\bf either} there exist non-zero rational integers $r_0$ and $s_0$ such that
\begin{equation} \label{rups2}
   r_0b_2=s_0b_1
\end{equation}
with
\begin{equation} \label{rups3}
   |r_0|
   \le \frac{(R_1+1)(T_1+1)}{\mathcal{M}-T_1}
    \;  \mbox{ and } \; 
   |s_0| 
   \le \frac{(S_1+1)(T_1+1)}{\mathcal{M}-T_1},
\end{equation}
{\bf or}
there exist rational integers  $r_1$, $s_1$, $t_1$ and $t_2$, with
$r_1s_1\not=0$, such that
\begin{equation} \label{rups4}
   (t_1b_1+r_1b_3)s_1=r_1b_2t_2, \qquad \gcd(r_1, t_1)=\gcd(s_1,t_2 )=1,
\end{equation}
which also satisfy
$$
    |r_1s_1|
   \le \gcd(r_1,s_1) \cdot
   \frac{(R_1+1)(S_1+1)}{\mathcal{M}-\max \{ R_1, S_1 \}},
$$
$$
   |s_1t_1| \le \gcd(r_1,s_1) \cdot
   \frac{(S_1+1)(T_1+1)}{\mathcal{M}-\max \{ S_1, T_1 \}} 
  $$
and
$$
   |r_1t_2| % &
   \le \gcd(r_1,s_1) \cdot
  \frac{(R_1+1)(T_1+1)}{\mathcal{M}- \max \{ R_1, T_1 \}}.
$$
Moreover, when $t_1=0$ we can take $r_1=1$, and
when $t_2=0$ we can take $s_1=1$.
\end{thm}

We will apply this result to $\Lambda=\Lambda_2$. 
For simplicity, we will provide full details for the case $q=7$ and indicate later the various parameter choices for the other values of $q$ under consideration.
If $j=0$, then $\Lambda_2$ immediately reduces to a linear form in two logarithms and we may appeal to Theorem \ref{LFL2}, with (in the notation of that result) 
$$
c_2=n, \; \; \beta_2 = \epsilon_1 \gamma, \; \; c_1 = 2, \; \; \beta_1 = \frac{1}{\epsilon_2 \beta_q}, \; \; D=1, 
$$
whence we may choose
$$
\log B_2 = \frac{1}{2} \log y \; \; \mbox{ and } \log B_1 =1.
$$
We thus have, from (\ref{starter!}) and Lemma \ref{ybound},
$$
b^\prime = \frac{4}{\log y} + n < 1.001n.
$$
From Theorem \ref{LFL2} with $(m,C)=(10,32.3)$, it follows, again from (\ref{starter!}),  that
$$
\log \left| \Lambda_2 \right| \geq -64.6 \left( \log n + 0.211 \} \right)^2 \log y.
$$
Combining this with inequality (\ref{upsies})  contradicts (\ref{starter!}).

If $j \neq 0$, then we can choose
\begin{equation} \label{case1}
b_1=2, \; \;  \alpha_1 = \frac{1}{\epsilon_2 \beta_q}, \;  \; b_2=n, \; \;  \alpha_2 = \epsilon_1  \gamma, \;  \;   b_3=-j \;  \; \mbox{ and }  \;  \; \alpha_3=-1,
\end{equation}
whence, from the fact that $\mbox{Im}(\log \left( \epsilon_1 \gamma \right))$ and $\mbox{Im}(\log \left( \epsilon_2 \beta_q\right))$ have opposite signs, (\ref{queer}) is satisfied and we have
$$
d_1=d_3=1 \; \; \mbox{ and } \; \; b_2^\prime = b_2^{\prime\prime} =n.
$$
It follows that
$$
h(\alpha_1)= \frac{1}{2} \log (2), \; h(\alpha_2) = \frac{1}{2} \log (y) \; \mbox{ and } \;  h(\alpha_3) =0,
$$
and hence we can take
$$
a_1=\rho  \arccos(3/4)+\log (2), \; a_2=  \rho \pi  + \log (y) \; \mbox{ and } \; a_3 = \rho \pi.
$$
As noted in \cite{BMS2}, if we suppose that $m \geq 1$ and define
\begin{equation} \label{excellent1}
\begin{array}{c}
K=[m L a_1 a_2 a_3], 
\; R_1 = [c_1a_2a_3], \; S_1=[c_1a_1a_3], \; T_1=[c_1a_1a_2],  \; R_2 = [c_2a_2a_3], \\
\\
S_2=[c_2a_1a_3], \; T_2=[c_2a_1a_2], \;
R_3 = [c_3a_2a_3], \; S_3=[c_3a_1a_3] \; \mbox{ and } \; T_3=[c_3a_1a_2],  \\
\end{array}
\end{equation}
where
\begin{equation} \label{excellent2}
\begin{array}{c}
c_1=\max \{ (\chi m L)^{2/3}, (2mL/a_1)^{1/2} \}, \; c_2=\max \{ 2^{1/3} (m L)^{2/3}, (m/a_1)^{1/2} L \}  \\ \\
\mbox{ and } c_3=(6m^2)^{1/3} L, \\
\end{array}
\end{equation}
then conditions (i)-(v) are automatically satisfied. It remains to verify inequality (\ref{needed}).

Define
$$
R = R_1+R_2+R_3+1, \; S = S_1+S_2+S_3+1 \; \mbox{ and } \; T = T_1+T_2+T_3+1.
$$
We choose
$$
\rho = 6.2, \;  L=153, \; m=24 \; \mbox{ and } \; \chi=0.02,
$$
so that
$$
c_1 = (2mL/a_1)^{1/2}, \; \; c_2=(m/a_1)^{1/2} L,
$$
and we have
$$
K = [K_1+K_2\log (y)],
$$
where  
$$
K_1=7208096.205\ldots   \; \mbox{ and } \; K_2=370065.852\ldots.
$$
We thus have
$$
S_1=3796,  \; S_2=33209 \; \mbox{ and } \; S_3=233126.
$$
Since Lemma \ref{ybound} and (\ref{starter!}) together imply that
\begin{equation} \label{wanda}
\log y > 19.8,
\end{equation} 
we find, after a little work, that $\mathcal{M} = R_1+T_1+1$ and that $g < 0.2411$. 

We have
$$
\eta_0= \frac{1}{2} \left( R_1+R_2+R_3 \right) + \frac{1}{n} \left( S_1+S_2+S_3 \right)  < 26105 \log y + 508457  
$$
and, from $|j| \leq n$, 
$$
\zeta_0=\frac{1}{2} \left( T_1+T_2+T_3 \right) + \frac{|j|}{2n} \left( S_1+S_2+S_3 \right)  <6935 \log y + 270132.
$$
From Lemma 3.4 of \cite{Mig2}, we have the inequality
\begin{equation} \label{mess}
\log \left( \prod_{k=1}^{K-1} k! \right)^{\frac{4}{K(K-1)}} \geq 2 \log K -3+\frac{2 \log \left( 2 \pi K/e^{3/2} \right)}{K-1}-
\frac{2+6 \pi^{-2}+\log K}{3K(K-1)},
\end{equation}
whence, from $K > 10^6$,
$$
\log \left( \prod_{k=1}^{K-1} k! \right)^{\frac{4}{K(K-1)}} > 2 \log K -3.
$$
It follows, appealing to  Lemma \ref{ybound} and (\ref{starter!}), that
$$
b < e^3 n^2 \frac{\left( 26105 \log y + 508457 \right) \left( 6935 \log y + 270132 \right) }{\left(7208095.205 +  370065.852 \log y \right)^2}
< 0.0398 n^2 < 3.06 \times 10^{25},
$$
where the last inequality is a consequence of (\ref{first-bound}).

The right-hand-side of inequality (\ref{needed}) is thus bounded above by 
$$
4 \log (K)  +5.8228 \times 10^8 + 2.9895 \times 10^7  \log (y) + 58.6831 K 
$$
while the left-hand-side satisfies
$$
U >  139.5700 K + 67.96.
$$
If inequality (\ref{needed}) fails to hold, it follows that
$$
80.8869 K < 4 \log (K)  + 5.828 \times 10^8 + 2.9895 \times 10^7  \log (y),
$$
contradicting
$$
K> 7208095.205 +  370065.852 \log y 
$$
and (\ref{wanda}).

Note that we have
$$
\frac{LSe^{LS|\Lambda_2|/(2b_2)}}{2|b_2|} = \frac{20665098 \, e^{20665098 |\Lambda_2|/n}}{n} 
%<  \frac{ 22292930  \mbox{exp} \left(\frac{22292930 }{n y^{0.4993}}}}{n} \right)}
$$
and hence, from (\ref{upsies}), 
$$
\frac{LSe^{LS|\Lambda_2|/(2b_2)}}{2|b_2|} 
<  \frac{ 20665098 \,  \mbox{exp} \left( \frac{20665098}{n y^{0.499n}}\right)}{n} < 0.207,
$$
where the last inequality is a consequence of Lemma \ref{ybound} and (\ref{starter!}).
If we have inequality (\ref{smosh}), it thus follows that
$$
\log \left| \Lambda_2 \right| >  1.5 -280 K.
$$
Once again appealing to (\ref{upsies}), we find that
$$
0.499 \, n \log y < 280 K- 1.5 < 280 \left( 7208097 +  370066 \log y \right)
$$
and so
$$
n < 2.077 \times 10^8 + \frac{4.045 \times 10^9}{\log y},
$$
whence, from  (\ref{wanda}),
\begin{equation} \label{starter2}
n < 4.12 \times 10^8.
\end{equation}

If, on the other hand, inequality (\ref{smosh}) fails to be satisfied, from inequality (\ref{starter!}) and our choices of $S_1$ and $S_2$, necessarily ({\bf{C3}}) holds. We have
$\mathcal{M} = R_1+T_1+1$ and hence if (\ref{rups2}) holds then $n \mid s_0$, where 
$$
|s_0| \leq \frac{(S_1+1)(T_1+1)}{R_1+1} < \frac{3797 \, (3798 + 195 \log (y))}{14293 + 733 \log (y)}.
$$
Since the last quantity is monotonically increasing in $y$, the right-hand-side here is at most $1011$, whence, from (\ref{starter!}), we have that $s_0=0$, a contradiction. We thus have (\ref{rups4}). In particular,
\begin{equation} \label{two-log}
  (2t_1- j r_1)s_1=r_1t_2 n,
\end{equation}
for integers $r_1, s_1, t_1, t_2$ with $ r_1 \mid 2$ and
\begin{equation} \label{tee1}
   |s_1t_1| \le \gcd(r_1,s_1) \cdot
  \frac{(S_1+1)(T_1+1)}{R_1+1}   <  \gcd(r_1,s_1) \cdot 1011.
\end{equation}
Since $s_1$ is coprime to $t_2$ and $n > 10^8$ is prime, it follows from (\ref{two-log}) that either $|r_1|=2$ and $|s_1|=1$, or that $|r_1|=|s_1|=1$.
In the first case, we have that
$$
|t_2| n = |j| \pm |t_1|,
$$
while, in the second,
$$
|t_2| n = |j| \pm 2 |t_1|.
$$
In either case, from (\ref{tee1}), necessarily $|t_1| \leq 1010$. Inequality (\ref{j-bound}) thus implies that $t_2 \in \{ 0, \pm 1 \}$.

If $t_2 = \pm 1$, then, we can write $j = \pm (n-t_3)$, where $t_3 \geq 0$ is an integer with $t_3 \leq 2 |t_1| \leq 2020$.
We thus can rewrite $\Lambda_2$ as a linear form in two logarithms :
\begin{equation} \label{dunnoh}
\Lambda_2 = n \log \alpha - \log \beta,
\end{equation}
where
$$
\log \alpha = \log \left( \epsilon_1 \gamma \right) \pm  \pi i \; \mbox{ and } \; 
\log \beta = -2 \log (-\beta_7) \pm  t_3  \pi i.
$$
If, on the other hand, $t_2=0$, then 
$$
|j| \leq 2 |t_1| \leq 2 \cdot 1010 = 2020
$$
and we can again express $\Lambda_2$ as a linear form in two logarithms  (\ref{dunnoh}), where
$$
\alpha = \epsilon_1\gamma \; \; \mbox{ and } \; \; \log \beta = -2 \log \left( - \beta_7  \right) - j \pi i.
$$
For either pair $(\alpha,\beta)$, we have
$$
h(\alpha) = \frac{1}{2} \log y, \; \; h(\beta) = \frac{1}{2} \log 2  \; \mbox{ and } \; \left| \log \beta \right| \leq 2 \arccos (3/4) + 2020 \pi.
$$

We apply Theorem \ref{LFL2} with  
$$
D=1, \; c_2=n, \; \beta_2 =  \alpha, \; c_1 =1 \; \mbox{ and } \; \log \beta_1 = \log \beta, 
$$
so that we may choose
$$
\log B_1 = \frac{1}{2} \log y \; \; \mbox{ and }  \; \; \log B_2 = 2 \arccos (3/4) + 2020 \pi, \; \; 
$$
whence
$$
b^\prime = \frac{1}{2 \arccos (3/4) + 2020 \pi} + \frac{2n}{\log y}.
$$
Choosing $(m,C)=(10,32.3)$,
we conclude, in conjunction with Lemma \ref{ybound}, (\ref{starter!}) and   (\ref{upsies}), that
$$
0.499 \, n <  16.15 \, \log^2 (n)  \left( 2 \arccos (3/4) + 2020 \pi \right),
$$
contradicting (\ref{starter!}). It follows that inequality (\ref{starter2}) holds in all cases.

To sharpen this, we iterate our argument, in each instance using our new upper bound for $n$ in place of (\ref{first-bound}), choosing $K, R_i, S_i, T_i$ and $c_i$ as in (\ref{excellent1}) and (\ref{excellent2}), with, in each case $\chi=0.02$, and $\rho, L$ and $m$ as follows
$$
\begin{array}{|c|c|c|c|c|} \hline
q & \rho & L & m & \mbox{upper bound upon $n$} \\ \hline
7 & 6.2 & 153 & 24 &  4.12 \times 10^8 \\
7 & 5.7 & 106 & 25 &  1.61 \times 10^8 \\
7 & 5.7 & 103 & 24 & 1.46 \times 10^8 \\
7 & 5.6 & 103 & 25 &  1.43 \times 10^8 \\ \hline
\end{array}
$$
We thus conclude that
\begin{equation} \label{vision}
n < 1.43 \times 10^8=U_7.
\end{equation}

We argue similarly for $q \in \{ 23, 31, 47, 71, 79 \}$, taking the $\alpha_i$ and $\beta_i$ as in (\ref{case1}) (provided $j \neq 0$; the $j=0$ case reduces quickly to linear forms in two logarithms and better bounds), so that we have
$$
h(\alpha_1)= \frac{h_q}{2} \log (2), \; h(\alpha_2) = \frac{h_q}{2} \log (y) \; \mbox{ and } \;  h(\alpha_3) =0,
$$
and may choose
$$
a_1 =  |\log \left( \epsilon_2 \beta_q \right)| \rho + h_q \log 2, \; \; a_2 =\pi \rho  + h_q \log y \; \; 
\mbox{ and } \; \; a_3 = \pi \rho.
$$
For the other parameters, we make
 choices as follows; in all cases, $\chi=0.02$.
 \vskip0.1ex
$$
\begin{array}{|c|c|c|c|c|} \hline
q & \rho & L & m & \mbox{upper bound upon $n$} \\ \hline
23 & 6.6 & 166 & 17  & 1.42  \times 10^9 \\
23 & 6.3 & 107 & 20  &  6.12 \times 10^8 \\
23 & 6.8 & 98 & 18  & 5.51 \times 10^8 = U_{23} \\ \hline
31 & 8.4 & 139 & 15 &  7.18 \times 10^8 \\ 
31 & 8.2 & 89 & 15 &   2.75 \times 10^8 \\ 
31 & 7.8 & 87 & 16 &   2.49 \times 10^8 = U_{31} \\  \hline
47 & 7.4 & 160 & 14  &  2.74 \times 10^9 \\ 
47 & 7.7 & 98 & 15 &  1.19 \times 10^9 \\ 
47 & 7.8 & 96 & 14  &  1.10 \times 10^9 = U_{47} \\ \hline
71 & 9.4 & 139 & 11  &  2.96 \times 10^9 \\ 
71 & 9.4 & 92 & 11  &  1.27 \times 10^9 \\ 
71 & 8.9 & 90 & 12 &  1.18 \times 10^9 = U_{71} \\ \hline
79 & 7.3 & 157 & 15  &  2.80 \times 10^9 \\ 
79 & 7.7 & 98 & 15 &  1.22 \times 10^9 \\ 
79 & 7.8 & 96 & 14 &  1.13 \times 10^9 = U_{79} \\ \hline
\end{array}
$$
\vskip0.5ex
\noindent This completes the proof of Proposition~\ref{prop:upper}.

\subsection*{Proof of Theorem~\ref{thm:everything}}
We now finish the proof of Theorem~\ref{thm:everything}.
By the remarks at the beginning of the current section,
we are reduced to considering solutions $(x,y,k)$ to 
\eqref{eqn:summary}, where $q$ belongs to \eqref{qval}.
Thanks to Propositions~\ref{prop:primary}
and~\ref{prop:upper}, we may suppose
that the prime exponent $n$ belongs to the 
range $1000<n < U_q$ where $U_q$ is given by \eqref{eqn:U_q}.

\begin{lem}\label{lem:finaldescent}
Let $(x,y,k)$ be a solution to \eqref{eqn:summary}
where $q$ belongs to \eqref{qval}
and the exponent $n$ is a prime
belonging to the range
$1000<n<U_q$. 
	Let $M=\Q(\sqrt{-q})$.
Let $h_q$ and $\alpha_q$
be as in Table~\ref{table:hq}, and choose $i$ to be the unique 
integer $0 \le i \le h_q-1$
satisfying $ni \equiv -2 \mod{h_q}$.
Write $n^* = (-ni-2)/h_q$. 
Then, after possible changing the sign of $x$,
\begin{equation}\label{eqn:finaldescent}
	\frac{x+q^k \sqrt{-q}}{2} \; = \; \alpha_q^{n^*} \cdot \gamma^n,
\end{equation}
where $\gamma \in \OO_M$. Moreover,
	$\Norm(\gamma)=2^{i-1} y$.
\end{lem}
\begin{proof}
Recall that $h_q$ is the class number of $M$,
and that $\fP^{h_q}=\alpha_q \OO_M$,
where $\fP$ is one of the two prime ideals
of $\OO_M$  above $2$.
From \eqref{qval}, after possibly replacing $x$ by $-x$,
\[
\left(\frac{x+q^k \sqrt{-q}}{2} \right) \cdot \OO_M \; = \; \fP^{-2} \cdot \fA^n,
\]
where $\fA$ is an ideal of $\OO_M$ of norm $y/2$. Now, for the values
of $q$ we are considering, the class group
	is cyclic and generated by $[\fP]$.
Thus there is some $0 \le i \le h_q-1$
such that $\fP^{i} \fA$ is principal. However,
\[
\left(\frac{x+q^k \sqrt{-q}}{2} \right) \cdot \OO_M \; = \; \fP^{-ni-2} \cdot (\fP^{i} \cdot \fA)^n.
\]
We deduce that $\fP^{-ni-2}$ is principal.
As the class $[\fP]$ generates the class group,
we infer that $i$ is the unique 
integer $0 \le i \le h_q-1$
satisfying $ni \equiv -2 \mod{h_q}$.
Write $n^* = (-ni-2)/h_q$. As $\fP^{h_q}=\alpha_q$,
we have $\fP^{-ni-2}=\alpha_q^{n^*} \cdot \OO_M$.
Hence
\[
\frac{x+q^k \sqrt{-q}}{2} \; = \; \alpha_q^{n^*} \cdot \gamma^n,
\]
where $\gamma \in \OO_M$ is a generator for the principal
ideal $\fP^{i} \fA$. We note that $\Norm(\gamma)=2^{i-1} y$.
\end{proof}

The following lemma, inspired by ideas of Kraus \cite{Kraus33p}, provides a computational
framework for showing that \eqref{eqn:summary} has no solutions
for a particular exponent $n$.
\begin{lem}
Let $q$ belong to the list \eqref{qval} and
let $\beta_q=\overline{\alpha_q}/\alpha_q$.
Let $n$ be a prime
belonging to the range
$1000<n<U_q$. 
Let $E$ be the elliptic curve given in Table~\ref{table:qF}.
Let $\ell \ne q$ be a prime satisfying the following three conditions.
\begin{enumerate}
\item[(I)] $-q$ is a square modulo $\ell$;
\item[(II)] $\ell=tn+1$ for some positive integer $t$;
\item[(III)] $a_\ell(E)^2 \not \equiv 4 \mod{n}$.
\end{enumerate}
Let $\fL$ be one of the two prime ideals of $\OO_M$ above $\ell$, and 
write $\F_{\fL}=\OO_M/\fL \cong \F_\ell$.
Let $\beta \in \F_\fL$ satisfy $\beta \equiv \overline{\alpha_q}/\alpha_q \mod{\fL}$.
Choose $g$ to be a cyclic generator for $\F_\fL^*$, 
set $h=g^n$, and define
\[
	\mathfrak{X}_{\ell,n}=\{ \beta^{n^*} \cdot h^j \; : \; j=0,1,\dotsc,t-1  \} \subset \F_{\fL}.
\]
For $\mathfrak{x} \in \mathfrak{X}_{\ell,n}$ let
\[
E_\mathfrak{x} \; : \; Y^2=X(X+1)(X+\mathfrak{x}).
\]
Finally, define
\[
\mathfrak{Y}_{\ell,n}=\{ \mathfrak{x} \in \mathfrak{X} \; : \; a_\fL(E_{\mathfrak{x}})^2 \equiv a_\ell(E)^2 \mod{n}\}.
\]
If $\mathfrak{Y}_{\ell,n} = \emptyset$, then the equation \eqref{eqn:summary} has no solutions.
\end{lem}
\begin{proof}
Suppose that $(x,y,k)$ is a solution to \eqref{eqn:summary} for our particular pair $(q,n)$.
We change the sign of $x$ if necessary so that \eqref{eqn:finaldescent}
holds and let $x^\prime=\pm x$ so that $x^\prime \equiv 1 \mod{4}$.
By Lemma~\ref{bounders}, we know that
$\overline{\rho}_{G_{x^\prime,k},n} \sim \overline{\rho}_{E,n}$.
Observe that $G_{x^\prime,k}$ is either the same elliptic curve
as $G_{x,k}$ if $x^\prime=x$, or it is a quadratic twist
by $-1$ if $x^\prime=-x$. Hence $a_\ell(G_{x,k}) = \pm a_\ell(G_{x^\prime,k})$
for any odd prime $\ell$ of good reduction for either (and hence both) elliptic curves. 
We let $\ell$ be a prime satisfying conditions  (I), (II) and (III). From (III) and (II),
we note that $a_\ell(E) \not \equiv \pm (\ell+1) \mod{n}$.
It follows from Lemma~\ref{lem:cong} that $\ell \nmid y$,
and that $a_\ell(G_{x^\prime,k}) \equiv a_\ell(E) \mod{n}$.
Thus $a_\ell(G_{x,k})^2 \equiv  a_\ell(E)^2 \mod{n}$.
By Lemma~\ref{lem:finaldescent}, identity \eqref{eqn:finaldescent}
holds where $\Norm(\gamma)=2^{i-1} y$. In particular,
$\fL$ is disjoint from the support of $\gamma$
and $\alpha_q$. It follows from \eqref{eqn:finaldescent}
that
\[
\frac{x-q^k\sqrt{-q}}{x+q^k \sqrt{-q}} \; =\; \left(\frac{\overline{\alpha_q}}{\alpha_q}\right)^{n^*} \cdot \left(\frac{\overline{\gamma}}{\gamma}\right)^n .
\] 
As $g$ is a generator of $\F_{\fL}^*$ which is cyclic of order $\ell-1=tn$, and as $h=g^n$, there is some $0 \le j \le t-1$ such that
$(\overline{\gamma}/\gamma)^n \equiv h^j \mod{\fL}$. 
Hence 
\[
\frac{x-q^k\sqrt{-q}}{x+q^k \sqrt{-q}} \; \equiv\; \mathfrak{x} \mod{\fL},
\]
for some $\mathfrak{x} \in \mathfrak{X}_{\ell,n}$.
The Frey--Hellegouarch curve $G_{x,k}$ defined in Lemma~\ref{bounders} can be
rewritten as
\[
Y^2=X 
\left(X + 2 (x-q^k\sqrt{-q}) \right)
\left(X + 2 (x+q^k\sqrt{-q}) \right)
\]
and hence modulo $\fL$ is a quadratic twist of $E_{\mathfrak{x}}$.
We deduce that $a_\fL(E_{\mathfrak{x}})^2 =a_\ell(G_{x,k})^2 \equiv a_\ell(E)^2 \mod{n}$, whence
 $\mathfrak{x} \in \mathfrak{Y}_{\ell,n}$. This completes the proof.
\end{proof}

\bigskip

To finish the proof of Theorem~\ref{thm:everything}, we wrote
a \texttt{Magma} script which, for each $q$ in \eqref{qval}
and each prime $n$ in the interval $1000<n<U_q$,
found a prime $\ell$ satisfying conditions (I), (II) and (III), with, further,
 $\mathfrak{Y}_{\ell,n}=\emptyset$.
The following table gives the approximate time
taken for this computation, on a single processor.

\bigskip

\hskip28ex \begin{tabular}{|c|c|c|c|}
\hline
$q$ & Time & $q$ & Time \\ \hline
	7 & 41 hours & 47 & 435 hours\\
	23 & 200 hours & 71 & 467 hours\\ 
	31 & 82 hours & 79 & 445 hours\\ \hline
\end{tabular}

\vskip2ex
As one may observe from our proofs, for a given $q$, the upper bound $U_q$ upon $n$ in equation (\ref{eqn:summary}), coming from bounds for linear forms in logarithms, depends strongly upon the class number of $\mathbb{Q}(\sqrt{-q})$. It is this dependence which makes extending Theorem \ref{thm:everything} to larger values of $q$ an expensive proposition, computationally.

%---------------------------------------------
\section{Concluding remarks}
%---------------------------------------------

There are quite a few additional Frey-Hellegouarch curves at our disposal, that might prove helpful in completing the solution of equation (\ref{covid1}), for certain of our problematical values of $q$. A number of these arise from considering (\ref{covid1}) as a special case of
$$
x^2 - q^\delta z^\kappa = y^n,
$$
where, say, $\kappa \in \{ 3, 4, 6 \}$ and $0 \leq \delta < \kappa$. In each case, the dimensions of the spaces of modular forms under consideration grow quickly, complicating matters. 
This is particularly true if $\kappa \in \{ 4, 6 \}$, where our Frey-Hellegouarch curve will {\it a priori} be defined over $\mathbb{Q}(\sqrt{q})$, 
and so the relevant modular forms are Hilbert modular forms
which are more challenging to compute than classical modular forms.
%there are additional difficulties in even asserting a suitable modularity statement for our corresponding Galois representations.

In case $y$ is even in equation  (\ref{covid1}) (whence we are in the situation
where our bounds coming from linear forms in logarithms are weaker), we can
attach a Frey-Hellegouarch $\mathbb{Q}$-curve to a potential solution 
(which at least corresponds to  a classical modular form).
To do this, write
$M=\Q(\sqrt{q})$ and $\OO_M$ for the ring of integers of $M$.  Assuming that
$M$ has class number one (which is the case for, say, the remaining values $q
\in \{ 41,  89,  97 \}$), we have
\[
\frac{x+q^k \sqrt{q}}{2}=\delta^r \gamma^{n-2} \alpha^n
\]
for some $r \in \Z$ and $\alpha \in \OO_M$. Here, $\delta$ is a fundamental unit for $\OO_M$ and 
 $\gamma$ is a suitably chosen generator
for one of the two prime ideals above $2$ in $M$. From this equation,
$$
q^k \sqrt{q}=\delta^r \gamma^{n-2} \alpha^n-
\overline{\delta}^r \overline{\gamma}^{n-2} \overline{\alpha}^n.
$$
Treating this as a ternary equation of signature $(n,n,2)$, we can attach to such a solution a Frey--Hellegouarch $\Q$-curve; see for example
\cite[Section 6]{vanLangen}. We will not pursue this here.

\bibliographystyle{abbrv}
\bibliography{samir}
\end{document}